# COMINUSCULE SUBVARIETIES OF FLAG VARIETIES

BENJAMIN McKAY

ABSTRACT. We show that every flag variety contains a naturally defined homogeneous cominuscule subvariety. From the Dynkin diagram of the flag variety, we compute the Dynkin diagram of that subvariety. We also investigate the tangent bundles of flag varieties.

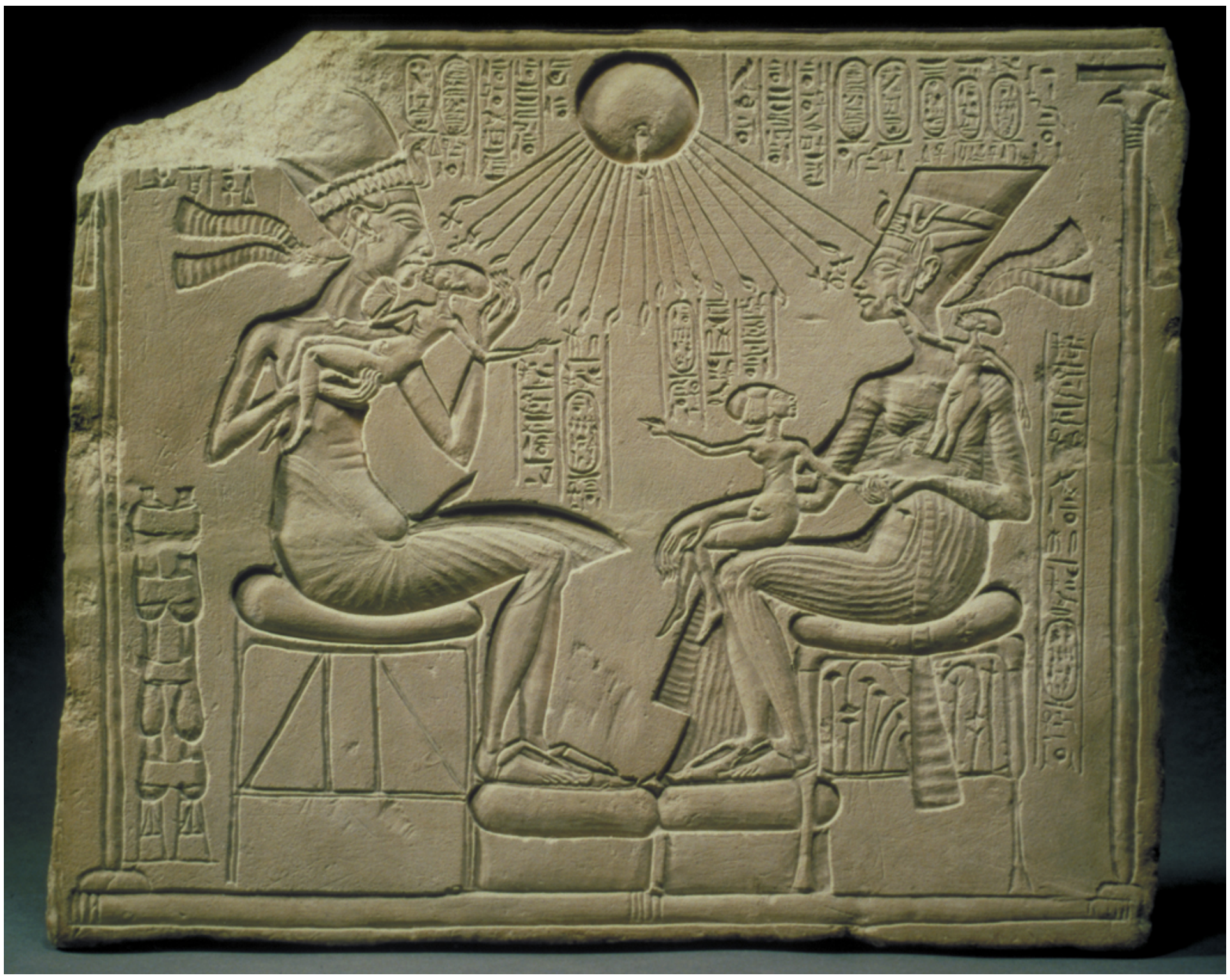

*Date*: November 26, 2026.

This research was supported in part by the International Centre for Theoretical Sciences (ICTS) during a visit to participate in the program Analytic and Algebraic Geometry (Code: ICTS/aag2018/03). I am grateful for the hospitality of the University of Catania, where this paper was mostly written, and in particular to Francesco Russo. Thanks to Indranil Biswas, Anca Mustaţă and Andrei Mustaţă for help with algebraic geometry. This article/publication is based upon work from COST Action CaLISTA, CA21109, supported by COST (European Cooperation in Science and Technology), www.cost.eu.

## Contents

## 1. Introduction

We cannot draw a flag variety, or even its associated root system (except in low dimensions), but we can always draw its Hasse diagram. The Hasse diagrams as drawn by Claus Ringel [50] are very clear, so we will follow his conventions. Our eyes immediately spot in that Hasse diagram its uppermost component. This component is always the Hasse diagram of a unique cominuscule subvariety. This suggests that each flag variety contains an associated homogeneous cominuscule subvariety. We see the root system of the subvariety inside the root system of the flag variety.

**Theorem 1.** *Take any flag variety. Draw its extended Dynkin diagram (defined in §4.5 on page 11). Erase every component not containing a crossed node. Erase every crossed node and its adjacent edges. Erase any component not connected to an affine node. Redraw every affine node as a crossed node. This yields the Dynkin diagram of the associated cominuscule subvariety (defined in §3 on page 8).*

To clarify the theorem, in table 6 on page 56, we draw the Dynkin diagram of every irreducible flag variety and the Dynkin diagram of its associated cominuscule subvariety.

### 1.1. The algorithm in examples.

*Example* 1. The extended Dynkin diagram of $E_8$ is

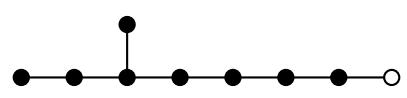

Take the $E_8$-flag variety with Dynkin diagram

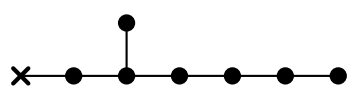

Its extended Dynkin diagram is

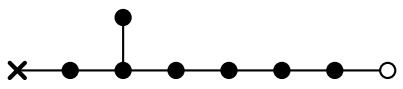

Let us apply the algorithm of our main theorem, theorem 1 on the previous page. Remove the crossed node, giving $D_8$ but with one root still drawn as a hollow $\circ$:

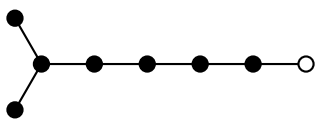

Change the hollow dot $\circ$ to a cross $\times$:

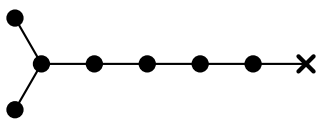

As we will prove, this is the Dynkin diagram of the associated cominuscule subvariety.

*Example* 2. Instead take the $C_7$-flag variety with Dynkin diagram

Its extended Dynkin diagram is

Let us apply the algorithm of our main theorem, theorem 1 on the preceding page. Remove the crossed node and its adjacent edges:

Drop all components not connected to a $\circ$.

Change the hollow dot $\circ$ to a cross $\times$:

As we will prove, this is the Dynkin diagram of the associated cominuscule subvariety.

1.2. **Why the associated cominuscule subvariety matters.** Flag varieties have few regular maps between them [2, 36, 46, 47, 55, 56, 57]. So they have few flag subvarieties and these subvarieties are surprising. Cominuscule varieties are simpler than other flag varieties in many ways; we hope that these cominuscule subvarieties will shed light on their ambient flag varieties.

We give some evidence for their significance in section 10 on page 40. We will see that the associated cominuscule subvariety of an irreducible flag variety $(X, G)$ satisfies an open condition on its tangent spaces, which we call *freedom*. We will see that the associated cominuscule subvariety of any irreducible flag variety is the unique free, connected, closed submanifold of maximal symmetry group, modulo $G$-action. We will see that all free smooth subvarieties are homogeneous. This supports the conjecture that every connected, closed, smooth free subvariety of an irreducible flag variety is a $G$-translate of its associated cominuscule subvariety.

*Example* 3. As in the image of Aten's rays, pick a point $p_0$ and a line $\ell_0$ in the projective plane $\mathbb{P}^2$, with $p_0$ not lying on $\ell_0$.

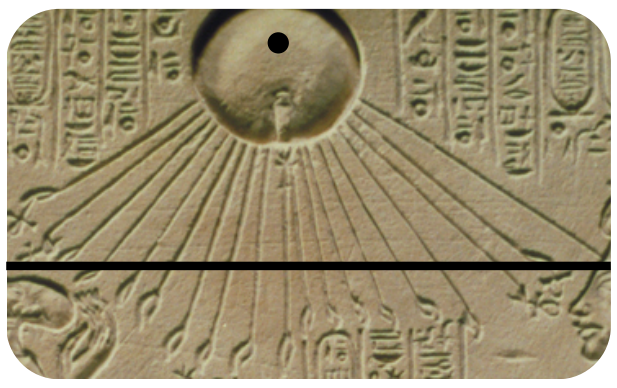

Each point $p$ of $\ell_0$ has an associated pointed line: the pair $(p, pp_0)$.

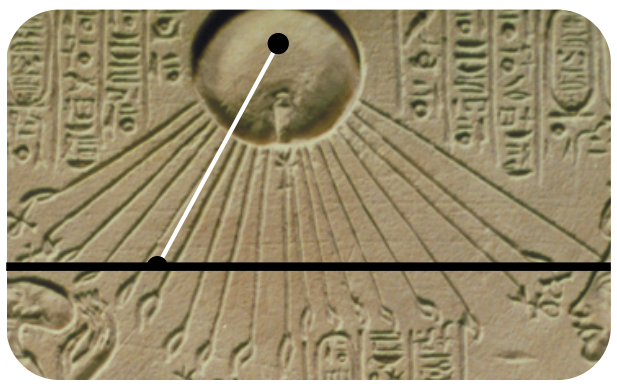

These pointed lines form a rational curve in the variety of pointed lines (*not* in $\mathbb{P}^2$).

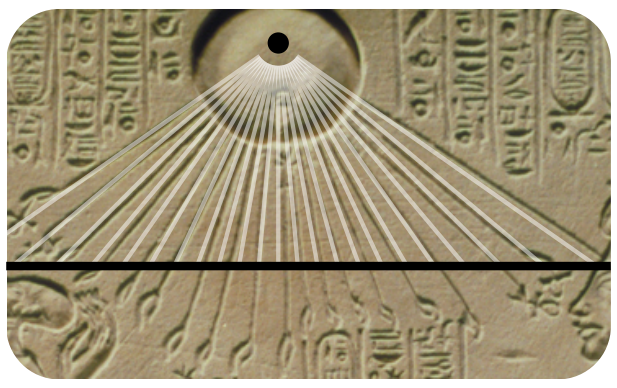

This rational curve is homogeneous under the projective transformations fixing $p_0$ and $\ell_0$; it is the associated cominuscule subvariety of the variety of pointed lines. Each Cartan subgroup of the projective transformations of the plane consists of the transformations which preserve three points in general position.

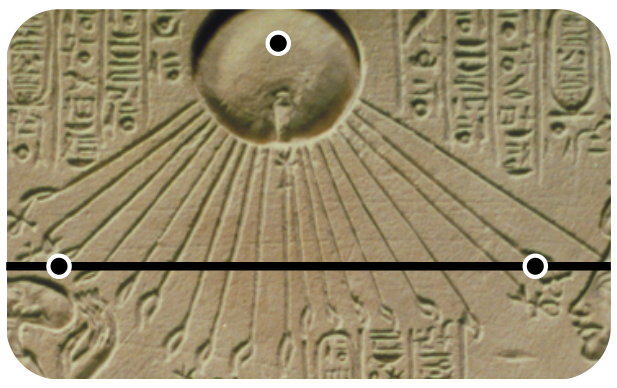

Hence the associated cominuscule subvariety is invariant under a Cartan subgroup. Conversely there are three associated cominuscule subvarieties invariant under any given Cartan subgroup. The projective transformations preserving the point $p_0$ and line $\ell_0$ act transitively on the associated cominuscule subvariety, moving the points $p$ of the line $\ell_0$.

Let us apply the algorithm of theorem 1 on page 3. The Dynkin diagram (defined in §4.2 on page 10) of the variety of pointed lines is ×—×. The extended Dynkin diagram is

.

Cut out the crosses and their edges, leaving only the affine node ∘. Turn the affine node ∘ into a cross ×, giving us ×. This is the Dynkin diagram of the projective line, our rational curve.

*Example* 4. Take a vector space $V$ and write it as the direct sum of linear subspaces $V_i \subseteq V$, say of dimension $n_i$, $i = 1, 2, \dots, k$. Let $G := \mathrm{SL}_V$. Let $P \subset G$ be the subgroup of linear transformations preserving the successive sums

$$V_1, V_1 \oplus V_2, \dots, V_1 \oplus \dots \oplus V_k = V.$$

So $X := G/P$ is the set of partial flags of dimensions

$$0, n_1, n_1 + n_2, \dots, n_1 + \dots + n_k = n.$$

Let $\breve{G} \subset G$ be the subgroup preserving $V_1 \oplus V_k$ and acting as the identity on every $V_i$, $i = 2, \dots, k-1$. Let $\breve{P} \subseteq \breve{G}$ be the subgroup preserving $V_1$. Then every element of $\breve{P}$ preserves $V_1, V_1 \oplus V_2, \dots$. So $\breve{P} \subseteq P$. So $\breve{X} = \breve{G}/\breve{P} \subseteq X = G/P$ is a Grassmannian subvariety of the partial flag variety $X$. The points of $\check{X}$ are precisely the partial flags

$$0 = W_0 \subset W_1 \subset \dots \subset W_k = V$$

obtained by taking any linear subspace $W_1 \subseteq V_1 \oplus V_k$ of the same dimension as $V_1$, and then taking

$$W_2 := W_1 \oplus V_2, W_3 := W_1 \oplus V_2 \oplus V_3, \dots, W_{k-1} := W_1 \oplus V_2 \oplus \dots \oplus V_{k-1}, W_k := V.$$

In other words,

$$(V_1 \oplus V_k) \cap W_1 = W_1, V_2 \cap W_1 = 0, V_3 \cap W_1 = 0, \dots, V_{k-1} \cap W_1 = 0.$$

Moreover

$$V_2 \subseteq W_2, V_2 \oplus V_3 \subseteq W_3, \dots, V_2 \oplus \dots \oplus V_{k-1} \subseteq W_{k-1}.$$

We can phrase this as

$$\dim W_1 = \dim((V_1 \oplus V_k) \cap W_1), 0 = \dim(V_2 \cap W_1) = \dots = \dim(V_{k-1} \cap W_1),$$

and

$$\dim(V_2 \cap W_2) \geqslant n_2, \dots, \dim((V_2 \oplus \dots \oplus V_{k-1}) \cap W_{k-1}) \geqslant n_2 + \dots + n_{k-1}.$$

In particular, $\check{X} \subseteq X$ is an intersection of Schubert varieties, defined relative to different reference flags.

## 2. Review of flag varieties

2.1. **Flag varieties.** A *flag variety* $(X, G)$ is a complex projective variety $X$ acted on transitively and holomorphically by a connected complex semisimple Lie group $G$ [23] pp. 134–135. Some authors prefer the term *generalized flag variety* or *rational homogeneous variety*. We will need to make use of ineffective flag varieties, i.e. $G$ might not act faithfully on $X$. It is traditional to denote the stabilizer $G^{x_0}$ of a point $x_0 \in X$ as $P$; the subgroup $P \subseteq G$ is a connected complex linear algebraic subgroup. A subgroup of $G$ is *parabolic* if it is the stabilizer of a point of a flag variety of $G$.

2.2. **Reducible flag varieties.** The center of any complex semisimple Lie group $G$ lies in every maximal torus, so in every Cartan subgroup [6] p. 220, so in every parabolic subgroup, and hence acts trivially on every flag variety. An *irreducible* flag variety is a flag variety $(X, G)$ with $G$ a simple Lie group, and with only the center of $G$ acting trivially. The *product* of flag varieties $(X_1, G_1), \dots, (X_s, G_s)$ is

$$\begin{aligned} X &= X_1 \times X_2 \times \dots \times X_s, \\ G &= G_1 \times G_2 \times \dots \times G_s. \end{aligned}$$

Every flag variety $(X, G)$, up to finite central extension of $G$, admits a factorization into a product of (i) irreducible flag varieties and (ii) at most one factor $(X_0, G_0)$ consisting of a point $X_0 = \{ x_0 \}$. This factorization is unique up to permutation of the $(X_i, G_i)$ for $i > 0$ and isomorphism. The flag variety $(X, G)$ is effective if and only if all its factors are effective, i.e. if and only if $G_0 = \{ 1 \}$ is trivial and $G_1, \dots, G_s$ are in adjoint form, and then $(X, G)$ is their product (not just up to finite central extension) [1] p. 74.

2.3. **Cominuscule varieties.** A flag variety is *cominuscule* if $\mathfrak{g}/\mathfrak{p} = T_{x_0}X$ is a sum of irreducible complex algebraic $P$-modules. This occurs precisely when there is a compact subgroup $K \subseteq G$ so that $(X, K)$ is a compact Hermitian symmetric space [35] p. 379 Proposition 8.2, [3] p. 26. Some authors prefer the term *compact Hermitian symmetric space*, *cominuscule Grassmannian*, or *generalized Grassmannian* to *cominuscule variety.* Every effective cominuscule variety is a product of the following irreducible effective cominuscule varieties [3] p. 26, Example 3.1.10, [34] theorem 1 p. 401:

| $G$ | $G/P$ | dim | description |
|---|---|---|---|
| $A_r$ | 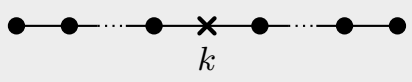 | $k(r+1-k)$ | Grassmannian of $k$-planes in $\mathbb{C}^{r+1}$ |
| $B_r$ |  | $2r-1$ | quadric hypersurface in $\mathbb{P}^{2r}$ |
| $C_r$ |  | $\frac{r(r+1)}{2}$ | space of Lagrangian $r$-planes in $\mathbb{C}^{2r}$ |
| $D_r$ | 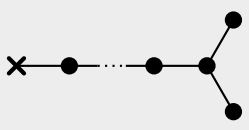 | $2r-2$ | quadric hypersurface in $\mathbb{P}^{2r-1}$ |
| $D_r$ | 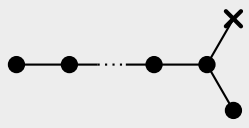 | $\frac{r(r-1)}{2}$ | component of space of null $r$-planes in $\mathbb{C}^{2r}$ |
| $D_r$ | 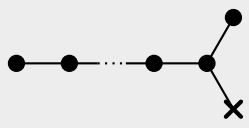 | $\frac{r(r-1)}{2}$ | component of space of null $r$-planes in $\mathbb{C}^{2r}$ |
| $E_6$ | 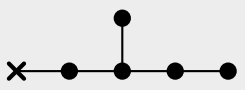 | 16 | complexified octave projective plane |
| $E_7$ | 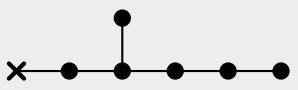 | 27 | space of null octave 3-planes in octave 6-space |

2.4. **Structure of linear algebraic groups.** A complex linear map is *unipotent* if its only eigenvalue is 1. A subgroup of a linear algebraic group is *unipotent* if it consists of unipotent linear maps. Every complex linear algebraic group $G$ has a *unipotent radical*, the unique maximal unipotent normal subgroup, which is a closed complex linear algebraic subgroup [6] p. 85 Theorem 4.5, p. 86 Theorem 4.7, p. 157 11.21.

A complex linear algebraic group is *reductive* if it contains a Zariski dense compact subgroup [61] p. 142 Theorem 2.7, Theorem 2.8. A *reductive Levi factor* of a complex linear algebraic group is a maximal reductive complex linear algebraic subgroup. Every complex linear algebraic group $G$ has a Levi factor, unique up to conjugacy. The group $G$ is a semidirect product of its reductive Levi factor and its unipotent radical [6] p. 158 11.22. Any Cartan subgroup is therefore a subgroup of the reductive Levi factor, after perhaps a conjugacy. Every compact subgroup lies in a maximal compact subgroup, which is unique up to conjugacy, so lies in the reductive Levi factor up to conjugacy [22] p. 531 Theorem 14.1.3. After perhaps extending by some finite group of order a power of 2, the Weyl group embeds in $G$ as a finite subgroup [60]. So this group lies in the reductive Levi factor up to conjugacy.

2.5. **Parabolic subgroups.** A Zariski closed subgroup $P \subseteq G$ of a connected linear algebraic group $G$ is *parabolic* if $X := G/P$ is a projective variety, and this occurs just when $(X, L)$ is a flag variety for a semisimple Levi factor $L \subseteq G$, and this occurs just when $P$ contains a Borel subgroup (i.e. a maximal connected solvable subgroup) [6] p. 148. Every parabolic subgroup is connected [6] p. 197 Proposition 14.18. The unipotent radical of $P$ is denoted $G_+ \subseteq P$, and a reductive Levi factor is denoted $G_0 \subseteq P$, so $P = G_0 \ltimes G_+$. (This notation is potentially confusing. The reader might expect to write these as $P_+$ and $P_0$ since they lie in $P$. But our notation is standard [15] p. 293 theorem 3.2.1. It is due to the presence of the grading of the Lie algebra of $G$ which we will define.) A flag variety is cominuscule just when $G_+$ is abelian [15] p. 296 §3.2.3. Denote the center of the unipotent radical by $Z := Z_{G_+}$.

2.6. **Opposite parabolic subgroups.** Two parabolic subgroups $P, P^{\text{op}} \subseteq G$ of a complex semisimple Lie group are *opposite* if $P \cap P^{\text{op}}$ is a reductive Levi factor of both $P$ and $P^{\text{op}}$. Every parabolic subgroup $P$ has an opposite, unique up to conjugation by an element of $G_+$ [6] p. 198 proposition 14.21.

## 3. Definition of the associated cominuscule subvariety

3.1. **The definition.** Take a flag variety $(X, G)$ and opposite parabolic subgroups $P, P^{\text{op}} \subseteq G$, so that $P$ is the stabilizer of $x_0 \in X$. As above, take their unipotent radicals $G_+, G_+^{\text{op}}$ and the centers $Z, Z^{\text{op}}$ of these. Let $\check{G} := \langle Z, Z^{\text{op}} \rangle \subseteq G$ be the subgroup generated by $Z \cup Z^{\text{op}}$, $\check{P} := \check{G} \cap P$, $\check{X} := \check{G}/\check{P}$. Then $(\check{X}, \check{G})$ is the *associated cominuscule subvariety* passing through the point $x_0 \in X$ determined by the chosen opposite parabolic subgroup.

**Lemma 1.** *The associated cominuscule subvariety of a flag variety $(X, G)$ is uniquely determined up to $G$-conjugacy.*

*Proof.* Since $G$ acts transitively on $X$, it suffices to prove the conjugacy of any two associated cominuscule flag varieties through the same point. Take a point $x_0 \in X$ and consider two associated cominuscule subvarieties passing through $x_0$. Let $P := G^{x_0}$. So it suffices to prove that any two associated cominuscule subvarieties passing through $x_0$ are $P$-conjugate. By definition, each associated cominuscule subvariety with a given parabolic subgroup $P$ is determined by the group generated by $Z, Z^{\text{op}}$. But $Z$ is the center of the unipotent radical of $P$, and thus determined by $P$. So each associated cominuscule subvariety through the given point $x_0$ is determined by the chosen opposite parabolic subgroup $Z^{\text{op}}$, the center of the unipotent radical of the choice of opposite parabolic subgroup $P^{\text{op}}$, and thus determined by $P^{\text{op}}$. But the opposite parabolic subgroup is uniquely determined up to $G_+$-conjugacy [6] p. 198 proposition 14.21, and hence up to $P$-conjugacy. □

The tangent space to $\check{X}$ at $x_0$ is $\mathfrak{z}^{\text{op}}$. Intuitively, the associated cominuscule variety arises by integrating $G$-translates of that subspace. Those $G$-translates are linear subspaces transverse to the largest $P$-invariant holomorphic distribution on the flag variety $(X, G)$. Thus $\check{X}$ is an integral manifold of those $G$-translates; it is not a priori obvious that the $G$-translates have integral manifolds. So our interest in associated cominuscule subvarieties fits into a larger story about tangent geometry of flag varieties [25, 26, 31].

3.2. **Example: the general linear flag variety.** We return to the study of the flag varieties of $A_{n-1} = \mathbb{P}\text{SL}_n$; see example 4 on page 5. We took a vector space $V$ of dimension $n$. We let $G := \mathbb{P}\text{SL}_V$ and $X$ the set of partial flags of dimensions

$$0, n_1, n_1 + n_2, \dots, n_1 + \dots + n_k = n.$$

We wrote $V$ as the direct sum of linear subspaces $V_i \subseteq V$, say of dimension $n_i$, $i = 1, 2, \dots, k$. Supposing that $V = \mathbb{C}^n$, our automorphism becomes transpose inverse. For concreteness take $k = 4$. Let $G' \subseteq G$ be the subgroup preserving $\check{X} \subseteq X$ and set $P' := G' \cap P$. In our table, each line is a pair $\Gamma\ M$, of a group $\Gamma$ and a matrix $M$ with some unspecified entries. In each pair, $\Gamma$ is the group of unimodular matrices of the given form $M$, modulo rescaling by the matrices $\lambda I$ of that form where $\lambda$ is an $n^{\text{th}}$ root of unity.

$$G \quad \begin{pmatrix} * & * & * & * \\ * & * & * & * \\ * & * & * & * \\ * & * & * & * \end{pmatrix}$$

$$P \quad \begin{pmatrix} * & * & * & * \\ 0 & * & * & * \\ 0 & 0 & * & * \\ 0 & 0 & 0 & * \end{pmatrix}$$

$$G_+ \quad \begin{pmatrix} I & * & * & * \\ 0 & I & * & * \\ 0 & 0 & I & * \\ 0 & 0 & 0 & I \end{pmatrix}$$

$$Z \quad \begin{pmatrix} I & 0 & 0 & * \\ 0 & I & 0 & 0 \\ 0 & 0 & I & 0 \\ 0 & 0 & 0 & I \end{pmatrix}$$

$$G_0 \quad \begin{pmatrix} * & 0 & 0 & 0 \\ 0 & * & 0 & 0 \\ 0 & 0 & * & 0 \\ 0 & 0 & 0 & * \end{pmatrix}$$

$$\breve{G} \quad \begin{pmatrix} * & 0 & 0 & * \\ 0 & I & 0 & 0 \\ 0 & 0 & I & 0 \\ * & 0 & 0 & * \end{pmatrix}$$

$$\breve{P} \quad \begin{pmatrix} * & 0 & 0 & * \\ 0 & I & 0 & 0 \\ 0 & 0 & I & 0 \\ 0 & 0 & 0 & * \end{pmatrix}$$

$$G' \quad \begin{pmatrix} * & 0 & 0 & * \\ 0 & * & 0 & 0 \\ 0 & 0 & * & 0 \\ * & 0 & 0 & * \end{pmatrix}$$

$$P' \quad \begin{pmatrix} * & 0 & 0 & * \\ 0 & * & 0 & 0 \\ 0 & 0 & * & 0 \\ 0 & 0 & 0 & * \end{pmatrix}$$

To be more precise, $G_+$ consists of the matrices

$$\begin{pmatrix} \lambda I & * & * & * \\ 0 & \lambda I & * & * \\ 0 & 0 & \lambda I & * \\ 0 & 0 & 0 & \lambda I \end{pmatrix}$$

with determinant 1, up to rescaling by $n^{\text{th}}$ roots of unity. But for each such matrix equivalence class, pick any representative and we can pick that root of unity uniquely to write it as an actual matrix with $\lambda = 1$, and similarly for $Z$. Clearly $\breve{G}$ has Lie algebra containing all root vectors of all $P$-maximal and $P$-minimal roots, and is generated by the 1-parameter subgroups of those root vectors. So it contains the root vectors of the root system generated by these, and the Cartan subgroup of that root system, hence a semisimple group. Note that $P'$ preserves $V_1$, $V_2$, $V_3$, and $V_1 \oplus V_4$, whereas $G'$ preserves $V_2$, $V_3$ and $V_1 \oplus V_4$. Hence $\breve{X} = \breve{G}/\breve{P} = G'/P'$ is the Grassmannian of linear subspaces of dimension $n_1$ inside $V_1 \oplus V_4$.

## 4. Dynkin diagrams

4.1. **Compact roots.** Denote the Lie algebras of $P \subseteq G$ by $\mathfrak{p} \subseteq \mathfrak{g}$. One can select a Cartan subgroup of $G$ lying inside $P$, whose positive root spaces all lie in $\mathfrak{p}$. A simple root $\alpha$ is $P$-*compact* (*compact* if $P$ is understood) if the root space of $-\alpha$ belongs to the Lie algebra of $P$. This occurs just when the $\mathfrak{sl}_2$ generated by the root vectors of $\alpha, -\alpha$ lies inside the Lie algebra $\mathfrak{p}$ of $P$, hence inside a semisimple Levi factor. So we can assume that all $P$-compact roots are roots of a semisimple Levi factor lying inside $G_0$.

4.2. **The Dynkin diagram of a flag variety.** The Dynkin diagram of a flag variety $(X, G)$ is the Dynkin diagram of $G$ decorated with $\bullet$ on each compact simple root and $\times$ on each noncompact simple root [6] p. 197 Proposition 14.18, [23] p. 197 Theorem′. Each flag variety is determined uniquely, up to finite central extension of $G$ and up to isomorphism, by its Dynkin diagram.

4.3. **Affine extensions.** Every Dynkin diagram sits inside its extended Dynkin diagram [8] p. 198, [14] p. 262, [30] lecture 17, [48] p. 164, 5°. The extended Dynkin diagram has one node added to each component of the original Dynkin diagram. We call this node the *affine node*, and draw it as a hollow dot $\circ$, with edges attached as:

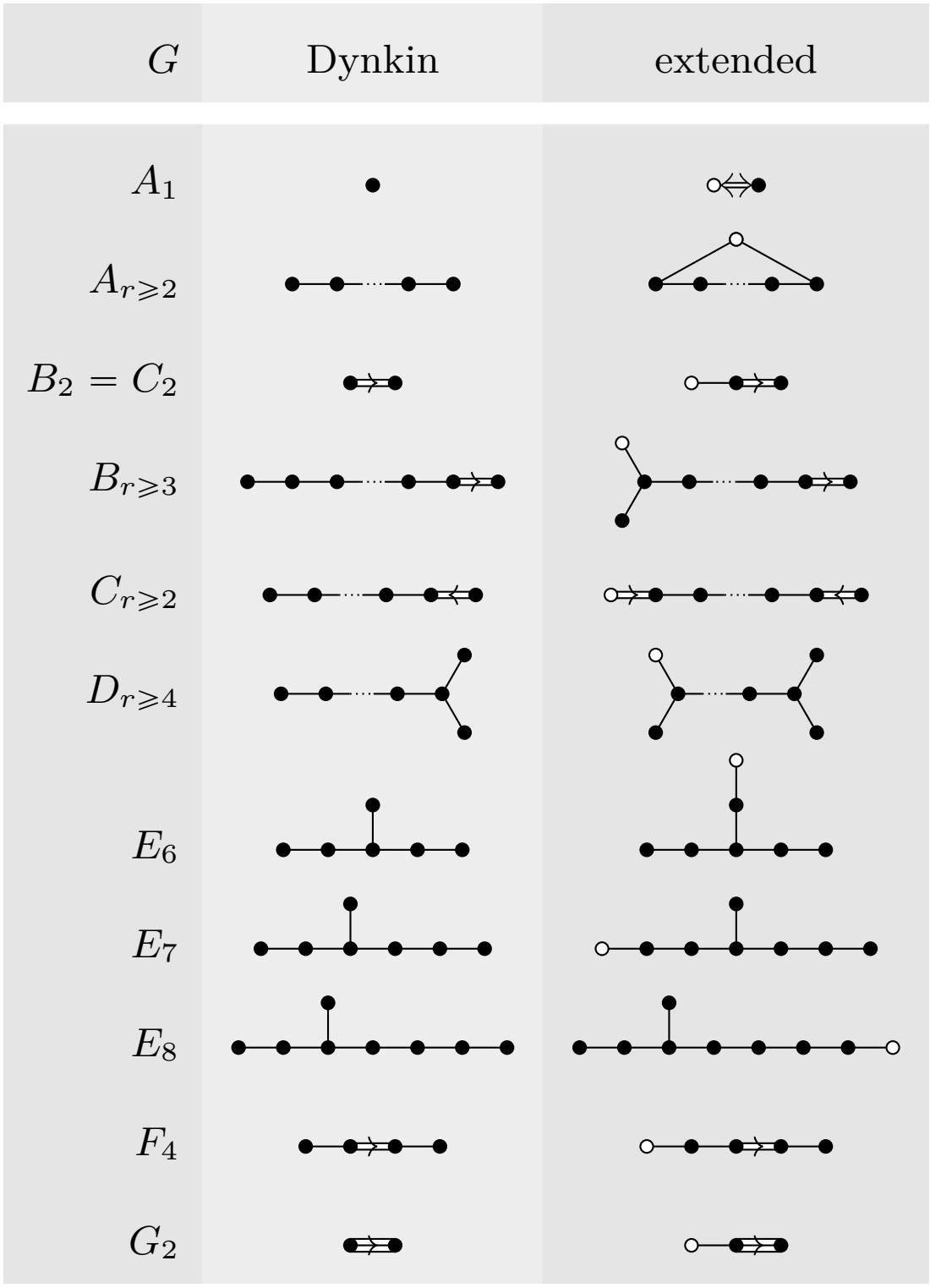

Recall that we add this node to the original Dynkin diagram to represent the negative $\alpha_0 := -\theta$ of the highest root $\theta$. We then compute the Cartan matrix as usual. Consider how we draw edges connecting that new node to each old node according to the usual Dynkin diagram rules. For any finite or extended root system, denote the Cartan matrix entry in row $i$ column $j$ as $n_{ij}$. If $i \neq j$ and $|n_{ij}| \geqslant |n_{ji}|$, we draw $|n_{ij}|$ edges from vertex $j$ to vertex $i$. If $|n_{ij}| > 1$ we draw an arrow pointing to node $i$.

4.4. **Rank 1.** For $A_1$, the Cartan matrix of the extended root system is

$$\begin{pmatrix} 2 & -2 \\ -2 & 2 \end{pmatrix}$$

so both off-diagonal Cartan matrix entries are $-2$. Hence we draw the extended Dynkin diagram of $A_1$ as 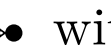 with a double edge with arrows pointing in both directions [29] p. 44 table Aff 1. The total number of edges of the extended $A_r$ Dynkin diagram is $r+1$, for any integer $r \geqslant 1$.

There is a different rule [48] p. 164 for forming Dynkin diagrams: we could instead make the edge between node $i$ and node $j$ have multiplicity $n_{ij}n_{ji}$, and make an arrow from node $j$ to node $i$ just when $|n_{ij}| > |n_{ji}|$. This rule leads to exactly the same Dynkin diagrams and extended Dynkin diagrams, except precisely for the extended Dynkin diagram of $A_1$, yielding a quadruple edge [61] p. 228 table 3, [48] p. 166 .

The choice of convention for drawing the extended Dynkin diagram of $A_1$ does not affect our algorithm for computing the associated cominuscule subvariety: either convention yields the same result.

4.5. **The extended Dynkin diagram of a flag variety.** Take the Dynkin diagram of a flag variety $(X, G)$. The *extended Dynkin diagram* of $(X, G)$ is the associated extended Dynkin diagram, with nodes drawn as × or • as for the flag variety, keeping the affine node drawn as ∘. So the flag variety 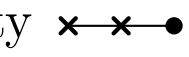 has extended Dynkin diagram

## 5. Reducing to root systems

5.1. **Gradings.** A root system with a basis of simple roots $\alpha_1, \ldots, \alpha_r$ is graded: each root $\sum n_i \alpha_i$ has grade $\sum_i n_i$. For a flag variety $X = G/P$, the root system is also *$P$-graded* by $\sum n_i$, but summing only over the noncompact (crossed) simple roots [15] p. 292. The $P$-grade is also called the *$P$-height* [15] p. 292. A root is *$P$-maximal*, also called a *box root*, if it has maximal and strictly positive $P$-grade in its irreducible factor. The *box* is the set of box roots, following what [13], [37] p. 57 might perhaps call the *maximal box*, by analogy with Young tableaux. We will see in lemma 5 on page 18 that the box generates the root system of $\breve{G}$.

5.2. **The box and the Lie algebra.** The unipotent radical $G_+ \subseteq P$ has Lie algebra $\mathfrak{g}_+ \subseteq \mathfrak{p}$ the direct sum of the root spaces of the positively graded roots [6] p. 197 Proposition 14.18, [32] p. 482. The zero graded roots are invariant under reflection in one another. The sum of the root spaces of the zero graded roots with the span of their associated coroots is the semisimple factor of the reductive Levi factor $\mathfrak{g}_0$ of $\mathfrak{p}$ [6] p. 197 Proposition 14.18. Denote by $\mathfrak{z}$ the Lie algebra of $Z$. For any root $\alpha$, denote by $\mathfrak{g}_\alpha$ the $\alpha$-root space of $\mathfrak{g}$. We will prove:

**Lemma 2.** *The abelian Lie algebra* $\mathfrak{z}$ *is the span of the root vectors of the roots of the box:*

$$\mathfrak{z} = \bigoplus_{\alpha \in \text{box}} \mathfrak{g}_\alpha .$$

*Take the box roots (i.e. roots of the box) as vertices of a graph. If the difference* $\alpha - \beta$ *of two box roots is a* $P$*-compact simple root, corresponding to node* $\ell$ *of the Dynkin diagram of* $(X, G)$*, draw an edge from* $\beta$ *to* $\alpha$ *labelled* $\ell$*. Thus the box becomes a graph. As a graph, the connected components of the box are precisely the boxes of the irreducible factors of the flag variety. Each component contains the highest root of the factor. In particular, the box of an irreducible flag variety is connected and contains the highest root.*

### 5.3. Associated cominuscule subvarieties in rank 2.

*Example* 5. Here we will explain how to read our drawings. We pick a basis of simple roots in each root system. We draw the simple roots as empty circles:

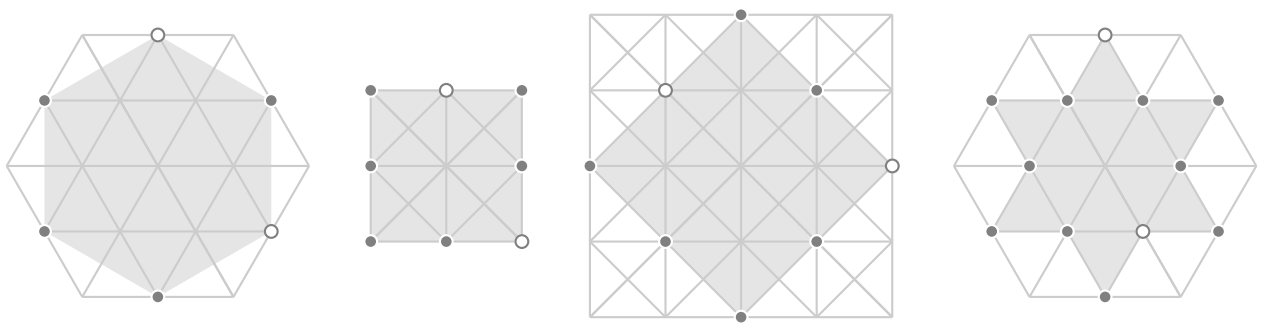

Start with the roots of $G_2$:

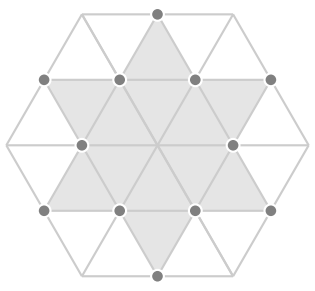

Each parabolic subgroup has Lie algebra consisting of the sum of the Cartan subalgebra and the root spaces of certain roots. To be precise, these roots are the ones which lie on or on one side of a hyperplane. The compact roots are those on the hyperplane. Conversely, draw any hyperplane through the origin, and pick a side of it, and it produces a parabolic subgroup. For example, here is the hyperplane of some parabolic subgroup.

Drawing both the hyperplane and roots together:

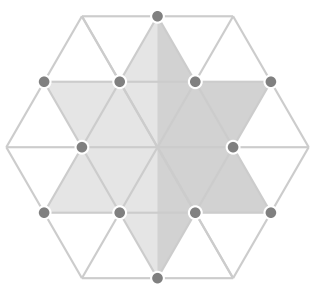

We can always pick a basis of simple roots so that every simple root lies on the hyperplane, or on the chosen side of the hyperplane. It is compact (so uncrossed) just when it is on the hyperplane:

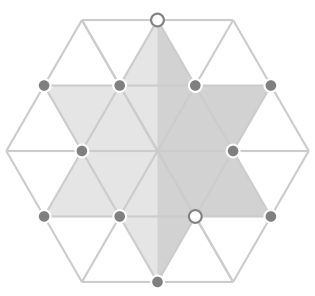

(When working with $G_2$, we will always pick our hyperplane to allow for the chosen bases of simple roots shown above.) We pick the hyperplane to be the zero locus of the real linear function vanishing on uncrossed simple roots and taking value one on crossed simple roots [15] p. 239 Proposition 3.1.2. Grade the roots by $P$-height, i.e. by sum of coefficients of crossed simple roots. Our hyperplane is thus always chosen so that we can see the heights, i.e. roots of a given height lie on a parallel hyperplane:

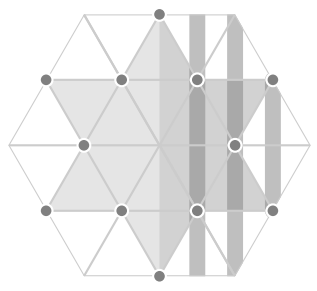

By definition of $\breve{G}$, its root system is the root system generated by the box, i.e. by the root of maximal $P$-grade:

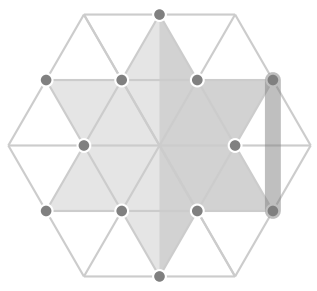

So the $\breve{G}$-roots form the smallest root subsystem containing these. Finally, draw dark lines through the $\breve{G}$-roots:

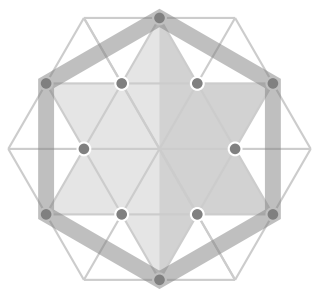

The roots of $\breve{P} = P \cap \breve{G}$ are the roots of $\breve{G}$ lying on or to the indicated side of the hyperplane: 4 of them in our picture:

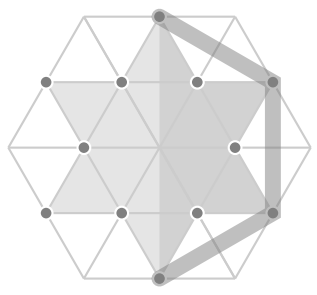

Hence $\breve{X} = \breve{G}/\breve{P}$ has dimension equal to the number of $\breve{G}$-roots not lying in $\breve{P}$, i.e. the dimension of $\breve{X}$ is the number of roots in the box. We can see in the picture that the root system of $\breve{G}$ in this example is that of $A_2 = \mathbb{P}\mathrm{SL}_3$, and that $\breve{X}$ has dimension 2, so must be $(\breve{X}, \breve{G}) = (\mathbb{P}^2, A_2)$. The elements of $G$ preserving $\breve{X}$ form a subgroup of $G$ which, as we will see, maps to $\breve{G}$, which we call the *automorphism group*. We will see that the automorphism group is generated by the flow through $1 \in G$ of various root vectors of $G$, associated to various roots in the root system of $G$. Among these roots are the roots of $\breve{G}$. Below we will colour these in; for our example, the roots of the automorphism group are precisely those of $\breve{G}$.

*Example* 6. We draw the gradings of the positive roots of the parabolic subgroups of the rank 2 simple groups. Under the heading $\breve{\mathfrak{g}}$, we draw the roots of the symmetry Lie algebra of the associated cominuscule subvariety, and under the heading $\mathfrak{g}'$ the roots of the automorphism Lie algebra.

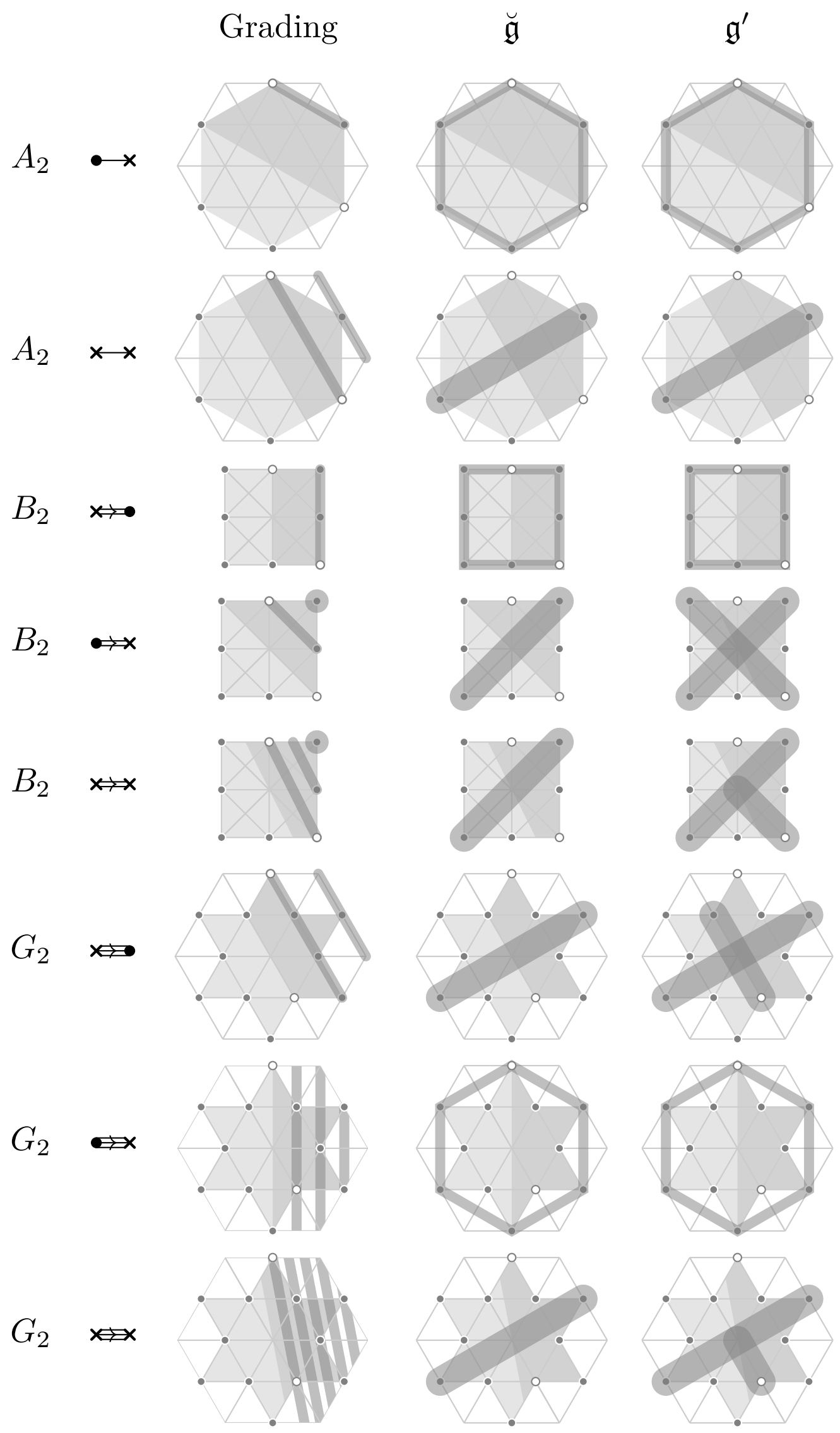

*Example* 7. Any irreducible complete flag variety $X = G/B$ has associated cominuscule subvariety $\times = (\mathbb{P}^1, \mathbb{P}\mathrm{SL}_2)$, reducing maximally. This occurs because $G$ has a unique highest root, whose root space generates $Z$. Note that the subgroup of $G$ acting on the associated cominuscule might be $\mathrm{SL}_2$, but its action factors through the action of $\mathbb{P}\mathrm{SL}_2$. More generally, the subgroup $\breve{G}$ of $G$ might act with a finite kernel on $\breve{X}$. If $X = G/B$ is a complete flag variety and $G$ has $s$ semisimple factors, then $\breve{X} = (\mathbb{P}^1)^s$.

*Example* 8. Similarly, the associated cominuscule subvariety of any adjoint variety is $\times = (\mathbb{P}^1, \mathbb{P}\mathrm{SL}_2)$. The adjoint varieties are precisely those with invariant contact structures. An invariant contact structure arises from a hyperplane in each tangent space, invariant under the parabolic. By invariance, it is a sum of root spaces. So a single root, negative for the parabolic, doesn't have its root vector in this hyperplane, i.e. a single root is at higher weight. So the associated cominuscule subvariety is complementary to the hyperplane, one dimensional, hence a projective line.

## 6. The associated cominuscule subvariety is cominuscule

Henceforth suppose that $G$ is a connected complex semisimple Lie group. Take a factorization

$$X = X_0 \times X_1 \times X_2 \times \cdots \times X_s,$$
$$G = G_0 \times G_1 \times G_2 \times \cdots \times G_s,$$

into irreducible flag varieties $(X_i, G_i)$, $i > 0$, and a point $X_0 = \{ x_0 \}$. Then

$$P = G_0 \times P_1 \times P_2 \times \cdots \times P_s,$$

with $P_i := P \cap G_i$. The unipotent radical of $P$ is obviously the product of the unipotent radicals of the $P_i$:

$$G_+ = \{ 1 \} \times G_{1+} \times \cdots \times G_{s+},$$

where $G_{i+}$ is the unipotent radical of $P_i$. Therefore

$$Z = \{ 1 \} \times Z_1 \times \cdots \times Z_s$$

where $Z_i$ is the center of $G_{i+}$. We prove lemma 2 on page 12.

*Proof.* Let $\mathfrak{z}'$ be the direct sum of the root spaces of the box roots, i.e. the $P$-maximal roots. Recall that $\mathfrak{z}$ is the center of the nilpotent radical $\mathfrak{g}_+$ of $\mathfrak{p}$. We have to prove that $\mathfrak{z}' = \mathfrak{z}$. Note that in our decomposition above,

$$\mathfrak{z}' = \{ 0 \} \oplus \mathfrak{z}'_1 \oplus \cdots \oplus \mathfrak{z}'_s.$$

So it suffices to prove the result for an irreducible flag variety $(X, G)$. By invariance under the Cartan subalgebra, $\mathfrak{z}$ is a sum of root spaces. So we let $B'$ be the box and $B$ the set of roots whose root spaces belong to $\mathfrak{z}$. Proving $\mathfrak{z} = \mathfrak{z}'$ is equivalent to proving that $B = B'$.

Recall that the bracket of root vectors of roots $\alpha, \beta$ is either zero or is a root vector of the root $\alpha + \beta$. If $\alpha \in B'$, i.e. $\alpha$ is $P$-maximal, its associated root vector brackets to zero with the root vector of any root $\beta$ of positive $P$-grade. Hence the root vector lies inside $\mathfrak{z}$, i.e. $\mathfrak{z}' \subseteq \mathfrak{z}$, i.e. $B' \subseteq B$.

Recall that every irreducible root system has a unique highest root, i.e. the highest weight of the adjoint representation [53] p. 61 Theorem 3. Picture a sequence of steps from root to root, by adding only simple roots. Every positive root can be brought to the highest root by such a sequence [19] pp. 330–331 §21.3, [24] p. 56, [33], [53] p. 46 2.12, [53] p. 58 Theorem 3. Adding $P$-compact simple roots preserves $P$-height. Adding $P$-noncompact simple roots increases $P$-height. When this sequence reaches a $P$-maximal root, it can no longer increase the $P$-height. So we can only add $P$-compact roots, moving through $P$-maximal roots, so moving along $B'$. So the highest root is $P$-maximal, so in $B'$, and is obtained from any root of $B'$ by adding only $P$-compact simple roots. Every $B'$-root can thus be brought to the highest root by adding only $P$-compact simple roots. So the box becomes a connected graph, with the highest root as one vertex. The vector space $\mathfrak{z}'$ is therefore an irreducible $G_0$-module [32] p. 332 proposition 5.105, and a $G_0$-submodule of $\mathfrak{z}$.

Pick a root $\alpha \in B$. It suffices to prove that $\alpha \in B'$, i.e. $P$-maximal. Pick a sequence of roots

$$\alpha_0 := \alpha, \alpha_1 := \alpha + \beta_1, \alpha_2 := \alpha + \beta_1 + \beta_2, \ldots, \alpha + \beta_1 + \cdots + \beta_n = \theta,$$

so that each $\beta_j$ is a simple root, travelling to the highest root $\theta$. If $\beta_1$ is compact, invariance under $\mathfrak{g}_0$ ensures that $\mathfrak{g}_{\alpha_1} \subset \mathfrak{z}$ so $\alpha_1 \in B$, and so on. If some $\beta_j$ is noncompact, then $0 \neq [\mathfrak{g}_{\alpha_{j-1}}, \mathfrak{g}_{\beta_j}] \subseteq \mathfrak{g}_{\alpha_j}$, contradicting the fact that $\mathfrak{g}_{\alpha_{j-1}}$ centralizes $\mathfrak{g}_+$ and $\mathfrak{g}_{\beta_j} \subseteq \mathfrak{g}_+$. So all $\beta_j$ are compact roots, so $\alpha$ and $\theta$ have the same $P$-grade.

Similarly, if $\alpha+\beta$ is a root, for some root $\beta$ of zero $P$-height, i.e. a compact root, then $0 \neq [e_\alpha, e_\beta] \in \mathfrak{g}_{\alpha+\beta}$, so $\alpha+\beta$ also belongs to the box, and its root space $\mathfrak{g}_{\alpha+\beta}$ also belongs to $\mathfrak{z}$. So $B$ is invariant under stepping through roots by compact roots. So $B = B'$ and $\mathfrak{z} = \mathfrak{z}'$. □

**Lemma 3.** *The group* $\breve{G}$ *is connected.*

*Proof.* Every parabolic subgroup of a connected reductive Lie group is connected, so $G$ and $P$ are connected [6] p. 154 theorem 11.16. By Langlands decomposition [32] p. 482, $G, P, X, G_+, G_0$ are connected. The subgroup $Z = Z_{G_+}$ is a Zariski closed subgroup of a unipotent linear algebraic group, so connected and isomorphic as an affine variety to complex Euclidean space [16] p. 36 corollary 15.1.11, [18] p. 499, Corollaire 4. The exponential map on the Lie algebra of a unipotent linear algebraic group is a biholomorphism and a regular algebraic morphism mapping subalgebras to subgroups, so taking $\mathfrak{z} \to Z$, an isomorphism of schemes [43] p. 289 Proposition 14.32. By definition of $\breve{G}$, $\breve{G} := \langle Z, Z^{\mathrm{op}} \rangle$ is a path connected complex Lie group, since $Z$ is path connected, and hence $Z^{\mathrm{op}}$ is also path connected. Moreover, $\breve{G}$ is connected in the algebraic Zariski topology by Chevalley's theorem on constructible sets [16] p. 38 Proposition 16.2.1. Hence it is a linear algebraic group [6] p. 57, Prop. 2.2. □

6.1. **Roots and opposite parabolics.** All Borel subgroups of $G$ are conjugate [6] p. 147 chapter IV 11.1, each containing a Cartan subgroup, hence the Lie algebra of $P$ is the sum of (1) the Cartan subalgebra with (2) various root spaces. Every automorphism of a root system arises from an automorphism of the associated adjoint form semisimple Lie group [19] p. 498 Proposition D.40, [24], p. 87, §16.5. Hence there is an automorphism $G \xrightarrow{a} G$ of $G$ which yields $\alpha \mapsto -\alpha$ in the root system. (We can define such an automorphism explicitly as $e_\alpha \mapsto -e_{-\alpha}$ on root vectors in a Chevalley basis; see §10.5 on page 42.) Our automorphism sends $P$ to an opposite parabolic subgroup $P^{\mathrm{op}} := aP$ with $P \cap P^{\mathrm{op}} = G_0$. Letting $G_- := aG_+$, $G_+ \cap G_- = \{1\}$. An open subset of $G$ consists of elements uniquely expressed as a product $p, q \in P \times G_- \mapsto pq \in G$ [15] p. 294, [6] p. 198 Proposition 14.21. Every root system also has an automorphism, traditionally called $w_0$, which belongs to the Weyl group and which interchanges the positive and negative roots of a root system [15] pp. 323–324; it is the unique element of the Weyl group of maximum length. Note that $w_0$ might not reverse the signs of simple roots [15] p. 324. The Weyl group lifts to a group of automorphisms of the Lie group $G$, after perhaps extension by some finite group of order a power of 2.

6.2. **Proof of the cominuscule property.** A Lie group $G$ acts *almost effectively* on a manifold $X$ if the elements of $G$ fixing every point of $X$ form a finite subgroup.

**Lemma 4.** *The complex homogeneous space* $(\breve{X}, \breve{G})$ *is a homogeneously embedded cominuscule subvariety of* $(X, G)$*. If* $X$ *has positive dimension then so does* $\breve{X}$*. If* $(X, G)$ *is almost effective then so is* $(\breve{X}, \breve{G})$*. The Dynkin diagram of* $(\breve{X}, \breve{G})$ *has precisely one connected component for each connected component of the Dynkin diagram of* $(X, G)$ *containing a cross.*

*Proof.* Take notation as above for a flag variety $X = G/P$. It suffices to assume that $(X, G)$ is effective. It also suffices to assume that $(X, G)$ is an irreducible flag variety, as otherwise it is a product of irreducibles and the subgroup $\breve{G}$ is the product of the associated subgroups. The Lie algebra $\mathfrak{z}$ of the center $Z$ of the unipotent radical $G_+$ of $P$ is the sum of the root spaces of the box, i.e. of the $P$-maximal roots. Let $a$ be an automorphism of $G$ which changes the sign of the $P$-grading of the roots; see §6.1 where we constructed one such automorphism. So $a\mathfrak{z}$ is the sum of the root spaces of the roots of minimal $P$-grade, the opposite box. Under bracket, these root spaces generate only root vectors and coroots, up to scaling. Hence the Lie algebra $\breve{\mathfrak{g}}$ of $\breve{G}$ is the sum of some such. A coroot arises only when we bracket the root vectors of the associated root and its negative. So $\breve{\mathfrak{g}}$ contains all of the root vectors of the root system generated by the box and its opposite. It also contains their brackets. Hence it contains the complex semisimple subalgebra with that root system. That subalgebra is generated in the same way, by the box root vectors. Every complex semisimple Lie group is generated from the root vectors of a basis of simple roots [53] p. 43 Theorem 1. So $\breve{\mathfrak{g}}$ contains, and therefore is, that subalgebra, so is complex semisimple. By lemma 3 on the preceding page, $\breve{G}$ is a complex semisimple Lie group. Since its root system lies inside that of $G$, its Cartan subgroup is a subgroup of the Cartan subgroup of $G$, generated by its coroots. Note that $P$ contains the Cartan subgroup of $G$, hence that of $\breve{G}$, and contains $Z$, so contains the parabolic subgroup of $\breve{G}$ generated by the box and its opposite. So this parabolic subgroup fixes the same point $x_0 \in X$ so lies in $\breve{P}$. The $\breve{G}$-orbit $\breve{X}$ of that point $x_0$ is a flag variety $(\breve{X}, \breve{G})$, with stabilizer $\breve{P} = P \cap \breve{G}$ a parabolic subgroup, so connected.

Since $\breve{\mathfrak{g}}$ is a complex subalgebra of $\mathfrak{g}$, $\breve{G}$ is a complex Lie subgroup of $G$, and so $\breve{X} \subseteq X$ is a compact complex submanifold, and $X$ is a smooth projective variety, so $\breve{X} \subseteq X$ is a smooth subvariety.

The vector spaces $\mathfrak{z}$ and $a\mathfrak{z}$ are irreducible $G_0$-modules [32] p. 332 proposition 5.105. Hence $\breve{\mathfrak{g}}$ is a $G_0$-module. As $G_0$ is reductive, $\breve{\mathfrak{g}}$ is a direct sum of $G_0$-invariant subspaces:

$$\breve{\mathfrak{g}} = \mathfrak{z} \oplus \breve{\mathfrak{g}}_0 \oplus a\mathfrak{z}.$$

If $d$ is the maximal $P$-grade, by lemma 2 on page 12, $\mathfrak{z}$ has grade $d$ and so $a\mathfrak{z}$ has grade $-d$, and brackets add grades, so a 3-step grading.

Suppose that $\alpha$ is a $P$-compact root, i.e. a root of $G_0$. Reflection in $\alpha$ is carried out by an element of the Weyl group of $G_0$, and so preserves the $P$-grading. So if $\beta$ is any $P$-maximal root, then so is its $\alpha$-reflection. In other words, reflection in $\alpha$ preserves the box. Recall that an *$\alpha$-root string* is a directed path of roots, each given by adding $\alpha$ to the previous. Reflection in $\alpha$ moves $\beta$ along an $\alpha$-root string.

If that root string has more than one root in it, then $\alpha$ is a difference of $P$-maximal roots: $\alpha = \beta - \gamma$. Their root vectors lie in $\breve{\mathfrak{p}}$. So

$$0 \neq [\mathfrak{g}_\beta, \mathfrak{g}_{-\gamma}] = \mathfrak{g}_{\beta-\gamma} = \mathfrak{g}_\alpha.$$

Both root spaces on the left belong to $\breve{\mathfrak{g}}$. So $\alpha$ is a $\breve{G}$-root, and a difference of $\breve{P}$-maximal roots.

Suppose that every such root string has a single root in it. So it is a root string of length 1. Reflection in $\alpha$ fixes all roots in the box, i.e. all $P$-maximal roots. Hence $\alpha$ is perpendicular to the box. Take $\breve{\Delta}$ to be the root system generated by the box and its opposite. Reflection in $\alpha$ therefore fixes every root in $\breve{\Delta}$. Hence $\alpha$ is perpendicular to every root in $\breve{\Delta}$.

Any decomposition into two root subsystems which are mutually perpendicular arises from factors of the associated semisimple Lie group. Hence the $P$-compact roots divide into (1) those in $\breve{\Delta}$ and (2) those perpendicular to $\breve{\Delta}$. The roots of $\breve{\Delta}$ are those which are differences of $P$-maximal roots, i.e. $\breve{P}$-compact roots. The roots perpendicular to $\breve{\Delta}$ form a root subsystem of the $P$-compact roots giving a complex semisimple subgroup of $G_0$. This subgroup acts trivially on $\breve{\mathfrak{g}}$ because it is generated by root vectors in a perpendicular root subsystem. The root system $\breve{\Delta}$ is graded into the $P$-minimal roots (grade $-1$), differences of $P$-maximal roots (grade 0) and $P$-maximal roots (grade 1). The Lie algebra $\breve{\mathfrak{g}}$ consists of the direct sum of the root spaces of the roots in $\breve{\Delta}$, and their coroots (grade 0). The subalgebra $\breve{\mathfrak{p}} := \mathfrak{p} \cap \breve{\mathfrak{g}}$ consists precisely of the 0 and 1 grades. Note that $\breve{\mathfrak{p}}$ acts on $\mathfrak{z}$ as $\breve{\mathfrak{g}}_0 := \mathfrak{g}_0 \cap \breve{\mathfrak{g}}$, i.e. as a sum of irreducible $\breve{\mathfrak{p}}$-modules, so $\breve{X} = \breve{G}/\breve{P}$ is cominuscule.

Recall that $\mathfrak{z}$ is an irreducible $G_0$-module. If we start at the highest root, we can pass from it via root strings to get to any $P$-maximal root, i.e. anywhere in the box. As we do this, we repeatedly pass between $P$-maximal roots $\alpha, \beta$ by subtracting a $\breve{P}$-compact positive root $\alpha - \beta$. Bracketing a root vector $e_{\beta-\alpha}$ of root $\beta - \alpha$ takes $e_\alpha$ to a nonzero multiple of $e_\beta$. Hence $\mathfrak{z}$ is an irreducible $\breve{\mathfrak{g}}_0$-module, hence an irreducible $\breve{\mathfrak{p}}$-module. As we have seen on page 6, the number of irreducible $\breve{\mathfrak{g}}_0$-modules of $\breve{\mathfrak{g}}/\breve{\mathfrak{p}}$ is the number of irreducible factors of its flag variety. Hence $\breve{X}$ is an irreducible flag variety.

We next prove that $(\breve{X}, \breve{G})$ is almost effective. Since $\breve{X}$ has positive dimension and $\breve{G}$ is simple, the kernel of the action of $\breve{G}$ on $\breve{X}$ is a proper normal subgroup of $\breve{G}$, so finite. □

We restate part of the proof of the previous result, as a separate lemma:

**Lemma 5.** *For any flag variety $(X, G)$, the box generates the $\breve{G}$-root system. The $P$-compact roots divide into (1) the $\breve{P}$-compact roots, which are precisely those roots which are differences of roots of the box and (2) those perpendicular to the $\breve{G}$-roots, forming a root subsystem of the $P$-compact roots giving a complex semisimple subgroup of $G_0$ acting trivially on $\breve{\mathfrak{g}}$.*

**Lemma 6.** *The associated cominuscule subvariety $(\breve{X}, \breve{G})$ of a flag variety $(X, G)$ has the same box as $(X, G)$.*

*Proof.* As above, we can assume that $(X, G)$ is an irreducible flag variety. The definition of associated cominuscule subvariety $(\breve{X}, \breve{G})$ of a flag variety $(X, G)$ has $\breve{G}$ the group generated by $Z, Z^{\text{op}}$, so $\breve{G}$ has Lie algebra containing the Lie algebra $\mathfrak{z}$ of $Z$ and $\breve{P}$ is defined to be $\breve{G} \cap P$, so the Lie algebra of $\breve{P}$ contains $\mathfrak{z}$. The roots of $G$ whose root spaces sum up to $\mathfrak{z}$ are precisely the box roots by lemma 2 on page 12. By lemma 4 on page 16, $\breve{P} \subset \breve{G}$ is a parabolic subgroup. By lemma 5, $(\breve{X}, \breve{G})$ has root system, marked into compact and noncompact, a subsystem of that of $(X, G)$. Moreover, as a graded Lie algebra, $\breve{\mathfrak{g}} = \mathfrak{z} \oplus \breve{\mathfrak{g}}_0 \oplus \mathfrak{z}^{\text{op}}$ is graded into $-1, 0, 1$ grades. So $\breve{\mathfrak{g}}_+ = \mathfrak{z}$ is abelian, so its own center: $\breve{\mathfrak{z}} = \breve{\mathfrak{g}}_+ = \mathfrak{z}$. □

### 6.3. Automorphisms of the associated cominuscule subvariety.

*Example* 9. The flag variety $(X, G) = (B_2/P, B_2)$, where $P \subseteq B_2$ is the Borel subgroup:

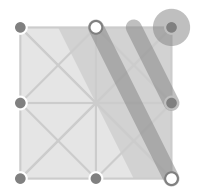

has associated cominuscule subvariety

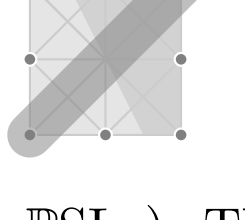

This is a rational curve $(\breve{X}, \breve{G}) = (\mathbb{P}^1, \mathbb{P}\text{SL}_2)$. This rational curve is invariant under not only its automorphism group $\breve{G} = A_1 = \mathbb{P}\text{SL}_2$, which we see in our diagram, but also, as we will see, under rescaling by this root:

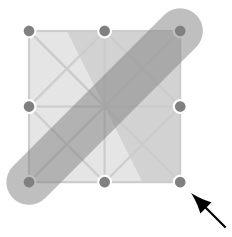

The root is perpendicular to the $\check{G}$ root system, so commutes with the $\check{G}$ root vectors. The vector field on $X$ associated to the root vector of that root vanishes at our chosen point $x_0 \in X$ stabilized by $P$, since the root lies in the Lie algebra of $P$. Since the vector field commutes with those of $\check{G}$, it is $\check{G}$-invariant. Since $\check{G}$ acts transitively on $\check{X}$, our vector field vanishes at all points of $\check{X}$. Hence, in this example, the automorphism group $G'$ of $\check{X}$ as a subvariety of $X$ is larger than the automorphism group $\check{G}$ of $\check{X}$ as a flag variety.

Conversely, consider the root

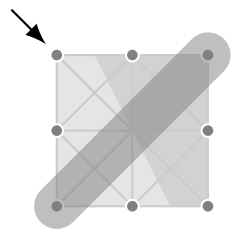

It doesn't belong to the Lie algebra of $P$, so doesn't vanish at $x_0$. Commuting with the root vectors of $\check{G}$, it is $\check{G}$-invariant, so it doesn't vanish at any point of $\check{X}$. If it were tangent to $\check{X}$ at some point of $\check{X}$, $\check{G}$-invariance would force it to be tangent at every point. But every holomorphic vector field on the projective line $\check{X} = \mathbb{P}^1$ has a zero. Indeed, every holomorphic vector field on any flag variety has a zero, since every automorphism has a fixed point [59]. Hence this root vector is a vector field on $X$, nowhere tangent to $\check{X}$.

6.4. **Computing the automorphism Lie algebra.** We return our thoughts to the general flag variety $(X, G)$. By definition of $(\check{X}, \check{G})$, $\check{X} \subseteq X$ is a smooth subvariety. So the subgroup $G' \subseteq G$ leaving $\check{X} \subseteq X$ invariant is a Zariski closed subgroup of $G$, hence linear algebraic. (As we have noted, we will find that, while $G'$ contains $\check{G}$, it is often larger than $\check{G}$ and is not always semisimple.) Let $P' := G' \cap P$, so

$$\check{X} = \check{G}/\check{P} = G'/P'.$$

Clearly $\check{X}, \check{P}, \check{G}$ are connected, so $G'$ is connected just when $P'$ is connected.

**Theorem 2.** *The automorphism group $G'$ of the associated cominuscule subvariety $\check{X} \subset X$ of a flag variety is a complex linear algebraic group with Lie algebra precisely the sum of*

- *the $P$-maximal and $P$-minimal root spaces,*
- *the maximal reductive $\mathfrak{g}_0$, and*
- *the root spaces of all positive roots perpendicular to all $P$-maximal roots.*

*Example* 10. In example 9 on the preceding page, we saw that the third case above is not redundant: there is a positive root perpendicular to all $P$-maximal roots in that example, and this root is not a root of the maximal reductive $\mathfrak{g}_0$, since in that example, $\mathfrak{g}_0$ is abelian.

*Proof.* Clearly $G' \subseteq G$ is a complex linear algebraic subgroup containing $\check{G}$, hence acts transitively on $\check{X}$. So we only need to prove that the Lie algebra $\mathfrak{g}'$ of $G'$ is the sum of the appropriate root spaces. To see that $\mathfrak{g}'$ is a sum of root spaces, we need to prove that it contains the Cartan subalgebra of $\mathfrak{g}$.

Since $\check{X}$ is a projective variety and $G'$ is a linear algebraic group acting transitively on $\check{X}$, the stabilizer $P' := P \cap G'$ of the point $x_0 \in X$ is a parabolic subgroup of $G'$ [6] p. 148, so contains a Borel subgroup of $G'$.

Let $G'_0 \subset G'$ be the connected component of the identity. Since $\check{X}$ is a projective variety and $G'_0$ is a connected linear algebraic group acting transitively on $\check{X}$, the stabilizer $P'_0 := P' \cap G'_0$ of the point $x_0 \in X$ is a parabolic subgroup of $G'_0$ [6] p. 148, hence is connected. Since $\check{X} = G'/P'$ is connected, $P'$ intersects every component of $G'$.

Since $\check{G}$ acts transitively on $\check{X}$, every element of $G'$ is a product of an element of $\check{G}$ and an element of $P$.

Claim: $P'$ is the normalizer of $\check{G}P$ in $P$. Proof: by definition $p \in P'$ if and only if $p \in P$ and $p\check{X} = \check{X}$. The points of $\check{X}$ have the form $gx_0$ for $g \in \check{G}$, unique up to multiplying by an element of $P$, i.e. $pgP = g'P$, for some $g' \in \check{G}$. So $p \in P'$ if and only if $p \in P$ and $p\check{G}P \subseteq \check{G}P$. Since $P'$ is a group, $p \in P'$ just when $p \in P$ and $p\check{G}P = \check{G}P$, or equivalently,

$$p\check{G}Pp^{-1} = \check{G}P.$$

Thus $P'$ contains all elements of $P$ which normalize $\check{G}$. It contains $\check{P}$. Hence it contains $Z = Z_{G_+} \subseteq \check{P}$. The maximal reductive $G_0 \subseteq P$ normalizes $P$ hence $G_+$ hence $Z$, and is invariant under the automorphism $a$. Thus $G_0 \subseteq P'$. Hence $P'$ contains the Cartan subgroup. So the Lie algebra $\mathfrak{p}'$ of $P'$ is a sum of the Cartan subalgebra together with various root spaces. Since $P'$ lies inside $G'$, the Lie algebra of $G'$ is also a sum of the Cartan subalgebra together with various root spaces. We will say that a root $\alpha$ of $G$ is a *$G'$-root* if its root space lies in $\mathfrak{g}'$, a *$P'$-root* if its root space lies in $\mathfrak{p}'$.

The remainder of the proof consists in identifying the $G'$-roots. First, we find some roots which are obviously $G'$-roots. Take a root $\alpha$ perpendicular to all roots in the box, i.e. to all $P$-maximal roots. It is then also perpendicular to their negatives, so to the root system they generate. If $\alpha$ has positive $P$-grade, and $\beta$ is a box root then $\alpha + \beta$ is not a root, by maximality of $\beta$. But then $\beta - \alpha$ is also not a root, or else it lies in a longer root string. So $[\mathfrak{g}_\alpha, \mathfrak{g}_{\pm\beta}] = 0$. So its associated root vector $e_\alpha$ brackets to zero with all of the root vectors $e_\beta$ of all of those roots, and hence of all roots $\beta$ of $\check{G}$. As we argued in §6.3 on page 18, $\alpha$ is $P$-nonnegative if and only if $e_\alpha$ is tangent everywhere to $\check{X}$, hence generates an automorphism of $\check{X}$. So the Lie algebra $\mathfrak{g}'$ of $G'$ contains the root vector $e_\alpha$ of a root $\alpha$ perpendicular to the box roots just when $\alpha$ is $P$-nonnegative, so either $P$-compact or $P$-positive. If $\alpha$ is $P$-compact, $\alpha$ is a root of $\mathfrak{g}_0$, so we may assume that $\alpha$ is $P$-positive, hence positive.

So far, we have affirmed that the roots which belong to the box or its negative or to the maximal reductive or are positive and perpendicular to the box are $G'$-roots. We have also affirmed that the roots which are $P$-negative and perpendicular to the box have root spaces disjoint from $\mathfrak{g}'$, so are not $G'$-roots.

We are left to consider a root $\alpha$ which is not in the box, and $-\alpha$ is not in the box, and $\alpha$ is not a root of the maximal reductive, and $\alpha$ is not perpendicular to the box. We need precisely to prove that $\alpha$ is not a $G'$-root, i.e. the root vector $e_\alpha$ of $\alpha$ is not in $\mathfrak{g}'$.

By definition $P'$ preserves $\check{\mathfrak{g}} + \mathfrak{p}$. Pick a root $\beta^+$ in the box, not perpendicular to our $P$-submaximal $P$-positive root $\alpha$. Note that if two roots have more than a right angle between them, then their sum is a root and their root vectors have nonzero Lie bracket [53] p. 29. The root system generated by $\alpha$ and $\beta^- := -\beta^+$ is of rank 2. Our $P$-grading grades that rank 2 root system. Our pictures above classify every parabolic subalgebra of every rank 2 simple Lie algebra. In those pictures, we see that the root vectors of $\alpha$ and $\beta^-$ have more than a right angle between them and that their sum is a root. Their $P$-grading is at least that of $\beta^-$, $P$-negative. Hence the Lie bracket of their root vectors is not in $\mathfrak{p}$, but also is not $P$-minimal, so not in $\check{\mathfrak{g}}$. Hence their root vectors bracket out of $\check{\mathfrak{g}} + \mathfrak{p}$. Hence the Lie algebra $\mathfrak{g}'$ of $G'$ does not contain the root space of $\alpha$. In other words, $\alpha$ is not a $G'$-root.

On the other hand, if $\alpha$ is $P$-negative, the same argument applies with $\beta^-$ replaced by $\beta^+$. So neither $\alpha$ nor $-\alpha$ are $G'$-roots.

Note that the $P'$-roots are precisely the $P$-maximal roots and the $P$-positive roots perpendicular to them and the $P$-compact roots. Thus

$$\mathfrak{p}' = \mathfrak{g}_0 \oplus \mathfrak{z} \oplus \bigoplus_{\deg_P \alpha > 0}^{\alpha \perp \mathrm{box}} \mathfrak{g}_\alpha .$$

□

**Lemma 7.** *For any flag variety $(X, G)$, the group $G'$ of automorphisms of the associated cominuscule subvariety $\check{X} \subseteq X$ is connected. The subgroup $P' \subseteq G'$ fixing a point of $\check{X}$ is a connected parabolic subgroup.*

*Proof.* Since $\check{G} \subseteq G'$, $G'$ acts transitively on the connected variety $\check{X}$. We saw in the proof of theorem 2 on page 19 that the stabilizer $P' \subset G'$ of any point $x_0 \in \check{X}$ intersects every component of $G'$. Flag varieties are connected and simply connected [1] p. 69 Prop. 7. By exact sequence in homotopy, inclusion $P' \to G'$ is an isomorphism on $\pi_0$. It suffices to prove that $P'$ is connected. The subgroup $G'_+ := G_+ \cap P' \subset G'$ is unipotent. It is therefore connected and isomorphic as an affine variety to complex Euclidean space [16] p. 36 corollary 15.1.11, [18] p. 499, Corollaire 4.

Since $P = G_0 \ltimes G_+$, we can write each element $p \in P$ as a product of elements of $G_0$ and $G_+$. Since $G_0 \subseteq P'$, we can write each element $p \in P'$ as product of elements of $G_0$ and $G'_+$ so $P' = G_0 \ltimes G'_+$ so $P'$ is connected. In the proof of theorem 2 on page 19, we saw that $P'$ is therefore a parabolic subgroup. □

**Theorem 3.** *For any flag variety $(X, G)$, the group $G'$ of automorphisms of the associated cominuscule subvariety $\check{X} \subseteq X$ acts on $\check{X}$ as precisely the biholomorphisms of $\check{X}$ arising from elements of $\check{G}$: the morphisms*

$$\check{G} \to \mathrm{Aut}_{\check{X}}, G' \to \mathrm{Aut}_{\check{X}}$$

*have the same image. The first morphism has kernel a finite central subgroup of the reductive Levi factor $\check{G}_0 \subseteq \check{P}$, central in $\check{G}$. The kernel of the second morphism has Lie algebra generated by*

- *the root vectors of the nonnegative roots perpendicular to the box and*
- *the linear subspace of the Cartan subalgebra on which all of the box roots vanish.*

*Proof.* We can assume that $G$ acts almost effectively on $X$. The automorphism group $G'$ contains $\check{G}$, so its image in the biholomorphisms of $\check{X}$ contains the image of $\check{G}$, which is a quotient $\check{G}/\Gamma$ by a finite subgroup $\Gamma \subset \check{G}_0$ of the reductive Levi factor $\check{G}_0 \subseteq \check{P}$, central in $\check{G}$. The Lie algebra $\mathfrak{g}'$ maps to vector fields on $\check{X}$. The kernel of that map contains all root spaces of the positive roots perpendicular to the box. The kernel intersects $\check{\mathfrak{g}}$ trivially. By theorem 2 on page 19, we need only investigate which elements of $\mathfrak{g}_0$ give rise to a zero vector field.

The box contains the highest root in each factor of $\mathfrak{g}$. The Lie algebra $\mathfrak{g}_0$ contains the Cartan subalgebra of $G$. So $\mathfrak{g}_0$ acts on $\check{\mathfrak{g}}/\check{\mathfrak{p}} = T_{x_0}\check{X}$ as a sum of irreducibles. Each simple factor of $\mathfrak{g}_0$ acts in an irreducible representation. From lemma 5 on page 18, we see how to decoppose $\mathfrak{g}_0$ into factors acting effectively and those acting trivially, and we see the action of the center by bracketing.

Quotienting out the root spaces of the roots perpendicular to the box, and the associated Cartan subalgebra, yields a morphism of Lie algebras $\mathfrak{g}' \to \check{\mathfrak{g}}$, an isomorphism on $\check{\mathfrak{g}} \subseteq \mathfrak{g}'$. Since $G'$ is connected, this map has connected image. □

## 7. The Hasse diagram of a root system

Recall the *Hasse diagram* of a root system [50]. Take an irreducible reduced root system. Pick a choice of basis of simple roots, and an ordering of the simple roots. A *successor* of a positive root $\alpha$ is a positive root of the form $\alpha + \beta$ for a positive simple root $\beta$. The *Hasse diagram* draws dots on the plane, one for each positive root and a line from each root upward to each of its successors, labelled by the number of the simple root by which they differ. The *grade* of a positive root is the total number of simple roots needed to add up to it. Positive roots of equal grade are drawn on the same horizontal line. (Note that the Hasse diagrams of this paper are *not* the Hasse diagrams of flag varieties described in [3] chapter 4. We draw the roots while they draw the Weyl group.)

*Example* 11. The Hasse diagram of $C_4$ is

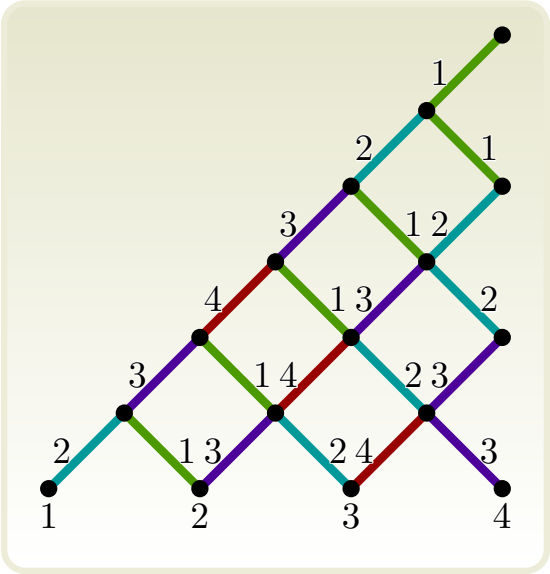

In our pictures, the Hasse diagrams look three dimensional. We won't make any mathematical description of the third dimension. The reader can see that the vertices of the Dynkin diagram form the intersection of a horizontal plane with the three-dimensional object we draw. Imagine the Hasse diagram as "growing" out of the Dynkin diagram. The Hasse diagrams are

$A_r$ Write the simple roots of $A_r$ as $\alpha_1, \alpha_2, \dots, \alpha_r$. Each positive root is a sum of a positive number of *successive* simple roots: $\alpha_i + \alpha_{i+1} + \dots + \alpha_{i+g-1}$ [8] p. 250. This positive root is drawn in the Hasse diagram as a point in the plane. With usual $(x, y)$ Cartesian coordinates, this point is $(x, y) = (2i + g, g)$. An edge labelled $i - 1$ goes up to the left, if $i > 1$. An edge labelled $i + g$ goes up to the right, if $i + g \leqslant r$.

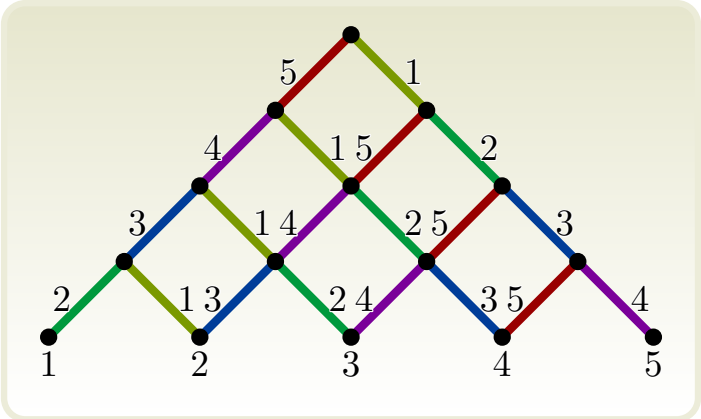

$B_r$ The union of an $A_r$ Hasse diagram and its reflection along the upper right side [8] p. 252.

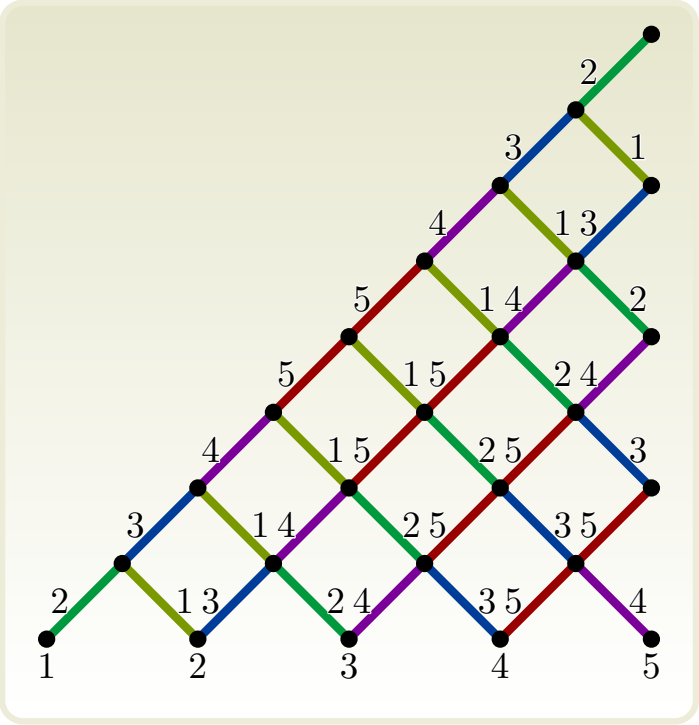

$C_r$ The same as the $B_r$, but all rightward edge labels in the upper copy of $A_r$ diminished by one.

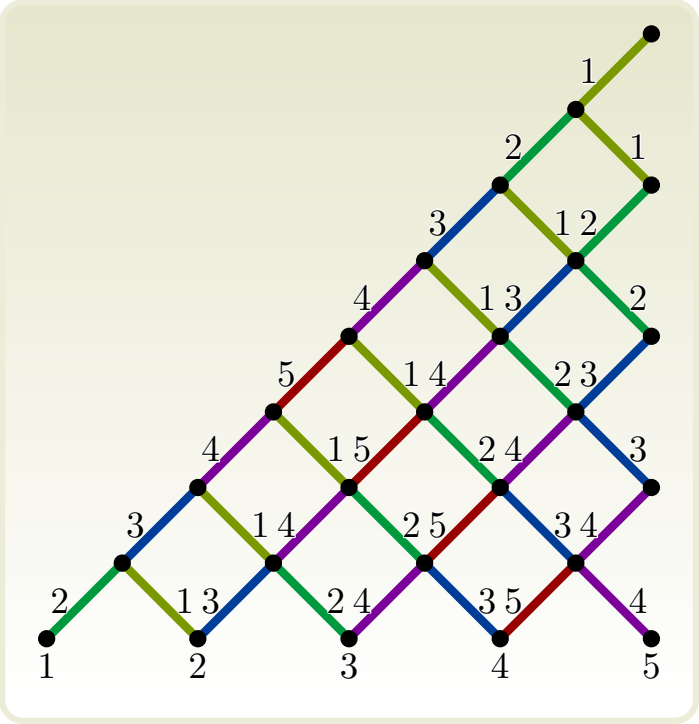

$D_r$ Take two $A_{r-1}$ Hasse diagrams, with one reflected as above. Instead of gluing together, for each matching pair of vertices along the two edges, add a pair of points. Connect one vertex to each with an edge labelled $r-1$ and an edge labelled $r$, to make a square with opposite sides having the same label [8] p. 256.

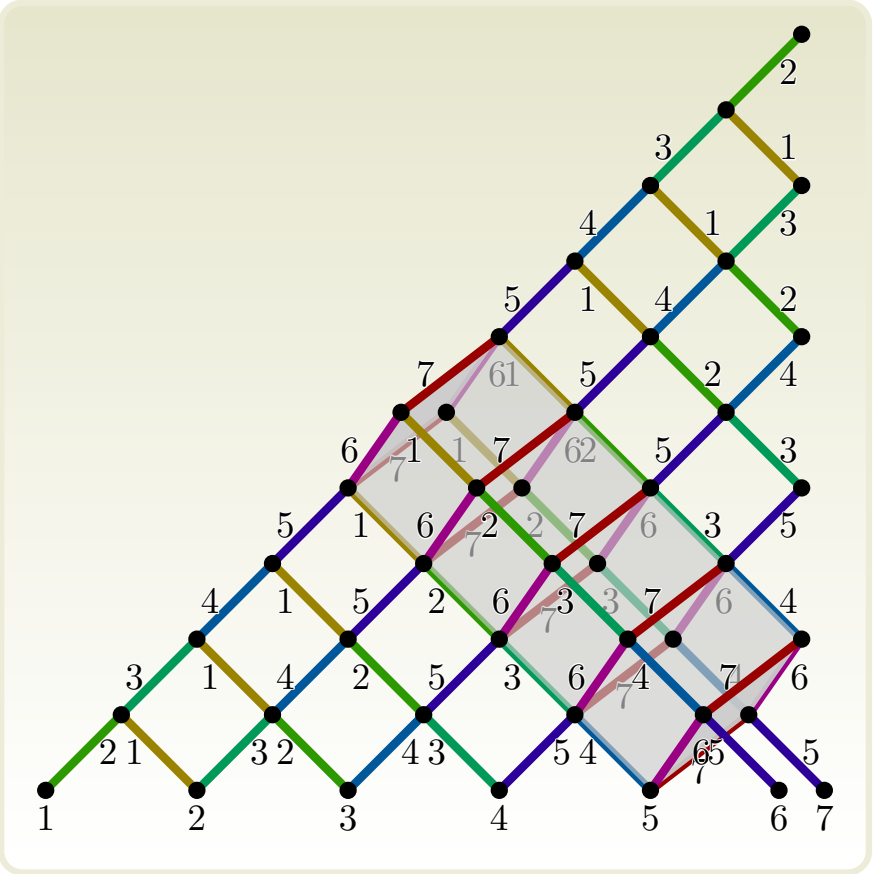

We draw the Hasse diagrams of all of the exceptional Lie algebras in table 7 on page 60.

## 8. Homogeneous vector bundles

Take a flag variety $(X, G)$. Recall that a *homogeneous vector bundle* on $X$ is a holomorphic vector bundle $\mathbf{V} \to X$ equipped with a lift of the $G$-action to vector bundle isomorphisms (a *$G$-linearization* [44] p. 30 §3). Denote the fiber of $\mathbf{V}$ over a point $x \in X$ as $\mathbf{V}_x$ and each vector in $\mathbf{V}_x$ as a pair $(x, v)$.

**Lemma 8** ([15] p. 52 Proposition 1.4.3). *Take a complex homogeneous space* $(X, G)$ *and a point* $x_0 \in X$*, and let* $P := G^{x_0}$*. Every homogeneous vector bundle on* $X = G/P$ *is obtained as an associated vector bundle* $G \times^P V$ *from the* $P$*-module* $V := \mathbf{V}_{x_0}$*. Given two homogeneous vector bundles, every* $G$*-equivariant vector bundle morphism between them arises from a unique* $P$*-module morphism.*

*Proof.* The map

$$(g, v) \in G \times V_{x_0} \mapsto (gx_0, gv) \in \mathbf{V}$$

is $G$-equivariant, and invariant under the right $P$-action

$$(g, v)p := (gp, p^{-1}v).$$

It thus descends to a holomorphic map $G \times^P V \to \mathbf{V}$, linear on fibers. It descends to the identity map on $X$. Hence it is a vector bundle isomorphism, as the reader can easily check. For details and generalizations see [15] p. 52 Proposition 1.4.3, [58] p. 17. ☐

Each filtration of $P$-modules induces a filtration of homogeneous vector bundles. The associated graded $P$-module yields the associated graded vector bundle. Hence we study $P$-modules in some detail.

8.1. **Filtration and grading.** We recall the filtrations and gradings from [15] pp. 238-244; we will refine these in §A on page 47, where we provide definitions of filtered and graded objects as needed. Take a flag variety $(X, G)$. The root system of $G$ is graded: grade each root by writing it as a linear combination of simple roots and taking as grade the sum of the coefficients of the noncompact simple roots, i.e. the number of noncompact simple root summands. The Lie algebra $\mathfrak{g}$ of $G$ is graded: $\mathfrak{g}_i$ is the direct sum of the root spaces of grade $i$, together with the Cartan subalgebra, if $i = 0$. The Lie algebra $\mathfrak{g}$ of $G$ is also filtered:

$$\mathfrak{g}^i = \bigoplus_{j \geqslant i} \mathfrak{g}_j,$$

a finite sum. Note that $\mathfrak{p} = \mathfrak{g}^0$ and $\mathfrak{g}_+ = \mathfrak{g}^1$. The Lie algebra $\mathfrak{p}$ of $P$ is graded and filtered in the same way, as it is invariant under the Cartan subgroup so is also a direct sum of root spaces. This induces a filtered $P$-module structure on $\mathfrak{g}, \mathfrak{p}, \mathfrak{p}^\vee, \mathfrak{g}/\mathfrak{p}, \mathfrak{g}_+$, hence on the vector bundles associated to these $P$-modules.

8.2. **The Hasse diagram of a flag variety.** The *Hasse diagram* of a flag variety $(X, G)$ is the Hasse diagram of the Lie group $G$, but erase the lines labelled by noncompact simple roots.

*Example* 12. Compare 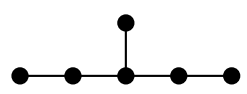 to 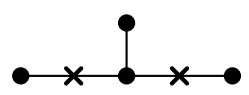:

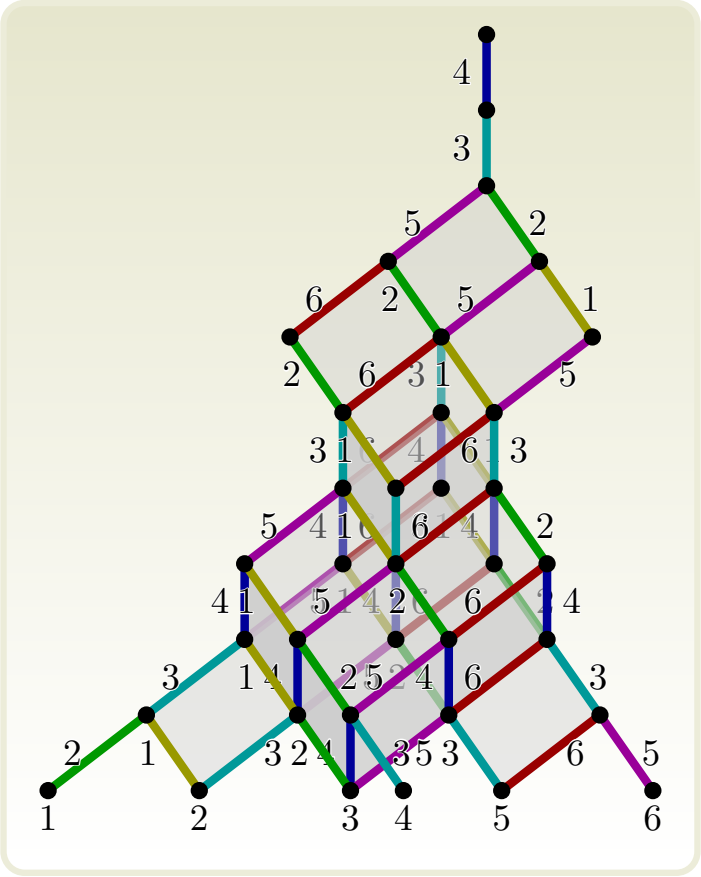

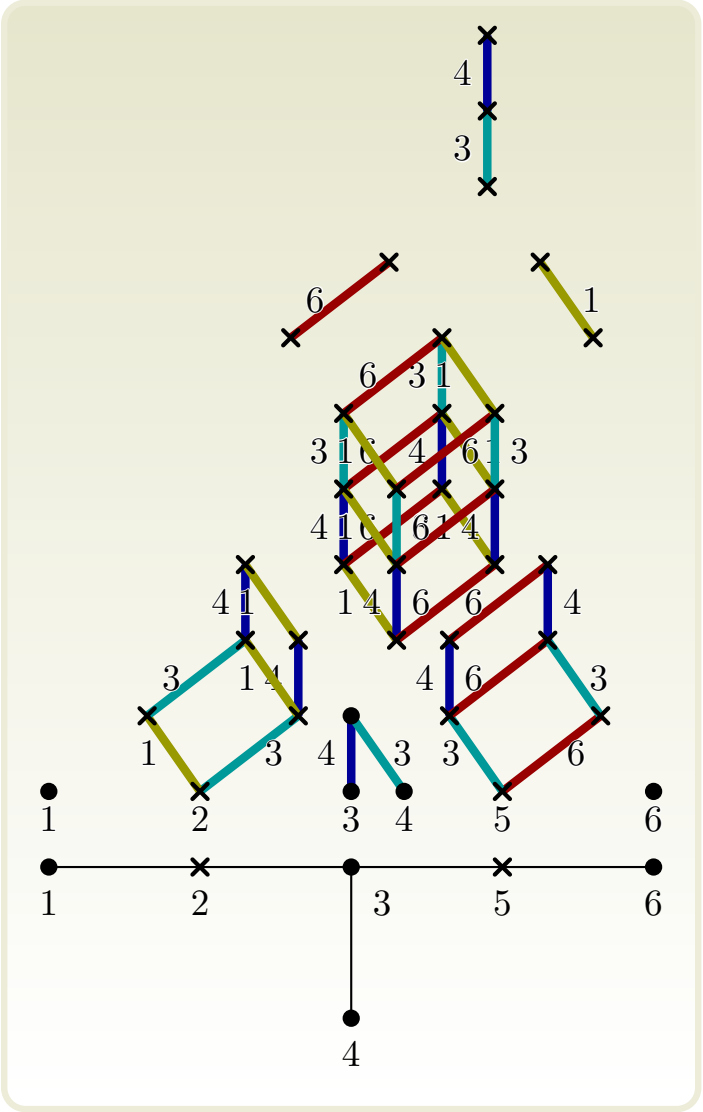

The compact roots (simple or not) we draw as dots, and all other roots we draw as crosses. The reason we draw these diagrams is to unveil as much as we can about the tangent bundles $TX$ of flag varieties $(X, G)$. We will see that we are drawing the decomposition of the associated graded vector bundle of the tangent bundle into invariant subbundles. Again we stress that this Hasse diagram is *not* the Hasse diagram associated to the Weyl group [3] chapter 4.

The tangent bundle arises, as does every homogeneous vector bundle on a flag variety by lemma 8 on the preceding page, as the homogeneous vector bundle associated to a $P$-module. For the tangent bundle, this $P$-module is $\mathfrak{g}/\mathfrak{p}$ [54] p. 188, theorem 3.15. This $P$-module is the dual $P$-module to $\mathfrak{g}_+$ [15] p. 239 Proposition 3.1.2. Its $P$-invariant subspaces are complicated, forming an elaborate maximal $P$-invariant filtration. We pass to the associated graded $P$-module to simplify the filtration to a grading. The associated graded $P$-module $\operatorname{gr}_\bullet(\mathfrak{g}/\mathfrak{p})$ is isomorphic to $\mathfrak{g}/\mathfrak{p}$ as a $G_0$-module, since $G_0$ acts on every $P$-module as a direct sum of $G_0$-irreducibles, so the filtration becomes trivial. Since $\mathfrak{g}_+$ acts by moving into the next filtered subspace, it acts trivially on the associated graded $P$-module. Since $G_+$ is connected, it also acts trivially on the associated graded $P$-module. So the associated graded $P$-module of $\mathfrak{g}/\mathfrak{p}$ is precisely the decomposition of $\mathfrak{g}/\mathfrak{p}$ into the direct sum of $G_0$-modules. The tangent bundle $TX$ has associated graded vector bundle arising from this decomposition of $\mathfrak{g}/\mathfrak{p}$ into its direct sum of $G_0$-modules. But $\mathfrak{g}/\mathfrak{p}$ is the sum of the negative root spaces with nonzero noncompact component. We break this sum into a sum of $G_0$-modules, i.e. broken up by weights of the $P$-compact roots. So in our drawings, the $P$-negative roots are connected by lines just when they differ by a $P$-compact root, so lie in a $P$-compact root string, hence in the same irreducible $G_0$-module. To each connected component of the Hasse diagram, we associate the $P$-module which is the direct sum of its root spaces. Drawing only the components whose roots are crossed, this is precisely the decomposition of $\mathfrak{g}/\mathfrak{p}$ into irreducible $G_0$-modules. The components whose roots are uncrossed form the Hasse diagram of the compact roots, i.e. of the semisimple Levi factor of $P$, as we will prove in corollary 1 on page 55. Hence the crossed root components in our pictures are the $G$-invariant decompositions of the associated graded vector bundle of $TX$.

Changing sign of each root takes positive roots to negative roots. It is an isometry, arising from our automorphism $a$ as above. It is a mirror reflection of the Hasse diagram. So when we draw the Hasse diagram, the reader can think that we are drawing positive or drawing negative roots, as required.

*Example* 13. Compare  to :

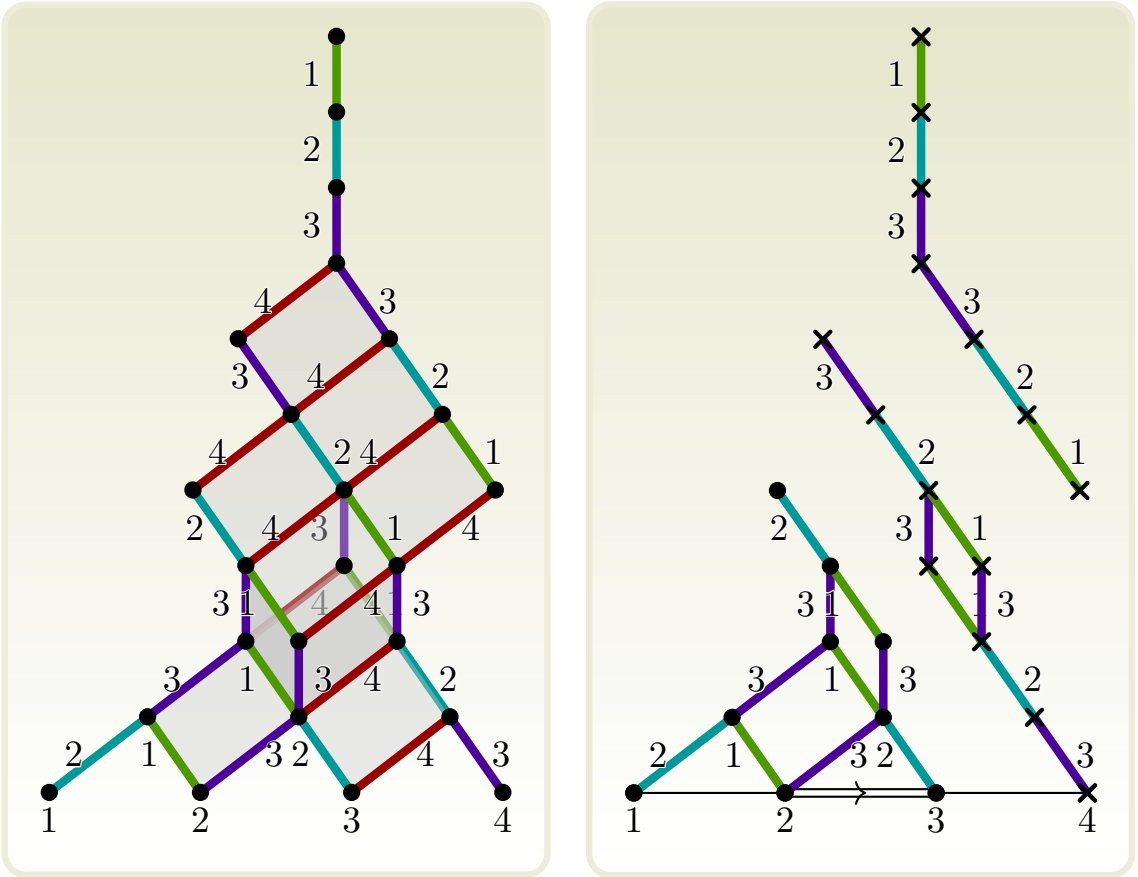

We can see the box: it consists of the 7 roots attached to the highest root:

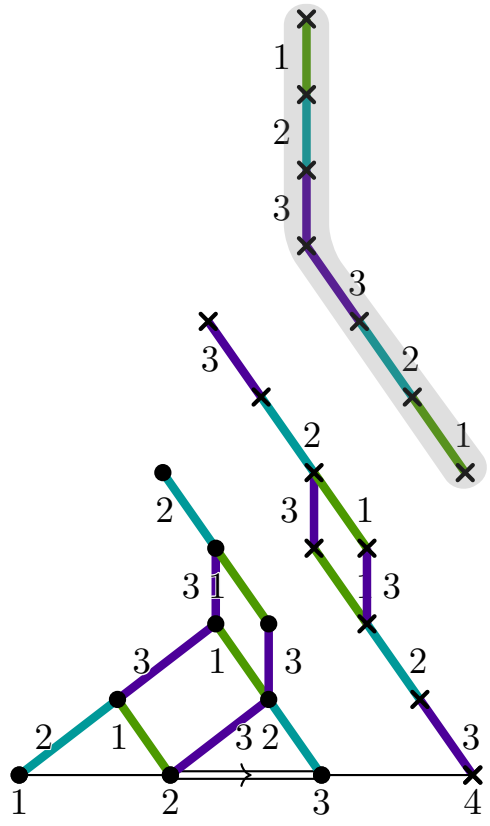

8.3. **The Hasse diagram of a cominuscule variety.** In table 8 on page 62, we draw the Hasse diagrams of the cominuscule varieties.

**Lemma 9** ([3] p. 27, table 3.2, [4] p. 120)**.** *Every irreducible cominuscule variety has a single crossed node in its Dynkin diagram. A positive root, expressed as a sum of simple roots, has coefficient of the associated crossed simple root either* 0 *or* 1*.*

We spot the crossed node immediately in the list of Dynkin diagrams of cominuscule varieties in table 8 on page 62.

*Proof.* Suppose that $(X, G)$ is an irreducible flag variety whose Dynkin diagram has more than one noncompact root. Each connected component of the Hasse diagram of $(X, G)$ is a set of roots whose root spaces sum to a vector space. This vector space is a $G_0$-module by invariance under the Cartan subgroup and the root vectors of the positive simple roots [3] p. 26. This module is irreducible by connectivity.

Recall that every positive $G$-root is obtained by repeated sum of simple roots, stepping up the Hasse diagram of $G$ [53] p. 32 proposition 5. Walk around the Hasse diagram of $G$. Every upward step we take along a noncompact root creates another step in the filtration of $\mathfrak{g}_+$. Dually it creates another step in the dual filtration of $\mathfrak{g}/\mathfrak{p}$. It passes out of a component of the Hasse diagram of $(X, G)$. The highest root contains a positive linear combination of all simple roots. We can only reach it by taking at least one step from each $P$-noncompact root. So the presence of two or more noncompact $P$-roots ensures that $\mathfrak{g}/\mathfrak{p}$ is not an irreducible $P$-module. Similarly, if we can step up the Hasse diagram taking two steps upward that use the same root, then $\mathfrak{g}/\mathfrak{p}$ is not an irreducible $P$-module. □

**Lemma 10.** *Up to possible relabeling of the roots, the boxes of the irreducible cominuscule varieties are (as drawn in table 8 on page 62):*

$A_r$ *Crossed node $k$: the rectangular box of points of the $A_r$ Hasse diagram $\geqslant$ the crossed root in the Hasse diagram ordering. Labels decrease $k-1, k-2, \dots, 1$ along the lower left side, and increase $k+1, k+2, \dots, r$ along the lower right side.*

$B_r$ *The line segment of points above $1$ in the Hasse diagram ordering, $2r-1$ points in all, with labels $2, 3, \dots, r-1, r, r, r-1, \dots, 3, 2$.*

$C_r$ *The triangle of points above $r$ in the Hasse diagram ordering. Thus a copy of the $A_r$ reflected Hasse diagram with rightward labels diminished by $1$.*

$D_r$ 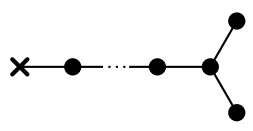 *Above $1$ in the Hasse diagram ordering, a line segment with labels $2, 3, \dots, r-2$, then a square with labels $r-1, r$ on opposite sides, then a line segment with labels $r-2, r-3, \dots, 3, 2$.*

$D_r$ 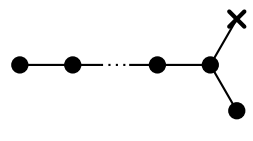 *The triangle above $r$ in the Hasse diagram ordering, i.e. a copy of the $A_{r-1}$ reflected Hasse diagram.*

$D_r$ 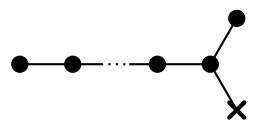 *Isomorphic to the previous.*

$E_6$ *As drawn in table 8 on page 62.*

$E_7$ *As drawn in table 8 on page 62.*

*In particular, each box, as a labelled Hasse diagram, up to label permutations, uniquely determines its cominuscule variety $(X, G)$.*

*Proof.* Billey & Lakshmibai [4] p. 120 classify the irreducible cominuscule varieties by their Dynkin diagrams, as above. The classification of compact Hermitian symmetric spaces is due to Cartan; see [35] p. 379 Proposition 8.2, [3] p. 26. It requires some effort to see that it coincides with the classification of cominuscule varieties.

From the classification, we see that each irreducible cominuscule variety has a unique crossed (i.e. $P$-noncompact) root. When the highest root is written as a linear combination of simple roots, Billey & Lakshmibai prove that an irreducible flag variety is cominuscule just when that unique crossed root appears in that linear combination with coefficient 1 [4] p. 120. The Lie algebra of the unipotent radical is the sum of the root spaces. Except for $A_1$, the Hasse diagram has at least two components: it has some components from the uncrossed roots, and some with that coefficient equal to 1, including the highest root, so including the box. In particular, the crossed simple root itself has that coefficient equal to 1. We know the Hasse diagrams of the simple Lie groups, and we know which root is crossed, so we can cut it out of the diagrams of the simple Lie groups, to obtain table 8 on page 62. □

*Example* 14. We illustrate in this example that the associated cominuscule subvariety does not map functorially. The inclusion of root systems of $E_7$ to $E_8$ yields a unique complex Lie group inclusion from a form of $E_7$ to the adjoint form of $E_8$. It therefore maps flag varieties

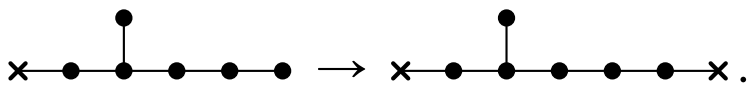

In section C on page 60, we can see the Hasse diagrams of $E_7$ and $E_8$. By removing the relevant edges, we see that the boxes of either of these flag varieties consist of a single ×. Hence both flag varieties have associated cominuscule subvariety a projective line ×. In the Hasse diagrams, we see that the box of the $E_7$-variety maps to another × lower down, not to the box of the $E_8$-variety. So the associated cominuscule subvariety of the $E_7$-variety maps to a rational curve tangent to that $E_8$-invariant subsheaf of the tangent bundle, not to the associated cominuscule subvariety of the $E_8$-variety. So the associated cominuscule subvariety is not functorial.

## 9. The algorithm is correct

We restate our main theorem, theorem 1 on page 3:

**Theorem.** *Take any flag variety. Draw its extended Dynkin diagram. Erase every component not containing a crossed node. Remove the crossed nodes and the edges adjacent to them. Erase all components of this Dynkin diagram which are not connected to an affine node. Change every affine node, drawing it as a crossed node, i.e. instead of drawing it as a hollow node ∘, draw it as ×. This is the Dynkin diagram of the associated cominuscule subvariety.*

We will now prove this theorem, i.e. that our algorithm for computing the Dynkin diagram of the associated cominuscule subvariety works.

### 9.1. Reducing to the irreducible.

If $(X, G)$ is a reducible flag variety, then the Dynkin diagram of $(X, G)$ is the disjoint union of the Dynkin diagrams of the factors. The associated cominuscule subvariety $(\check{X}, \check{G})$ of $(X, G)$ is the product of associated cominuscule subvarieties of the factors. So the Dynkin diagram of $(\check{X}, \check{G})$ is the disjoint union of the Dynkin diagrams of those associated cominuscule subvarieties of those factors, by lemma 4 on page 16. So, to prove the correctness of our algorithm, we can assume that $(X, G)$ is an irreducible flag variety.

### 9.2. Finding the box.

**Lemma 11.** *Take an irreducible flag variety $(X, G)$. Take the Hasse diagram of $G$. Knock out the edges which are labelled by crossed roots, drawing the Hasse diagram of $(X, G)$. In the remaining graph, the box is precisely the connected component of the highest root.*

*Proof.* By lemma 2 on page 12, the box is the set of $P$-maximal roots. Knocking out those edges, again by lemma 2 on page 12, the box is connected and contains the highest root. □

*Example* 15. Look at the Hasse diagram of the group $B_8$

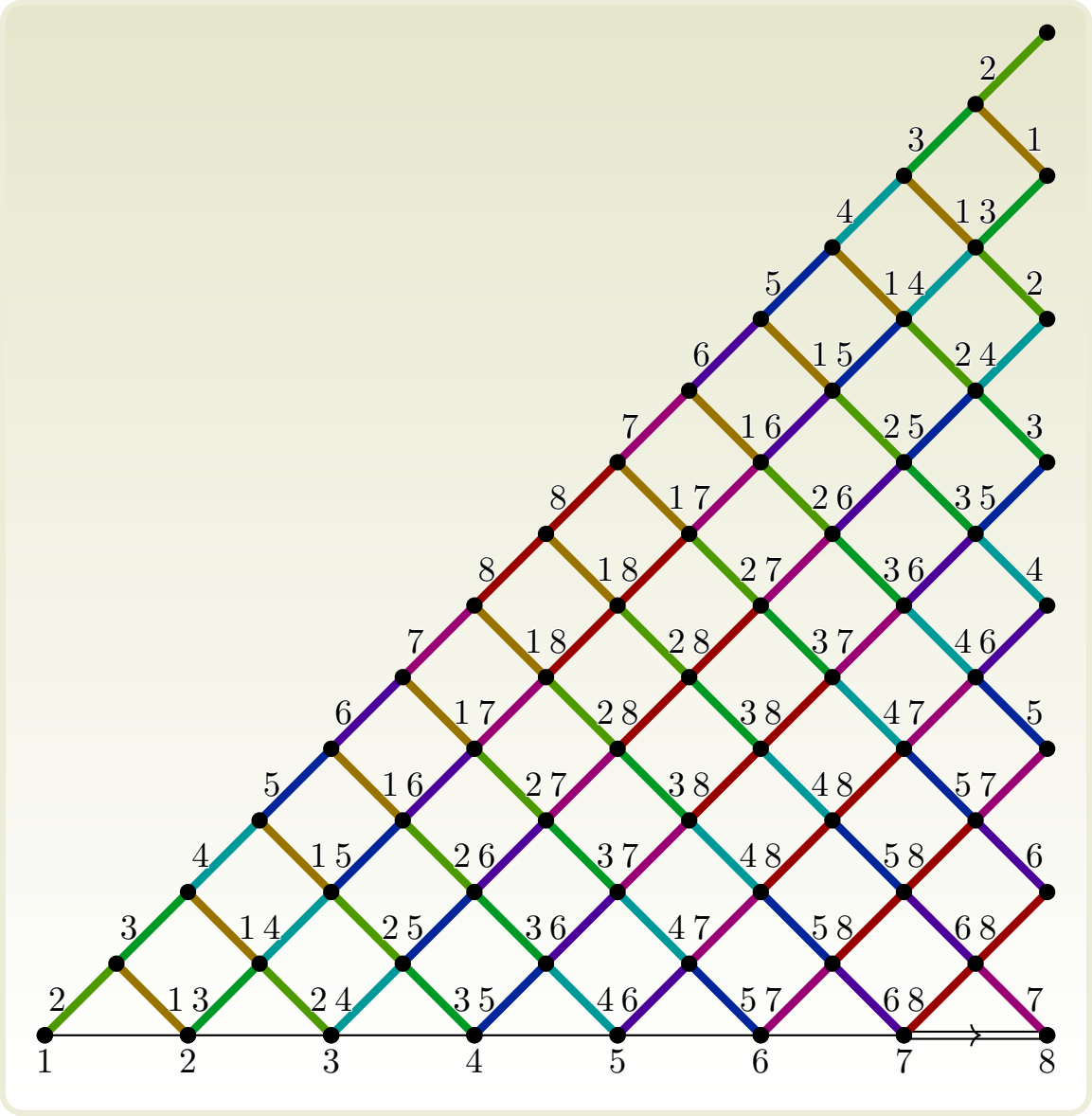

Compare the Hasse diagram of its flag variety :

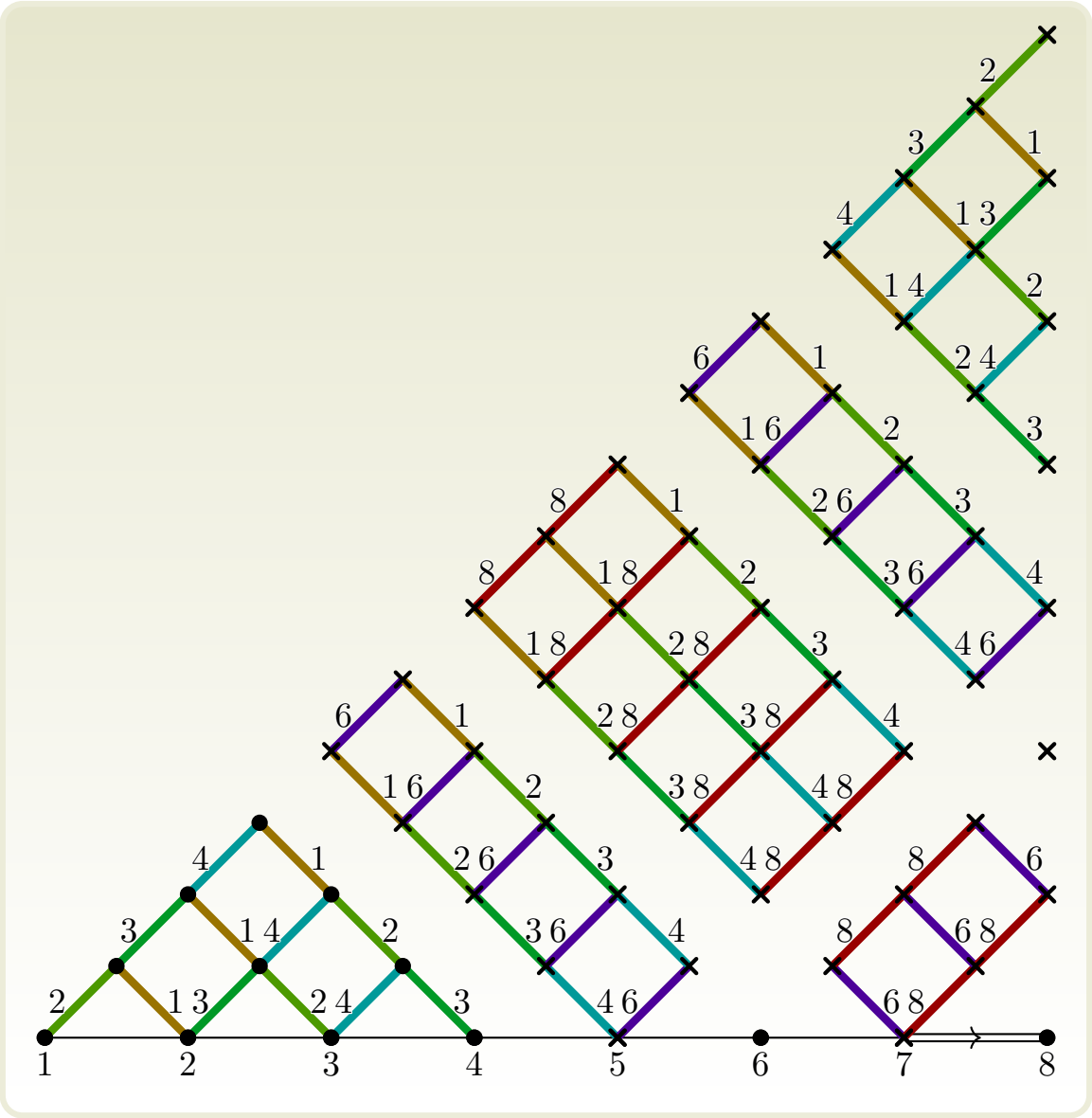

Knock out the edges labelled by crossed roots to pass from the first picture to the second. We see the box: the roots attached to the highest root:

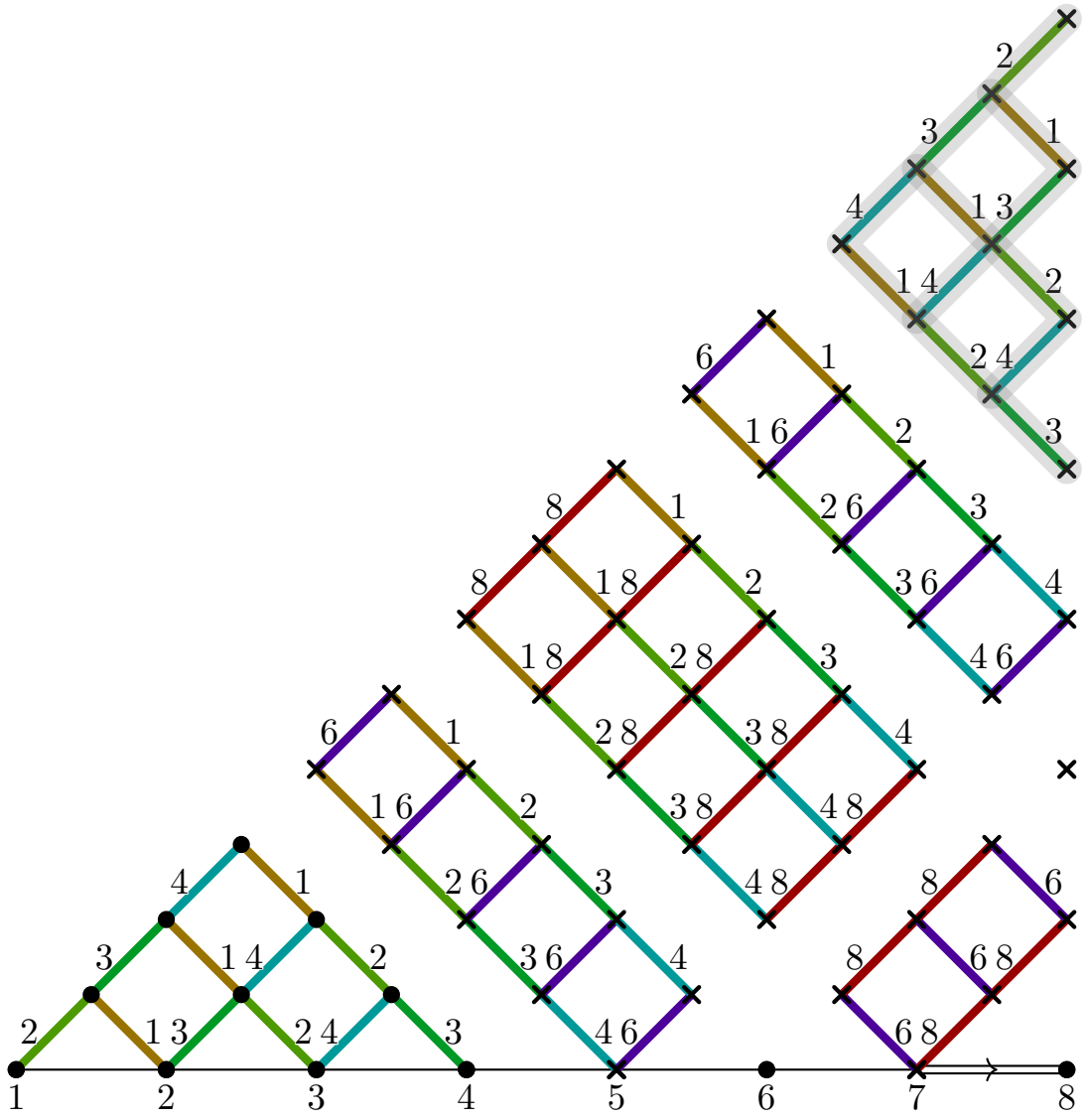

Each edge of the box is labelled by a node $\ell$ of the Dynkin diagram corresponding to a compact root. Consider the nodes which arise as these labels. In our example, these nodes are $1, 2, 3, 4$. Each node $\ell$ represents a $P$-compact simple root $\alpha_\ell$, which is the difference between two roots of the box, the vertices which this edge meets.

By lemma 5 on page 18, the $G$-simple roots arising as the edges of the box of an irreducible flag variety $(X, G)$ are precisely the $\breve{P}$-compact simple roots. In our example, these are the roots $\alpha_1, \alpha_2, \alpha_3, \alpha_4$. The roots $\alpha_6, \alpha_8$ generate a subgroup of $G$ which commutes with $\breve{G}$.

9.3. **The lowest box root.** Let $\breve{\alpha}$ be the lowest box root: in our example,

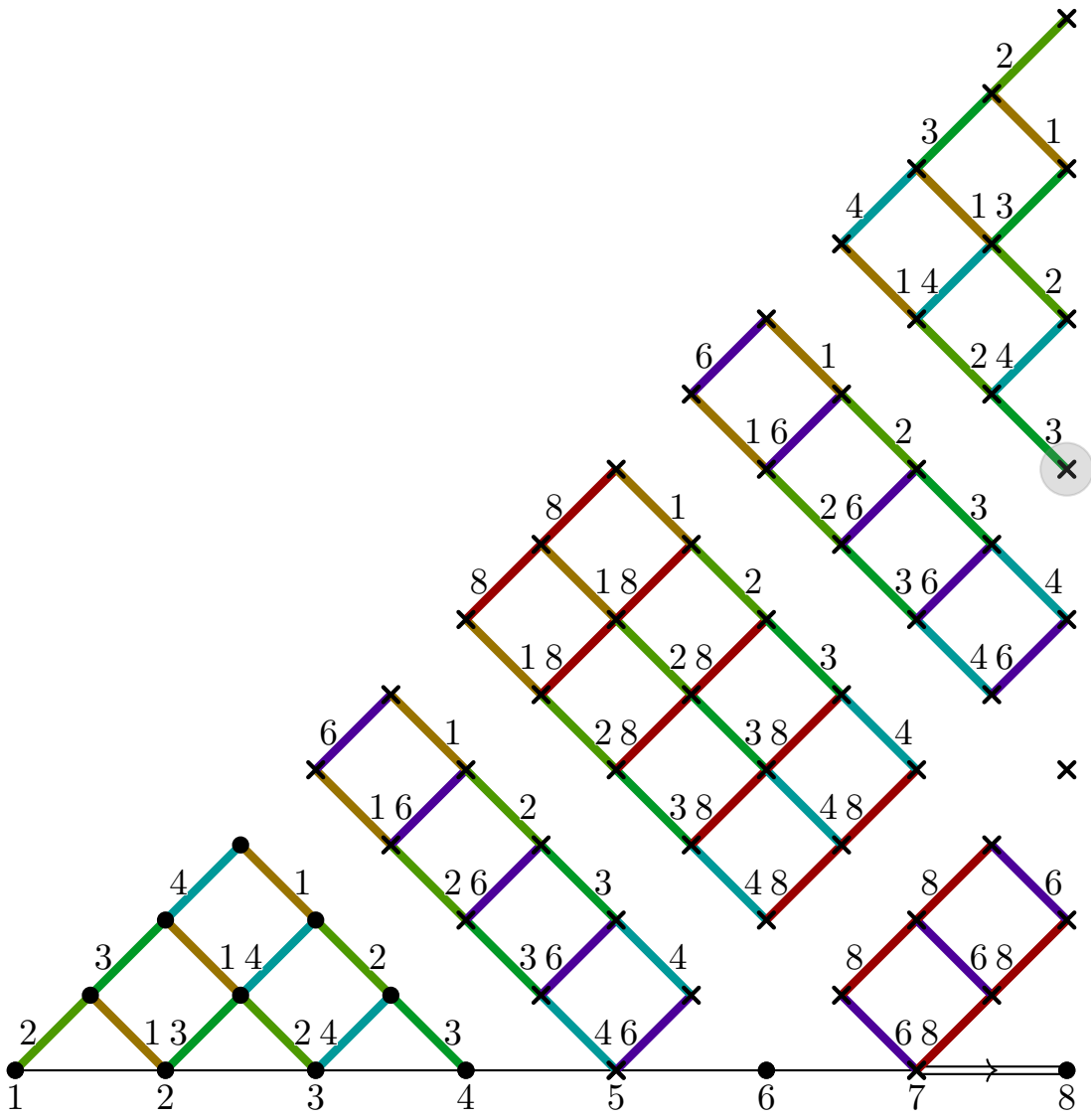

**Lemma 12.** *For any irreducible flag variety $(X, G)$, the lowest box root $\breve{\alpha}$ is the noncompact simple root of $\breve{P}$.*

*Proof.* Recall that the $\breve{G}$-root system is the subsystem of the $G$-root system generated by the $P$-maximal roots. The $\breve{P}$-compact roots are among the $P$-compact roots: precisely those not perpendicular to all $\breve{G}$-roots. For any root system, we can move downwards from the highest root via steps, each step subtracting a simple root. The steps that stay inside the box are $\breve{P}$-compact simple roots, since the $\breve{P}$-compact simple roots are differences of the box roots by lemma 5 on page 18. So we can move around the box from the highest root downwards to $\breve{\alpha}$ along a sequence of $\breve{P}$-compact roots. This is true either in the Hasse diagram of $(\breve{X}, \breve{G})$ or in the Hasse diagram of $(X, G)$. So $\breve{\alpha}$ is the lowest box root in the Hasse diagram of $(\breve{X}, \breve{G})$. In the $(\breve{X}, \breve{G})$-diagram, we can always go downward along compact roots until we hit the crossed root, since the box is connected and consists of all of the roots with positive coefficient of the crossed root. So $\breve{\alpha}$ is the crossed root of $(\breve{X}, \breve{G})$. □

9.4. **The Levi Dynkin diagram.** Take the Dynkin diagram of any irreducible flag variety. Label its nodes distinctly, with any set of labels. Construct the Hasse diagram from the Dynkin diagram, labelling the edges according to the labelling of the Dynkin diagram. Remove the crossed roots from the Dynkin diagram.

*Example* 16.

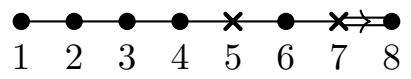

becomes

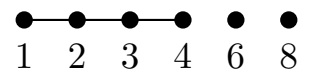

The resulting Dynkin diagram is the Dynkin diagram of the semisimple Levi factor $G_{0,ss}$ of the parabolic subgroup $P$, call it the *Levi Dynkin diagram.*

**Lemma 13.** *Take an irreducible flag variety $(X, G)$. The Dynkin diagram of its associated cominuscule subvariety $(\breve{X}, \breve{G})$ is drawn by taking that of $(X, G)$, cutting out all crossed nodes and adding one new crossed node (which is then the only crossed node), and perhaps adding some edges from the new crossed node to some other nodes. That crossed node has an edge (of multiplicity to be determined in section 9.5) drawn to an uncrossed node just when that uncrossed node has its label appearing on an edge of the box going up from the lowest box root $\breve{\alpha}$. The resulting Dynkin diagram has precisely one component with a crossed node. That component is the Dynkin diagram of the associated cominuscule subvariety. It consists precisely of the new crossed node and the uncrossed nodes whose labels arise along the edges of the box.*

*Proof.* Every Dynkin diagram is drawn by taking the nodes to represent the simple roots. Edges appear precisely when two roots have root spaces with nonzero bracket.

By lemma 5 on page 18, the $\breve{P}$-compact roots are precisely those $P$-compact roots which are differences of box roots, i.e. which appear as labels in the edges of the box. All other $P$-compact roots are perpendicular to the box roots, forming a perpendicular root system, acting trivially on $\breve{X}$, hence perpendicular to the $\breve{G}$-roots.

By lemma 4 on page 16, $(\breve{X}, \breve{G})$ is an irreducible cominuscule flag variety. So it has precisely grades $-1, 0, 1$. So it has at least one $\breve{G}$-simple root not of grade zero, so of grade one.

If there are two $\breve{G}$-simple roots of grade one, the sum of the root spaces of the grade one roots (i.e. of the box roots) is reducible over the reductive Levi factor of $P$ by theorem 6 on page 53. By lemma 2 on page 12, the box graph is connected, i.e. this representation is irreducible.

So there is a unique $\breve{G}$-simple root of grade one. Call that $\breve{G}$-simple root $\breve{\alpha}$. By lemma 5 on page 18, the remaining $\breve{G}$-roots of grade one are obtained from $\breve{\alpha}$ by repeatedly adding compact simple roots which are differences of box roots, i.e. repeatedly adding $\breve{P}$-compact simple roots, i.e. $P$-compact simple roots whose labels lie along the edges of the box, i.e. adding $P$-compact simple roots which are not perpendicular to the box roots. Hence in the $(\breve{X}, \breve{G})$-Dynkin diagram, there is one node for each $P$-compact simple root not perpendicular to the box roots, and there is one node for $\breve{\alpha}$. These nodes sit in a connected Dynkin diagram because $(\breve{X}, \breve{G})$ is irreducible. The $P$-compact roots perpendicular to all of the box roots generate a root system whose root vectors act trivially on $\breve{X}$, so we have the complete $(\breve{X}, \breve{G})$-Dynkin diagram.

The crossed $(X, G)$-nodes do not appear in the $(\breve{X}, \breve{G})$-Dynkin diagram, because we have all of the $(\breve{X}, \breve{G})$-nodes drawn. Consider the $(\breve{X}, \breve{G})$-nodes which are not crossed. These are the simple $(X, G)$-roots which arise as labels in the box. By the standard rules of construction of Dynkin diagrams, in the $(\breve{X}, \breve{G})$-Dynkin diagram, the new crossed node (corresponding to $\breve{\alpha}$) is connected to another node $\alpha_j$ (hence to an uncrossed node) precisely when their root vectors have nonzero bracket. This occurs precisely when the root $\breve{\alpha}$ has an edge labelled by $j$. □

*Example* 17. Suppose that the Dynkin diagram of $(X, G)$ is:

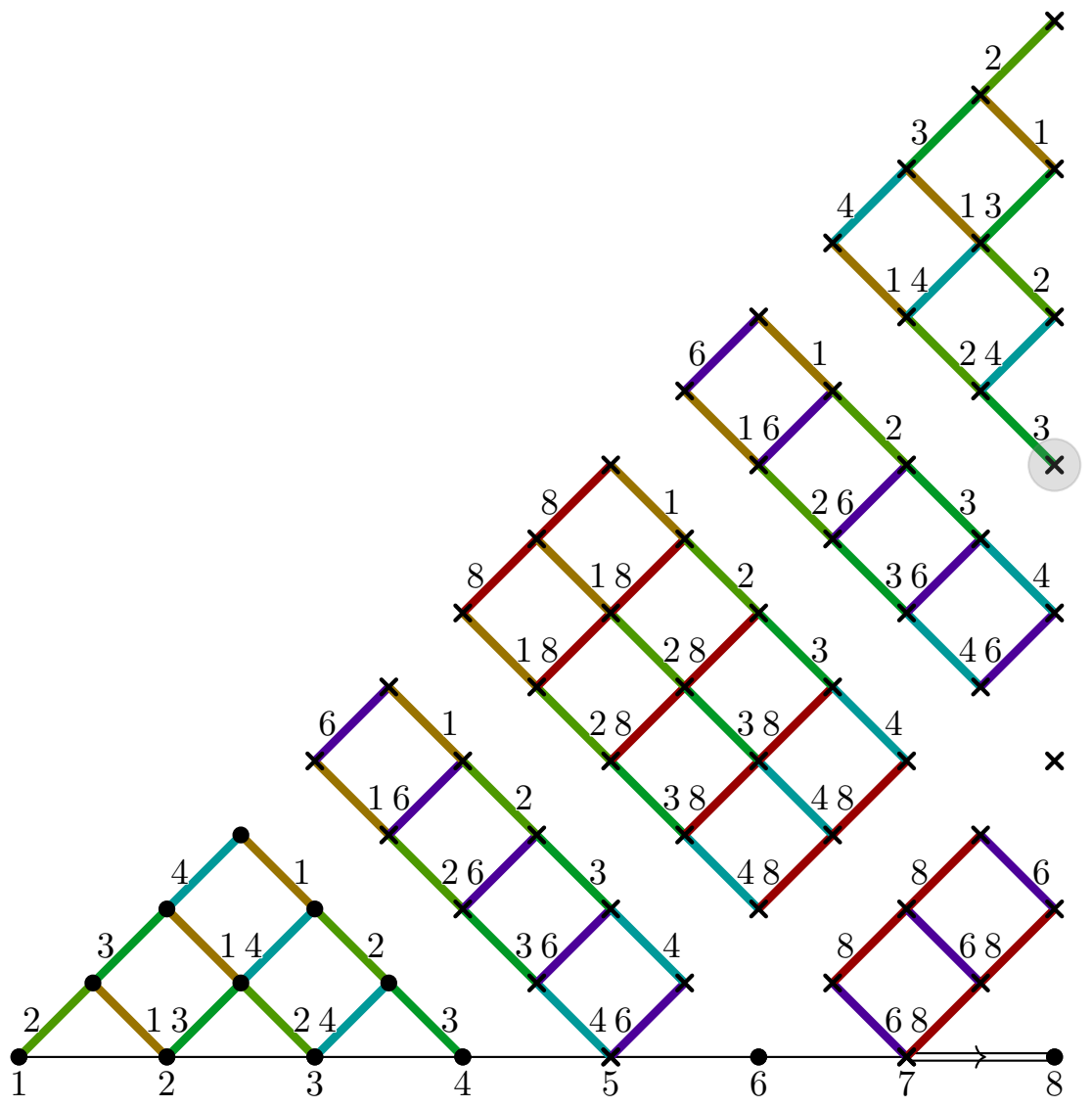

Then the Dynkin diagram of $(\breve{X}, \breve{G})$ is with the edge connecting the new × to node 4 having as yet unknown multiplicity.

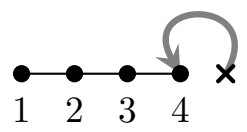

9.5. **Edge multiplicities.** To finish our picture of the Dynkin diagram of the associated cominuscule subvariety of a flag variety, we need to decide multiplicities of the edges we have just drawn. Our aim is to prove that these edge multiplicities are precisely those of the affine root in the extended Dynkin diagram.

Pick a root $\alpha$ in a Hasse diagram of a flag variety. Take an edge adjacent to it, say labelled $b$, with associated simple root $\alpha_b$. The *b-multiplicity* of that root is the length of the $\alpha_b$-root string starting at $\alpha$, i.e. the largest number of edges labelled $b$ which lie in a sequence through $\alpha$.

Start at the lowest box root $\breve{\alpha}$. For any $b$, the *$b$-multiplicity* of $(\breve{X}, \breve{G})$ is the $b$-multiplicity of $\breve{\alpha}$.

*Example* 18. In the Hasse diagram of :

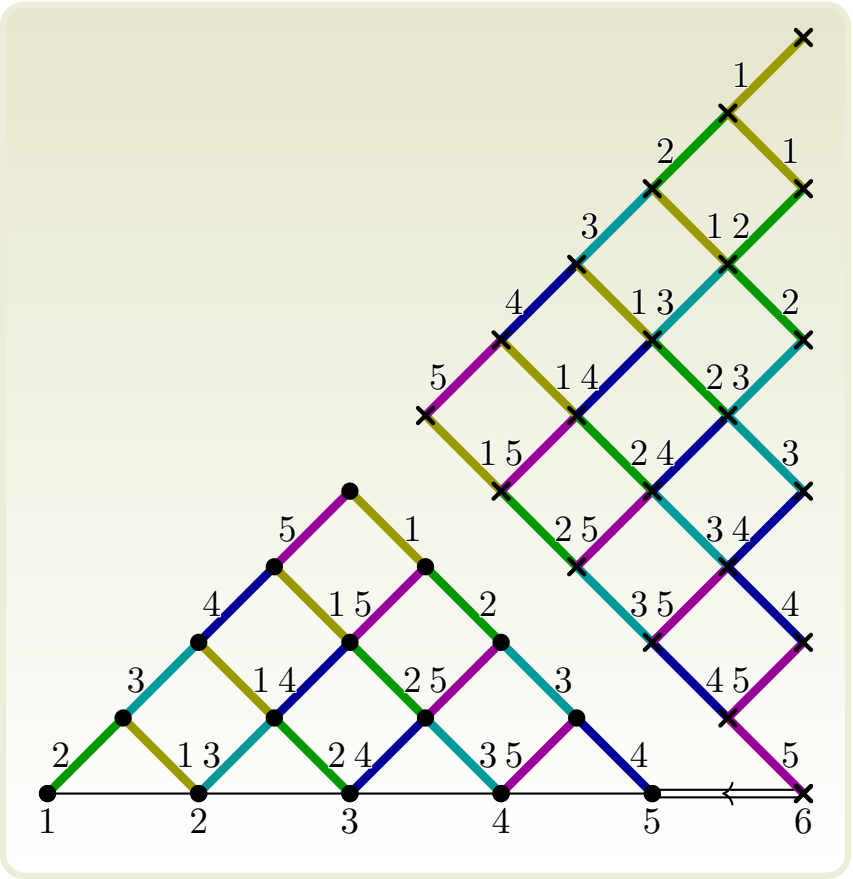

we can move up from the × node (node 6) by following root 5 twice, so the 5-multiplicity is 2. Since there is only one edge from the × node, all other multiplicities are zero.

**Lemma 14.** *Take an irreducible flag variety $(X, G)$. Draw the Dynkin diagram of its associated cominuscule subvariety $(\breve{X}, \breve{G})$, and the box. The Dynkin diagram has a single crossed node, corresponding to the lowest box root $\breve{\alpha}$. For each uncrossed root $\alpha_b$, the multiplicity in that Dynkin diagram of the edge from the crossed node to the uncrossed node labelled $b$ is the $b$-multiplicity, which is either $0, 1$ or $2$. If that multiplicity is $2$ then the arrow on that edge points from the crossed node to the uncrossed node.*

*Proof.* The $\alpha_b$-root string through $\breve{\alpha}$ is the same in the $G$-root system, the $\breve{G}$-root system, or the rank 2 root system generated by $\breve{\alpha}, \alpha_b$ [24] p. 145 §25.1, Proposition. The Cartan integers involved are the same, and any Chevalley basis for the $G$-root system extends to a unique Chevalley basis for the $\breve{G}$-root system, which extends a unique Chevalley basis for the rank 2 root system generated by $\alpha_b, \breve{\alpha}$. Draw all of the rank 2 root systems (which we have done in §5.3 on page 12). We see directly that the length of the root string is precisely the multiplicity of the edge, and see the direction of the arrow. The multiplicity of an edge in a Dynkin diagram can be $0, 1, 2$ or $3$, and $3$ occurs only for $G_2$. But the $G_2$-flag varieties are not cominuscule, by the classification of cominuscule varieties. □

9.6. **Symmetry of the box.** Soon we will finish proving the correctness of our algorithm to find the Dynkin diagram of the associated cominuscule subvariety of a flag variety. The reader can spot immediately in the pictures that there is a graph isomorphism taking the highest root (at the top of the box) to the lowest box root, and interchanging the labels on the Dynkin diagram. We will prove this. Our algorithm will take advantage of that symmetry of the box.

*Example* 19. Taking the $B_8$-flag variety

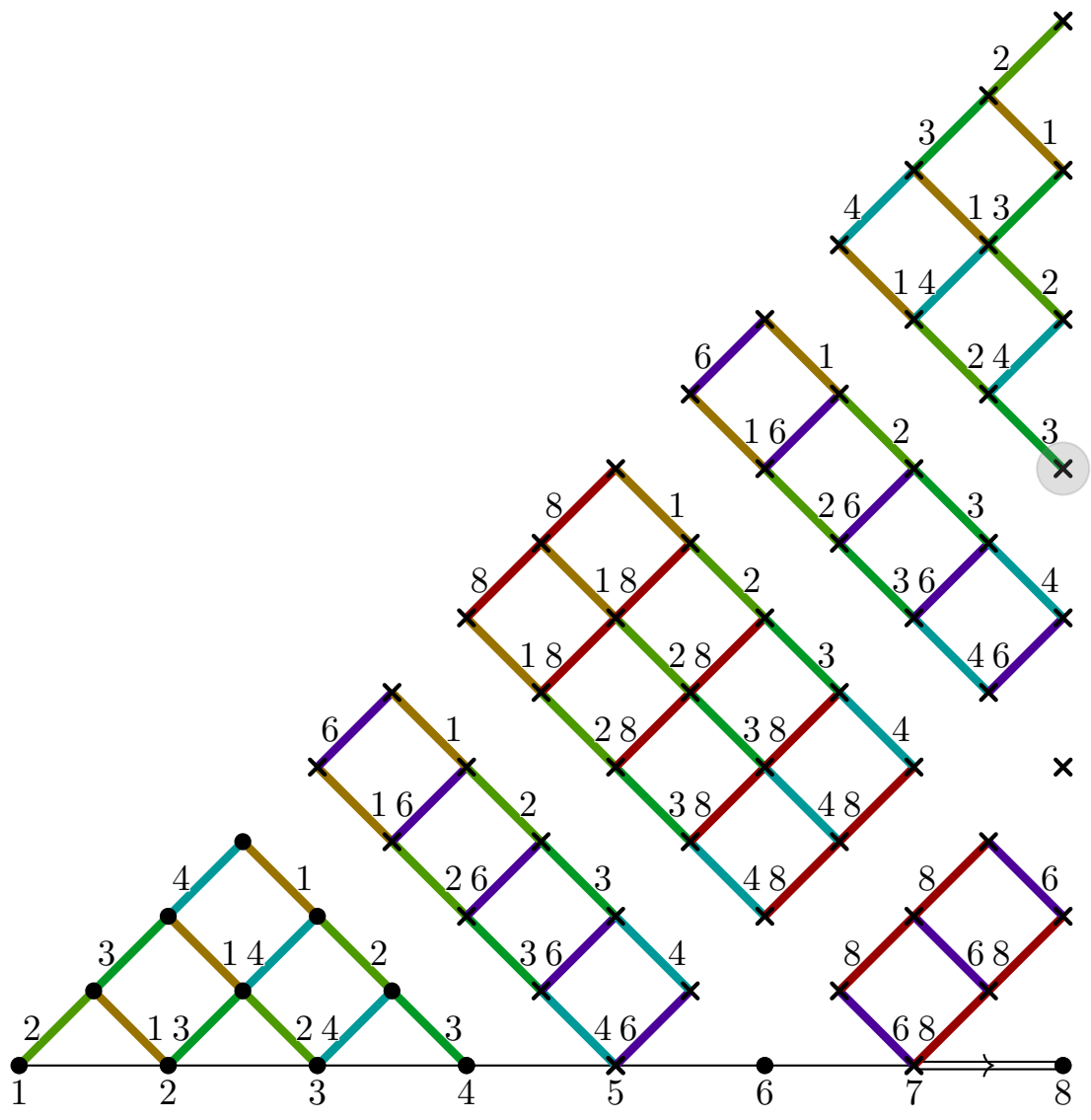

Turn the box upside down, interchanging labels according to the symmetry reversing the direction:

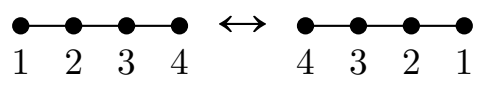

of the $\check{P}$-compact roots.

*Example* 20. For the $A_8$-flag variety ,

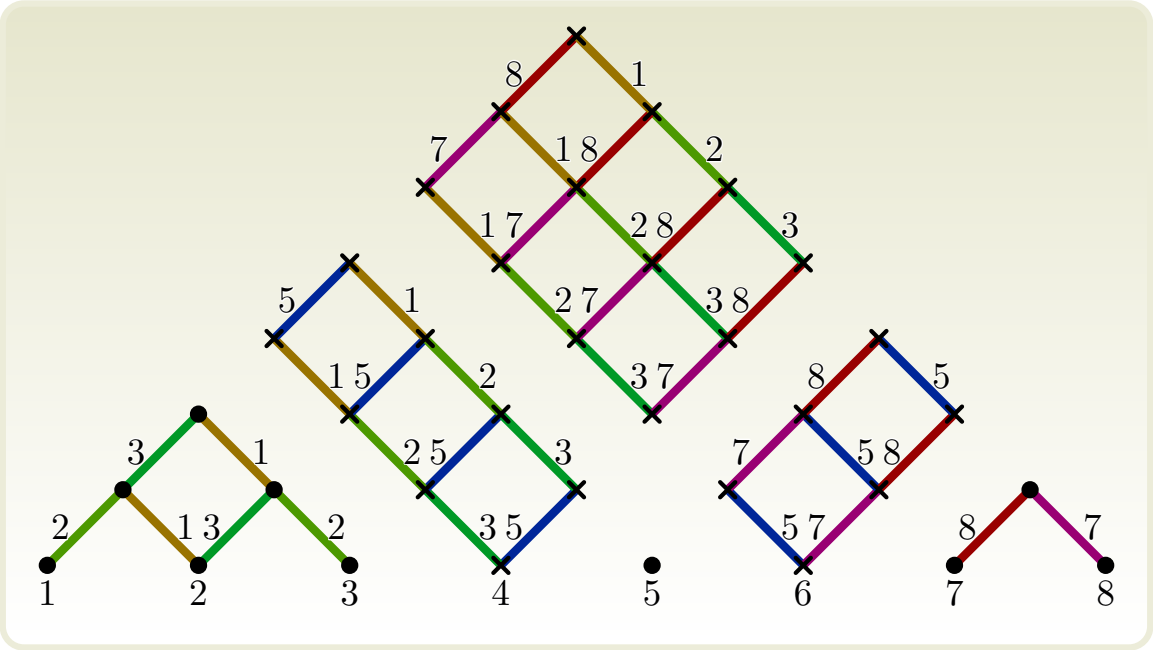

the symmetry of the box, turning it around 180°, comes from the left-right symmetry of the two triangles along the bottom. The associated cominuscule subvariety is the Grassmannian, given by removing roots $4, 5, 6$ and the components of the Hasse diagram that contain them.

To compute the Dynkin diagram of the associated cominuscule subvariety $(\check{X}, \check{G})$ of a flag variety $(X, G)$, we need to remove edges from the Hasse diagram of $G$. This requires computation of the Hasse diagram. We will use this symmetry to avoid calculating the Hasse diagram, working directly from the Dynkin diagram. The result will become obvious by eye without calculation.

9.7. **Outer automorphisms.** Recall that each of the Dynkin diagrams $A_r, D_r, E_6$ have graph automorphisms (Bourbaki [7] Chapter VIII, section 5):

- $A_r$ has automorphism group generated by permuting its two ends,

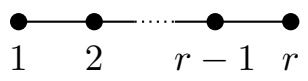

- $D_4$ has automorphism group the six permutations of its 3 ends,
- $D_{r \geqslant 5}$ has automorphism group generated by permuting $r-1, r$, when we order the roots of $D_r$ according to Bourbaki [8] pp. 265–290, plate IV:

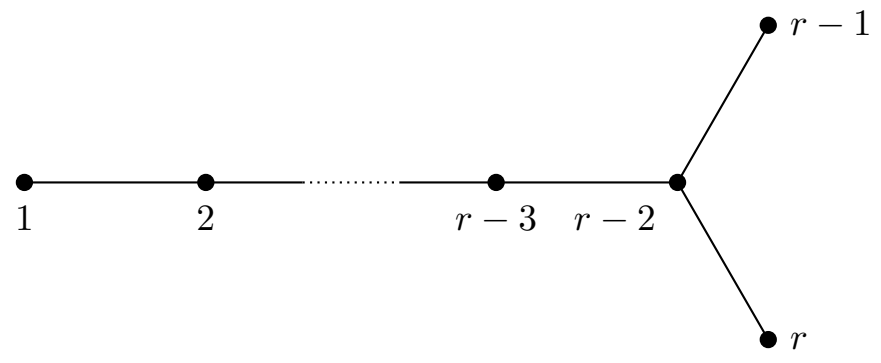

- $E_6$ has automorphism group generated by permuting nodes as

No other connected Dynkin diagram has any automorphism other than the identity [7] Chapter VIII, section 5.

Every group acts on itself by conjugation, its *inner automorphisms.* The automorphisms of the Dynkin diagram of a semisimple Lie group $G$ form the quotient of the automorphism group modulo inner automorphisms [7] Chapter VIII, section 5, [19] p. 498, proposition D.40. We can view the automorphisms as geometrically transforming the Dynkin diagram, or just as permutations of the labels of a labelled Dynkin diagram, leaving the nodes of the diagram in fixed locations on the paper. We will view them as relabelling.

From the Hasse diagram, as above, remove the edges corresponding to those crossed roots. Each connected component is clearly an irreducible module of the semisimple Levi factor.

**Lemma 15.** *Take the Hasse diagram of a cominuscule flag variety. The box of that Hasse diagram, as a graph, admits an involution that maps the highest root to the lowest root, and that maps edges to edges (while turning them upside down, i.e. reversing their arrows as directed edges), and that changes the labels of the edges according to some (possibly trivial) involution of the Levi Dynkin diagram. That automorphism, as applied to the Dynkin diagram of a reducible flag variety, is the product of the automorphisms of the components.*

*Proof.* We can assume that the cominuscule flag variety is irreducible. The longest Weyl group element $w_0, c$ of the compact root subsystem preserves the box, sending its highest weight to its lowest weight, sending each compact simple root to the negative of some (perhaps different) compact simple root [15] pp. 323–324. So it reverses the edges on the box.

By lemma 10 on page 27, we have an explicit description of every box, so we can look at each family of classical simple group, and at every exceptional simple group, just to see this reversal in action.

The reader is asked to look at all of the Hasse diagrams of the cominuscule flag varieties in table 8 on page 62.

$E_6$: We can see directly from the Dynkin diagram

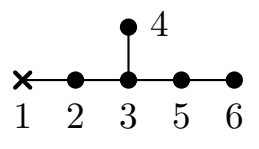

that, removing the crossed node at root 1, we have $D_5$,

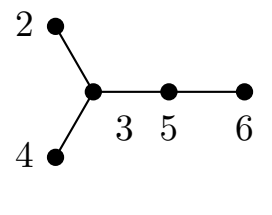

which admits a unique involution interchanging roots 2 and 4 of $E_6$, which also interchanges roots 3 and 6. Root 2 is the unique root going up from the lowest box root. Root 4 is the unique root going down from the highest root of the box:

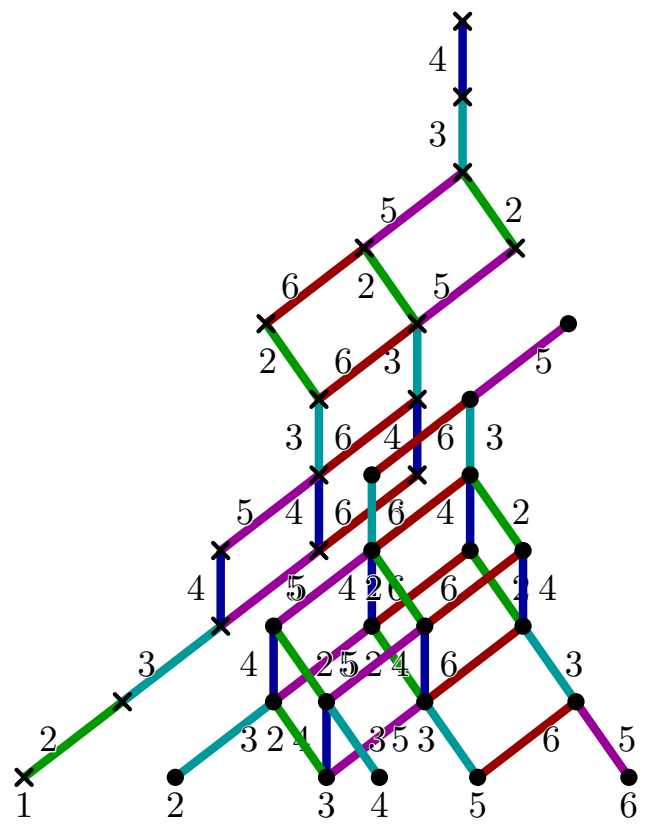

$E_7$: We can see directly from the Dynkin diagram

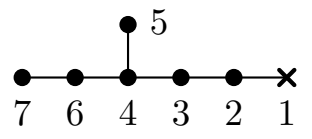

that, removing the crossed node at root 1, we have $E_6$:

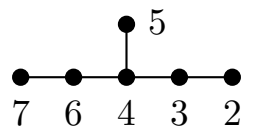

Clearly $E_6$ admits a unique involution interchanging roots 2 and 7 of $E_7$.

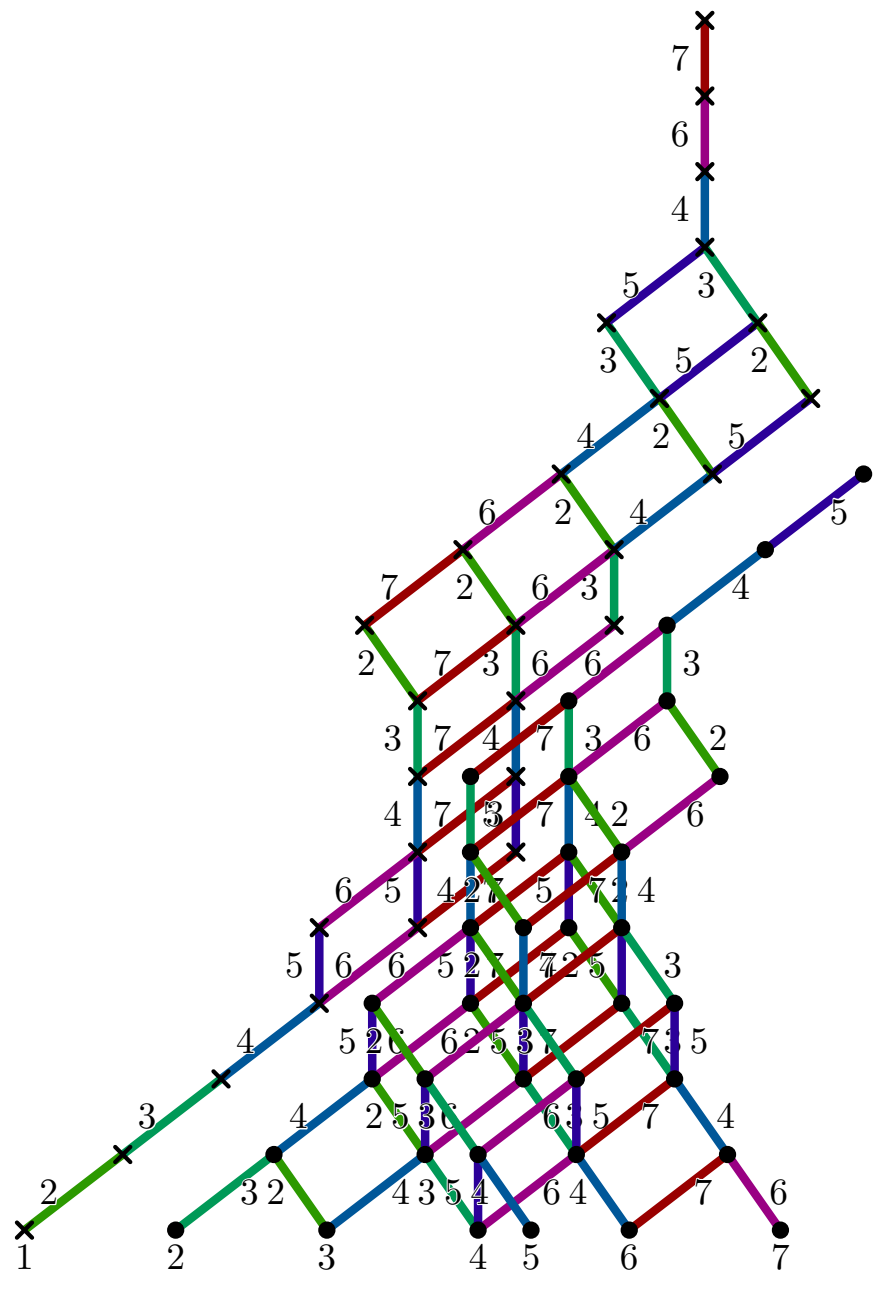

Root 2 is the unique root going up from the lowest box root. Root 7 is the unique root going down from the highest root of the box.

$B_r$: By lemma 10 on page 27, the box is a single path from lowest box root to highest root. Its labels are already symmetric under the graph isomorphism permuting highest and lowest box roots. We can apply the identity automorphism.

$D_r$: For

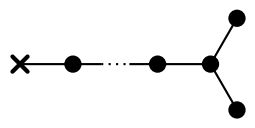

by lemma 10 on page 27, the box sits above 1 in the Hasse diagram ordering, a line segment with labels $2, 3, \dots, r-2$, then a square with labels $r-1, r$ on opposite sides, then a line segment with labels $r-2, r-3, \dots, 3, 2$. The automorphism permuting $r-1, r$ interchanges the sides of the square. But the box is again already symmetric under the graph isomorphism permuting the highest and lowest box roots. We can apply the identity automorphism.

For

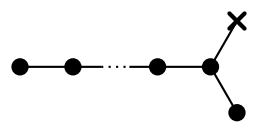

by lemma 10 on page 27, the box is a triangle, isomorphic as a labelled graph to a copy of $A_{r-1}$, and the unique nonidentity automorphism of $A_{r-1}$ interchanges the lowest and highest roots of the box.

For

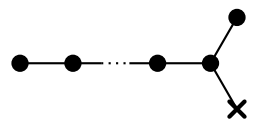

the obvious automorphism identifies it with

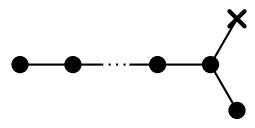

$A_r$: The $A_r$ cominuscule flag varieties: the Grassmannians. If the roots are ordered as above, following Bourbaki, each Grassmannian is given by crossing out one root, say root number $k$. Applying the obvious automorphism of the Dynkin diagram (*not* the Levi Dynkin diagram), we can assume that $k \leqslant r/2$.

Suppose first that $r = 1$, i.e. the projective line, and hence $k = 1$. The box is just that one root and the highest and lowest box roots in the box are the same.

Suppose next that $k = 1$ and $r > 1$, i.e. the Grassmannian is projective space $\mathbb{P}^r$. Then the Levi Dynkin diagram consists of roots $2, \dots, r$. Since root 1 is connected to root 2, and so on, the box contains roots

$$\alpha_1, \alpha_1 + \alpha_2, \dots, \alpha_1 + \dots + \alpha_r.$$

The first of these is the lowest box root of the box, and the last of these is the highest root, so also the highest root of the box. All other positive roots are compact. So this is the box. The edge 2 comes up from the lowest box root. The edge $r$ comes down from the highest root of the box. So the automorphism interchanges the ordering of the edges along the box "turning it upside down".

Suppose next that $k > 1$. The Levi Dynkin diagram splits into two components: roots $1, 2, \dots, k-1$ and roots $k+1, \dots, r$. The automorphisms of the Levi Dynkin diagram include reversing both of these components to $k-1, k-2, \dots, 2, 1$ and to $r, r-1, \dots, k+1$ respectively.

The box is given by starting at $k$, adding on roots with edges $k-1, k-2, \dots, 1$, and then adding to each of these successively with $k+1, k+2, \dots, r$, to make a rectangle of roots

$$\alpha_i + \dots + \alpha_j$$

for any $i \leqslant k \leqslant j$. The highest root is at $\alpha_1 + \cdots + \alpha_r$. The automorphism interchanges label $i$ with label $k-i$, so turning one side of that rectangle backwards, while it interchanges label $j$ with label $r+k+1-j$, turning the other side backwards. The roots move according to

$$\alpha_i + \cdots + \alpha_j \mapsto \alpha_{k-i} + \cdots + \alpha_{r+k+1-j}.$$

$C_r$: Recall [8] p. 254 that, if we let

$$\alpha_{ij} := \alpha_i + \cdots + \alpha_j$$

then row $g$ of the $C_r$-Hasse diagram contains roots

$$\begin{cases} \alpha_{i,i+g-1}, & i = 1, 2, \dots, r-g, \\ \alpha_{i,2r-g-i} + 2\alpha_{2r+1-g-i,r-1} + \alpha_r, & i = \max\left\{1, r+1-g\right\}, \dots, r - \left\lceil \frac{g}{2} \right\rceil, \\ 2\alpha_{r+(1-g)/2,r-1} + \alpha_r, & \text{if } g \text{ is odd}, \end{cases}$$

for rows parameterized by $g = 1, 2, \dots, 2r-1$. The crossed node is $r$. So the box consists of the roots with an $\alpha_r$-coefficient of 1: the roots

$$\begin{cases} \alpha_{i,2r-g-i} + 2\alpha_{2r+1-g-i,r-1} + \alpha_r, & i = \max\left\{1, r+1-g\right\}, \dots, r - \left\lceil \frac{g}{2} \right\rceil, \\ 2\alpha_{r+(1-g)/2,r-1} + \alpha_r, & \text{if } g \text{ is odd}, \end{cases}$$

The highest root is $2\alpha_{1,r-1} + \alpha_r$ [8] p. 254, on row $g = 2r-1$. The reader can check that the unique root on row $g = 2r-2$ is $\alpha_1 + 2\alpha_{2,r-1} + \alpha_r$, subtracting $\alpha_1$, and that there are two roots on row $g = 2r-3$:

$$\alpha_{1,2} + 2\alpha_{3,r-1} + \alpha_r, 2\alpha_{2,r-1} + \alpha_r,$$

subtracting $\alpha_2, \alpha_1$ respectively. Hence the Hasse diagram near the highest root is

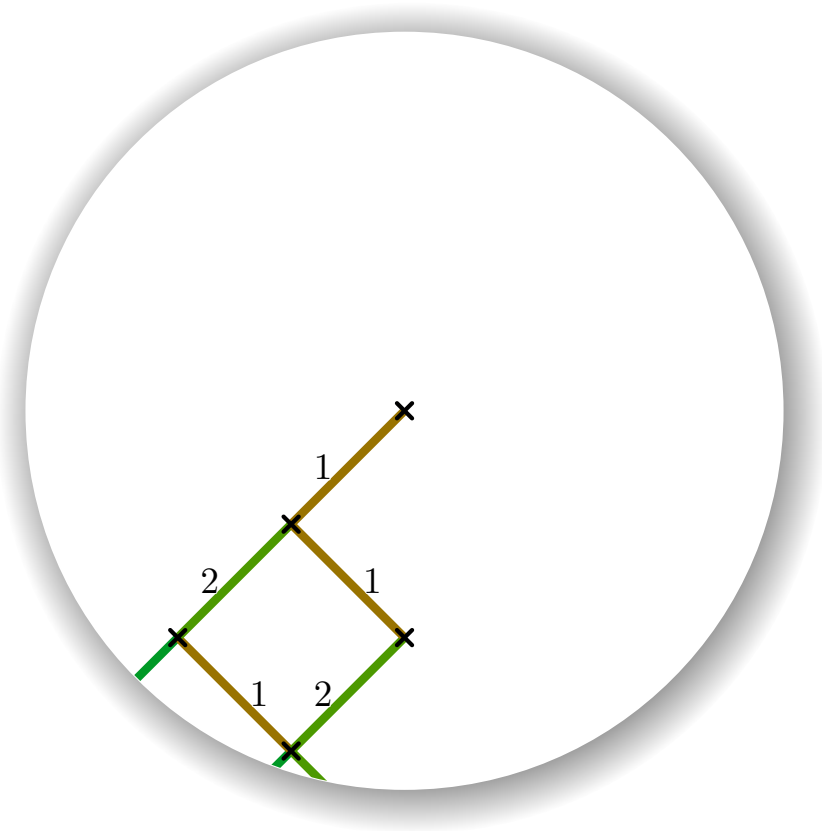

The only root in the Hasse diagram of $C_r$ in which 1-multiplicity is 2 is the highest root, since only the highest root has a coefficient of 2 in $\alpha_1$.

Similarly, the reader can check that the unique box root on row $g = 1$ is $\alpha_r$, the unique box root on row $g = 2$ is $\alpha_{r-1,r}$, and the unique box roots on row $g = 3$ are

$$\alpha_{r-2,r}, 2\alpha_{r-1} + \alpha_r$$

Hence the Hasse diagram near the lowest box root is

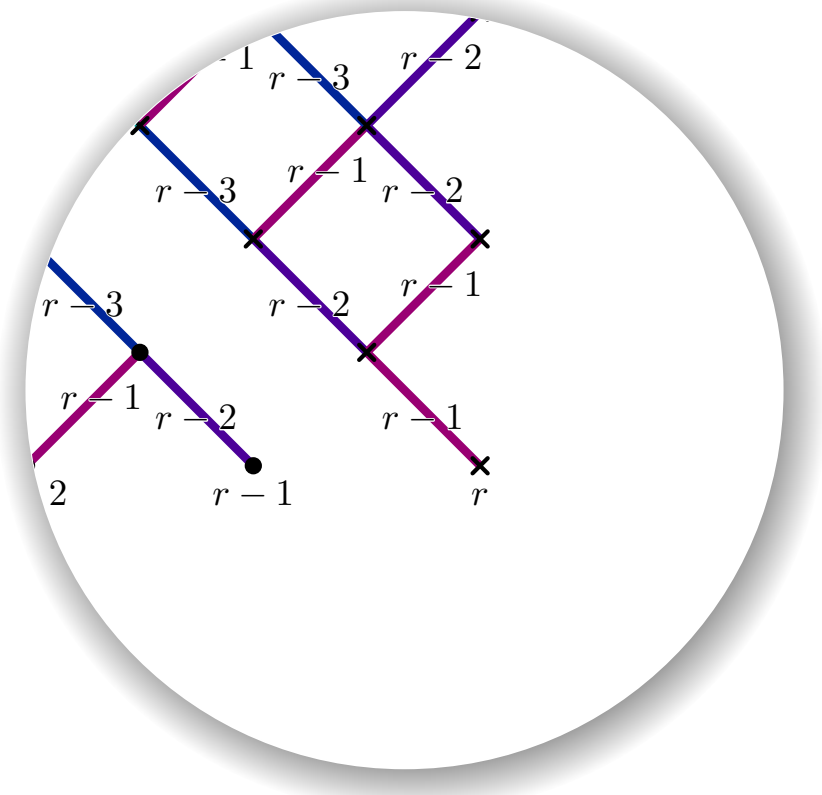

The only root in the Hasse diagram of $C_r$ in which $(r-1)$-multiplicity is 2 is the lowest box root $\alpha_r$, since $2\alpha_{r-1}+\beta$ is a positive root precisely for $\beta=\alpha_r$.

The Levi Dynkin diagram is given by cutting out the crossed node, i.e. is $A_{r-1}$. The nonidentity automorphism of $A_{r-1}$ interchanges labels by the permutation

$$\begin{pmatrix} 1 & 2 & 3 & \dots & r-1 & r \\ r & r-1 & r-2 & \dots & 2 & 1 \end{pmatrix}$$

Hence the labels near the highest root are interchanged with those near the lowest box root.

Perhaps we can see this more clearly starting at the lowest box root, following the path

$$\alpha_r \xrightarrow{r-1} \alpha_{r-1,r} \xrightarrow{r-1} 2\alpha_{r-1}+\alpha_r \xrightarrow{r-2} \alpha_{r-2}+2\alpha_{r-1}+\alpha_r \xrightarrow{r-2} 2\alpha_{r-2,r-1}+\alpha_r \dots$$
$$\dots \xrightarrow{i} \alpha_i+2\alpha_{i+1,r-1}+\alpha_r \xrightarrow{i} 2\alpha_{i,r-1}+\alpha_r \xrightarrow{i-1} \dots$$
$$\dots \xrightarrow{1} \alpha_1+2\alpha_{1,r-1}+\alpha_r \xrightarrow{1} 2\alpha_{1,r-1}+\alpha_r,$$

where the labels of edges between these roots are

$$r-1, r-1, r-2, r-2, \dots, i, i, i-1, i-1, \dots, 2, 2, 1, 1.$$

The automorphism of the Dynkin diagram swaps the labels, and we can see that the labels are put down in reverse order.

□

**Lemma 16.** *Take the Hasse diagram of a flag variety. The box of that Hasse diagram, as a graph, admits an involution that maps the highest root to the lowest root, and that maps edges to edges (while turning them upside down, i.e. reversing their arrows as directed edges), and that changes the labels of the edges according to some (possibly trivial) involution of the Levi Dynkin diagram.*

*Proof.* Its root system contains the root system of its associated cominuscule subvariety, and the two have the same box by lemma 6 on page 18, so the result follows from the previous lemma. □

*Example* 21. Consider the Dynkin diagram of the associated cominuscule subvariety of an irreducible flag variety. As above, we still have not found out how to attach the crossed node to the Levi Dynkin diagram, if we just look at the Dynkin diagram of the original flag variety. Suppose that we somehow found that we needed to attach a double edge. By the symmetry in this lemma, we would have a root string of length 2 coming down from the highest root. Among all Hasse diagrams of simple Lie groups, only a Hasse diagram of type $C$ has a root string of length 2 coming down from its highest root. So the only irreducible flag varieties $(X, G)$ that have a $C$-series associated cominuscule subvariety are $C$-series flag varieties.

9.8. **Using the symmetry.** We need to see how to compute, for every Dynkin diagram, the length of the root string of each simple root through the highest root. This, as above, is the multiplicity of edge in the associated cominuscule subvariety. We can forget flag varieties entirely: these root string lengths are well defined integer invariants of any root system.

Take an irreducible Dynkin diagram $\Delta$. Over each node $b$, say with associated simple root $\alpha_b$, write the length of the $\alpha_b$ root string through the highest root. (We don't write anything if this is zero.) By the symmetry, this will provide the edge multiplicities of all associated cominuscule subvarieties $(\check{X}, \check{G})$ of all flag varieties $(X, G)$ for which $G$ has Dynkin diagram $\Delta$. So these multiplicities depend on $\Delta$, the Dynkin diagram of the simple Lie group $G$, not on the particular choice of $G$-flag variety $X$, although which of them we need to compute depends on $X$. So it is the same for $(X, G)$ as for any cominuscule $G$-variety, i.e. for any cominuscule flag variety whose group has Dynkin diagram $\Delta$, with some root crossed. In particular, $\Delta = A_r, B_r, C_r, D_r, E_6, E_7$ each occur as the Dynkin diagram of a cominuscule flag variety. Since we know the Dynkin diagrams of those cominuscule flag varieties, we get explicitly the multiplicities of the root strings. This leaves only $F_4$, $G_2$ and $E_8$ to compute by hand, directly from their Hasse diagrams.

The multiplicities in the extended Dynkin diagram of a simple Lie group $G$ are defined as the Cartan integers

$$\alpha_0 \cdot \check{\alpha}_b$$

where $\alpha_0 = -\theta$ is the negative of the highest root $\theta$ [8] ch. VI, §4, p. 198. The length of the $\alpha_b$-root string through $\theta$ is

$$|\alpha_0 \cdot \check{\alpha}_b| = \theta \cdot \check{\alpha}_b$$

by [24] p. 145 Proposition. So the multiplicities in the extended Dynkin diagram are the edge multiplicities of the edges reaching the highest root in the Hasse diagram. It is consistent with the definition of the affine root, since the extended Dynkin diagram has as its affine root the negative of the highest root. Summing up, we have proved the main theorem: theorem 1 on page 3.

## 10. Freedom

In this section, we prove that the associated cominuscule subvariety in any irreducible flag variety is preserved by a subgroup of the automorphisms of the flag variety. We prove that this subgroup is as large as possible subject to an open condition on tangent spaces. This is similar to the theory of the twisted cubic: nondegeneracy of a cubic space curve implies homogeneity, but for the associated cominuscule subvariety, we don't need to constrain its degree [21] p. 10–14.

10.1. **Defining freedom.** Suppose that $X = G/P$ is an irreducible flag variety. A linear subspace $V \subseteq T_x X$ in a tangent space is *subfree* if it does not contain any line lying in a $G$-invariant coherent proper subsheaf of $TX$. It is *free* if it is maximally subfree. (Note that every $G$-invariant coherent subsheaf of $TX$ is a homogeneous holomorphic vector bundle, hence a holomorphic distribution. Also, $TX$ contains at most a finite collection of $G$-invariant coherent subsheaves, as each is determined by a $P$-submodule of $\mathfrak{g}/\mathfrak{p}$.) If $X$ is cominuscule, there is no such invariant subsheaf; by maximality, $V = T_x X$. If $X$ is not cominuscule, $V$ is complementary to the sum of the proper $G$-invariant homogeneous subbundles. So the set of free linear subspaces at a point $x \in X$ is an affine open cell of a Grassmannian, so Zariski open. The tangent spaces of the associated cominuscule subvariety are free. Hence every free linear subspace $V$ has dimension equal to that of the associated cominuscule subvariety.

10.2. **Free morphisms.** Take a morphism $Z \xrightarrow{\varphi} X$ from a reduced complex space $Z$ to an irreducible flag variety $X$. A point $z_0 \in Z$ is *free* if $z_0$ is a smooth point of $Z$ at which $\varphi'(z_0)$ is injective with free image. Intuitively, think of each $G$-invariant distribution on $X$ as an invariant differential equation of first order. Near a free point of $Z$, any curve smooth and tangent to $Z$ at that point doesn't satisfy any such equation. We say $Z$ is a *free complex space* if every point is free, *generically free* if its open set of free points $Z^{\text{free}} \subseteq Z$ is dense in every component. Note that $Z^{\text{free}}$, if not empty, is a free immersed complex submanifold.

10.3. **Automorphisms.** An *automorphism* of a morphism $Z \xrightarrow{\varphi} X$ from a complex space $Z$ to a flag variety $(X, G)$ is a pair $(f, g)$ of automorphism $Z \xrightarrow{f} Z$ of complex spaces and $g \in G$ so that $\varphi(f(z)) = g\varphi(z)$ for all $z \in Z$. The automorphism group of $Z \xrightarrow{\varphi} X$ we denote $G_Z$ and adorn with the topology of uniform convergence on compact sets. If $Z$ has finitely many components, so does its smooth locus. The automorphisms of the smooth part of $Z$ then form a Lie group acting smoothly and properly on a sufficiently high jet bundle of maps $Z \to X$ [42] p. 62 theorem 25. The automorphisms of any generically free complex space with finitely many components form a closed subgroup of that Lie group, hence also a Lie group acting smoothly on each smooth stratum. The associated cominuscule subvariety is a natural object to study because:

**Theorem 4.** *Take a connected generically free complex manifold $Z \xrightarrow{\varphi} X$ in an irreducible flag variety $X$. The automorphism group of $Z$ has dimension*

$$\dim G_Z \leqslant \dim G_{\check{X}}.$$

*The following are equivalent:*

- $\dim G_Z = \dim G_{\check{X}}$.
- *$G_Z$ is conjugate to $G_{\check{X}}$ in $G$.*
- *$Z = g\check{X}$ for some $g \in G$, i.e. $Z$ is an associated cominuscule subvariety.*

**Theorem 5.** *If a connected compact complex manifold $Z$ with holomorphic map $Z \to X$ to a flag variety $(X, G)$ is free then it is homogeneous under automorphisms, and is a product of an abelian variety and a flag variety.*

Before we can prove theorem 4 and theorem 5, we need some additional notation.

10.4. **Root systems.** We use notation and terminology from [32]. Pick a point $x_0 \in X$ in a flag variety $(X, G)$ and let $P := G^{x_0}$. Then $P$ contains a Cartan subgroup $H \subset P \subset G$, unique up to $P$-conjugacy. Recall that a Cartan subgroup is a maximal connected abelian subgroup of $G$ for which $\mathfrak{g}$ splits

$$\mathfrak{g} = \mathfrak{h} \oplus \bigoplus_{\alpha \in \Delta} \mathfrak{g}_\alpha$$

into a sum of the Cartan subalgebra with various 1-dimensional holomorphic $H$-modules, the *root spaces*, each thus isomorphic to the $H$-module of a character $\alpha \in \check{H} \subset \mathfrak{h}^\vee$ of $H$. These characters $\alpha$ are the *roots*. Denote by $\Delta$ the set of roots.

A root $\alpha$ is *$P$-compact* if the root spaces of $\alpha$ and $-\alpha$ belong to $\mathfrak{p}$. The Cartan subalgebra, together with the root spaces of the $P$-compact roots, generate a subalgebra $\mathfrak{g}_0 \subset \mathfrak{p}$. There is a unique element $e \in \mathfrak{h}$, the *P-grading element*, in the centre of $\mathfrak{g}_0$ so that

$$\alpha_i(e) = \begin{cases} 0, & \alpha_i \ P\text{-compact}, \\ 1, & \alpha_i \ P\text{-noncompact}. \end{cases}.$$

[15] p. 239 proposition 3.1.2. We grade weights by inner product with $e$, and preorder weights by grade. Note that the $P$-positive roots are *not* necessarily the positive roots in the sense of the Cartan subgroup $H$ and a choice of base of simple roots; indeed they are the positive roots in that sense just when $P$ is a Borel subgroup.

10.5. **Chevalley bases.** Pick a complex semisimple Lie group $G$ and a Cartan subgroup $H \subset G$. Denote the Killing form on $\mathfrak{h}^\vee$ by $\alpha, \beta \mapsto \alpha \cdot \beta$. A *Chevalley basis* $e_\alpha \in \mathfrak{g}_\alpha, \check{\alpha} \in \mathfrak{h}$ is a spanning set of $\mathfrak{g}$, parameterized by roots $\alpha \in \Delta$, so that

(1) $[\check{\beta} e_\alpha] = 2\frac{\alpha\cdot\beta}{\beta^2} e_\alpha$
(2) $[\check{\alpha}\check{\beta}] = 0$,
(3)
$$[e_\alpha, e_\beta] = \begin{cases} \check{\alpha}, & \text{if } \alpha+\beta=0, \\ N_{\alpha\beta} e_{\alpha+\beta}, & \text{if } \alpha+\beta \text{ a root}, \\ 0, & \text{otherwise} \end{cases}$$
with
(a) $N_{\alpha\beta} = \pm(p+1)$, where $p$ is the largest integer for which $\beta - p\,\alpha$ is a root,
(b) $N_{-\alpha,-\beta} = -N_{\alpha\beta}$,
(c) we set $N_{\alpha\beta} = 0$ if $\alpha+\beta = 0$ or $\alpha+\beta$ is not a root.

It follows then that

$$\alpha(\check{\beta}) = 2\frac{\alpha\cdot\beta}{\beta^2}.$$

Every complex semisimple Lie group $G$ with a Cartan subgroup $H \subset G$ admits a Chevalley basis [53] p. 51. A Chevalley basis in this sense is really only a spanning set, not a basis; if we pick out a basis of simple roots $\alpha_i$ then $\{\check{\alpha}_i\}_{i=1}^r$ is a basis of $\mathfrak{h}$ and $\{\check{\alpha}_i\}_{i=1}^r \sqcup \{e_\alpha\}_{\alpha\in\Delta}$ is a basis of $\mathfrak{g}$. The *Cartan integers* are

$$n_{ji} := 2\frac{\alpha_i\cdot\alpha_j}{\alpha_j^2},$$

and the Cartan matrix is the integer matrix of the transposed entries $C = (n_{ij})$, with inverse matrix $C^{-1} = (C^{ij})$. The *fundamental weights* are

$$\varpi^i := C^{ji}\alpha_j.$$

10.6. **Dual bases.** The dual basis is defined by

$$\begin{aligned} \alpha(\check{\beta}) &= 2\frac{\alpha\cdot\beta}{\beta^2}, \\ \alpha(e_\beta) &= 0, \\ \omega^\alpha(\check{\beta}) &= 0, \\ \omega^\alpha(e_\beta) &= \delta^\alpha_\beta. \end{aligned}$$

The dual basis vectors extend uniquely to left invariant 1-forms on $G$, to which we apply

$$d\omega(X,Y) = \mathcal{L}_X(\omega(Y)) - \mathcal{L}_Y(\omega(X)) - \omega([X,Y]),$$

to left invariant 1-forms $\omega$ and vector fields $X, Y$. The first two terms vanish by left invariance:

$$d\omega(X,Y) = -\omega([X,Y]).$$

Let $\omega^{\beta\gamma} := \omega^{\beta} \wedge \omega^{\gamma}$, etc. and compute the structure equations of semisimple Lie groups:

$$
\begin{aligned}
d\omega^{\alpha} &= -\alpha \wedge \omega^{\alpha} - \frac{1}{2} \sum_{\beta+\gamma=\alpha} N_{\beta\gamma}\omega^{\beta\gamma} \\
d\alpha &= -\sum_{\beta} \frac{\alpha \cdot \beta}{\beta^2} \omega^{\beta,-\beta}
\end{aligned}
$$

with sums over all roots.

10.7. **Proving the freedom theorems.** We now prove theorem 4 on page 41:

*Proof.* We carry out the moving frame method on an arbitrary free immersed connected submanifold [20, 28]. Note that freedom is precisely that every tangent space is complementary to the largest proper $G$-invariant subbundle inside $TX$. We will see that various differential invariants vanish precisely for the associated cominuscule subvariety.

Pick a point $x_0 \in X$ so that we can identify $X = G/P$. Let $Z'$ be the set of pairs $(z, g)$ for $z \in Z$ and $g \in G$ so that $\varphi(z) = gx_0$. Each automorphism $(f, a)$ acts on $Z'$ by

$$
(z, g) \mapsto (f(z), ag).
$$

The map $Z' \to Z$ is a holomorphic principal right $P$-bundle: the pullback bundle

$$
\begin{array}{ccc}
Z' & \longrightarrow & G \\
\downarrow & & \downarrow \\
Z & \longrightarrow & X.
\end{array}
$$

We write $P$-maximal roots as $\alpha^{+}$, $P$-minimal as $\alpha^{-}$, and $P$-compact as $\alpha^{\circ}$. On $G$, the forms $\omega^{\beta}$ for $\beta$ $P$-negative are a basis for the semibasic 1-forms for $g \in G \mapsto gx_0 \in X$. The equations

$$
0 = \omega^{\beta^-} \text{ for all } \beta^-
$$

cut out the preimage in $TG$ of the largest $G$-invariant subbundle $V \subset TX$. On $Z'$, the Maurer–Cartan form restricts to have $\omega^{\beta^-}$ a basis for the semibasic forms, since the tangent space to $Z$ is complementary to $V$. So on $Z' \subset Z \times G$,

$$
\omega^{\alpha} = a^{\alpha}_{\beta^-}\omega^{\beta^-},
$$

for each $P$-negative root $\alpha$ which is not $P$-minimal, with the sum over $P$-minimal roots $\beta^-$. So on $Z'$ we have equations for all $\omega^{\alpha}$ for all $\alpha$ $P$-negative and not $P$-minimal.

Take exterior derivative:

$$
\begin{aligned}
0 &= d\omega^{\alpha} - da^{\alpha}_{\beta^-} \wedge \omega^{\beta^-} - a^{\alpha}_{\beta^-} d\omega^{\beta^-}, \\
&= -\alpha \wedge a^{\alpha}_{\beta^-}\omega^{\beta^-} \\
&\quad - N_{\alpha-\beta^-,\beta^-}\omega^{\alpha-\beta^-,\beta^-} \\
&\quad - \sum_{\gamma+\varepsilon=\alpha}{}' N_{\gamma\varepsilon}\omega^{\gamma\varepsilon} - \frac{1}{2}\sum_{\gamma+\varepsilon=\alpha}{}'' N_{\gamma\varepsilon}\omega^{\gamma\varepsilon} \\
&\quad - da^{\alpha}_{\beta^-} \wedge \omega^{\beta^-} \\
&\quad + a^{\alpha}_{\beta^-}\beta^- \wedge \omega^{\beta^-} \\
&\quad + a^{\alpha}_{\beta^-}\left( \sum_{\gamma^{\circ}+\varepsilon^{-}=\beta^{-}} N_{\gamma^{\circ}\varepsilon^{-}}\omega^{\gamma^{\circ}\varepsilon^{-}} + \frac{1}{2}\sum_{\gamma+\varepsilon=\beta^-}{}'' N_{\gamma\varepsilon}\omega^{\gamma\varepsilon} \right),
\end{aligned}
$$

where $\sum'$ is the sum over a nonmaximal positive or compact and nonminimal negative root and $\sum''$ is the sum over two nonminimal negative roots, so

$$
\begin{aligned}
0 &= -a^{\alpha}_{\beta^-}\alpha \wedge \omega^{\beta^-} \\
&\quad - N_{\alpha-\beta^-,\beta^-}\omega^{\alpha-\beta^-,\beta^-} \\
&\quad - \frac{1}{2}\sum_{\gamma+\beta=\alpha}^{\gamma,\beta<0} N_{\gamma\beta} a^{\gamma}_{\varepsilon^-} a^{\beta}_{\sigma^-}\omega^{\varepsilon^-\sigma^-} \\
&\quad - \sum_{\gamma+\beta=\alpha}^{\gamma\geqslant 0,\beta<0} N_{\gamma\beta} a^{\beta}_{\sigma^-}\omega^{\gamma\sigma^-} \\
&\quad - da^{\alpha}_{\beta^-} \wedge \omega^{\beta^-} \\
&\quad + a^{\alpha}_{\beta^-}\beta^- \wedge \omega^{\beta^-} \\
&\quad + a^{\alpha}_{\beta^-}\sum_{\gamma^\circ+\varepsilon^-=\beta^-} N_{\gamma^\circ\varepsilon^-}\omega^{\gamma^\circ\varepsilon^-} \\
&\quad + \frac{1}{2}a^{\alpha}_{\beta^-}\sum_{\gamma+\varepsilon=\beta^-} N_{\gamma\varepsilon} a^{\gamma}_{\sigma^-} a^{\varepsilon}_{\tau^-}\omega^{\sigma^-\tau^-} \\
&= -\left(da^{\alpha}_{\beta^-} + (\alpha-\beta^-)\wedge a^{\alpha}_{\beta^-} + N_{\alpha-\beta^-,\beta^-}\omega^{\alpha-\beta^-}\right)\wedge\omega^{\beta^-} \\
&\quad - \frac{1}{2}\sum_{\gamma+\sigma=\alpha}^{\gamma,\sigma<0} N_{\gamma\sigma} a^{\gamma}_{\varepsilon^-} a^{\sigma}_{\beta^-}\omega^{\varepsilon^-\beta^-} \\
&\quad - \sum_{\gamma+\sigma=\alpha}^{\gamma\geqslant 0,\sigma<0} N_{\gamma\sigma} a^{\sigma}_{\beta^-}\omega^{\gamma\beta^-} \\
&\quad + a^{\alpha}_{\varepsilon^-}\sum_{\gamma^\circ+\beta^-=\varepsilon^-} N_{\gamma^\circ\beta^-}\omega^{\gamma^\circ\beta^-}, \\
&\quad + \frac{1}{2}a^{\alpha}_{\tau^-}\sum_{\gamma+\varepsilon=\tau^-} N_{\gamma\varepsilon} a^{\gamma}_{\sigma^-} a^{\varepsilon}_{\beta^-}\omega^{\sigma^-\beta^-}
\end{aligned}
$$

So we let

$$
\begin{aligned}
\nabla a^{\alpha}_{\beta^-} &:= da^{\alpha}_{\beta^-} + a^{\alpha}_{\beta^-}(\alpha-\beta^-) + N_{\alpha-\beta^-,\beta^-}\omega^{\alpha-\beta^-} + \sum_{\gamma+\sigma=\alpha}^{\gamma\geqslant 0,\sigma<0} N_{\gamma\sigma} a^{\sigma}_{\beta^-}\omega^{\gamma} \\
&\quad - a^{\alpha}_{\varepsilon^-}\sum_{\gamma^\circ+\beta^-=\varepsilon^-} N_{\gamma^\circ\beta^-}\omega^{\gamma^\circ} - \frac{1}{2}a^{\alpha}_{\tau^-}\sum_{\gamma+\varepsilon=\tau^-} N_{\gamma\varepsilon} a^{\gamma}_{\sigma^-} a^{\varepsilon}_{\beta^-}\omega^{\sigma^-} + \frac{1}{2}\sum_{\gamma+\sigma=\alpha} N_{\gamma\sigma} a^{\gamma}_{\varepsilon^-} a^{\sigma}_{\beta^-}\omega^{\varepsilon^-}
\end{aligned}
$$

and we find that

$$0 = \nabla a^{\alpha}_{\beta^-} \wedge \omega^{\beta^-}$$

Hence

$$\nabla a^{\alpha}_{\beta^-} = a^{\alpha}_{\beta^-\gamma^-}\omega^{\gamma^-},$$

for unique functions $a^{\alpha}_{\beta^-\gamma^-} = a^{\alpha}_{\gamma^-\beta^-}$ on $Z'$.

Take a $G$-root $\gamma$ which is not a $G'$-root. If $\gamma$ is $P$-negative then it is negative, so is a $P$-nonminimal negative root, hence $\omega^{\gamma} = a^{\gamma}_{\beta^-}\omega^{\beta^-}$ is solved for in terms of $G'$-roots $\beta^-$. Suppose that $\gamma$ is $P$-null, i.e. $P$-compact, so a root of the maximal reductive Levi factor $G_0 \subseteq P$. Then $\gamma$ is a $G'$-root, a contradiction. So we can suppose that $\gamma$ is $P$-positive, not $P$-maximal. Since $P'$ contains all roots perpendicular to the $P$-maximal roots, we can assume that $\gamma$ is not perpendicular to some $P$-maximal root, which we write as $-\beta^-$ for some $P$-minimal root $\beta^-$.

The roots $\beta^-, \gamma$ form a basis for a rank 2 root system; looking over all of the possibilities above, since $\gamma$ is $P$-positive and $\beta^-$ is $P$-negative, the angle between them is more than a right angle. So the sum $\alpha := \gamma + \beta^-$ is also a $G$-root [53] p. 29. Conversely, by reversing the same steps, if $\alpha$ is a $P$-nonminimal root which can be written as $\alpha = \gamma + \beta^-$ with $N_{\beta^-\gamma} \neq 0$, i.e. $\beta^-, \gamma$ not perpendicular, then $\gamma$ is a $P$-root which is not a $P'$-root.

Take any such $\alpha, \beta^-, \gamma$. Going back to the equations of $\nabla a^\alpha_{\beta^-}$, move in the direction of the root vector of $\gamma$. This root vector is a complete vector field as it lies in the fiber of $Z' \to Z$, so it is given by the action of a 1-parameter subgroup of $P$. We find $da^\alpha_{\beta^-} = -N_{\gamma\beta}$, a constant. So over every point of $Z$, we can find a point of $Z'$ at which $a^\alpha_{\beta^-} = 0$. Just as in the model, we can arrange that $a^\alpha_{\beta^-} = 0$, for every root $\gamma$ which is a $G$-root and not a $G'$-root. We reduce the structure group down from $P$ to have $\omega^\gamma$ semibasic. At each step when we do this, we lose one dimension, cutting $Z'$ down to a hypersurface, cut out by the equation $a^\alpha_{\beta^-} = 0$. The differentials of these equations are linearly independent, so these hypersurfaces intersect transversely.

We can see this from another perspective: if we let $p \in P$ be the element $p := \exp(te_\gamma)$, for some $t \in \mathbb{C}$, we can compute from the right action on the Maurer–Cartan form $r_p^*\omega = \mathrm{Ad}_p^{-1}\omega$, so that

$$
\begin{aligned}
0 &= r_p^* 0, \\
&= r_p^*(\omega^\alpha - a^\alpha_{\beta^-}\omega^{\beta^-}), \\
&= \omega^\alpha - (tN_{\gamma\beta^-} + r_p^* a^\alpha_{\beta^-})\omega^{\beta^-}, \\
&= a^\alpha_{\beta^-}\omega^{\beta^-} - (tN_{\gamma\beta^-} + r_p^* a^\alpha_{\beta^-})\omega^{\beta^-}, \\
&= (a^\alpha_{\beta^-} - tN_{\gamma\beta^-} - r_p^* a^\alpha_{\beta^-})\omega^{\beta^-},
\end{aligned}
$$

so that

$$
r_p^* a^\alpha_{\beta^-} = a^\alpha_{\beta^-} - tN_{\gamma\beta^-},
$$

which shows that we can set

$$
t := \frac{a^\alpha_{\beta^-}}{N_{\gamma\beta^-}},
$$

above each point of $Z$, to find a point where $a^\alpha_{\beta^-} = 0$. Above each chosen point of $Z$, the linear relations among the various $\omega^\alpha$ transform under $P$-action to other linear relations. One easily sees that these equations cut out a smooth complex submanifold $Z'' \subseteq Z'$ which is a fiber bundle over $Z$.

We see the bound on the dimension of automorphism group: counting equations, we see that $\dim Z'' = \dim G'$. Since automorphisms of $Z$ act freely on $Z'$ preserving $Z''$, they act freely on $Z''$, so the automorphism group $G_Z$ has orbits which are copies of itself, hence $G_Z$ has dimension at most that of $G'$.

We arrange $a^\alpha_{\beta^-} = 0$ on $Z''$ for one $a^\alpha_{\beta^-}$ for every $\gamma = \alpha - \beta^-$ which is $P$-root and not a $P'$-root. So we have $Z'' \subset Z'$ a smooth complex subvariety of complex codimension equal to that number.

Start this process with the smallest possible root $\alpha$ (in some root ordering), and proceed through all such, making some choice of $\alpha = \gamma + \beta^-$ as above. At each step, modulo the various $\omega^{\beta^-}$, we find

$$
0 = N_{\gamma\beta^-}\omega^\gamma + \sum_{\varepsilon+\sigma=\alpha}^{\varepsilon\geqslant 0, \sigma<0} N_{\varepsilon\sigma} a^\sigma_{\beta^-}\omega^\varepsilon - a^\alpha_{\varepsilon^-} \sum_{\tau^\circ+\beta^-=\varepsilon^-} N_{\tau^\circ\beta^-}\omega^{\tau^\circ}.
$$

But by induction, the second term is a linear combination of $P$-compact and $P$-minimal roots. We conclude that, for each $P$-root $\gamma$ which is not a $P'$-root, on $Z''$, $\omega^\gamma$ is a linear combination of $\omega^{\sigma^o}, \omega^{\beta^-}$ for $P$-compact root $\omega^{\sigma^o}$ and $P$-minimal roots $\omega^\beta$.

Suppose we find some $a^\alpha_{\beta^-} \ne 0$ on $Z''$. Let $\gamma := \alpha - \beta^-$; so $\gamma$ is a $P$-root, not a $P'$-root. By $P'$-action, looking again at $da^\alpha_{\beta^-}$, we can arrange that $a^\alpha_{\beta^-} = 1$, and on this smooth fiber subbundle, we have $\gamma$ now a linear combination of $\check{G}$-roots. Putting in

$$d\gamma = -\sum_\sigma \frac{\gamma \cdot \sigma}{\sigma^2} \omega^{\sigma,-\sigma}$$

apply Cartan's lemma: for each $P$-maximal root $\sigma^+$ not perpendicular to $\gamma$, $\omega^{\sigma^+}$ is a linear combination of $P$-compact and $P$-minimal roots. Since there is at least one such $\sigma^+$, $\dim G_Z \leqslant \dim G_{\check{X}} - 2$.

On the other hand, suppose that all $a^\alpha_{\beta^-}$ vanish at every point of $Z''$, i.e. $Z$ is tangent to an associated cominuscule subvariety at each point. So $\omega^\alpha = 0$ for all $\alpha$ $P$-negative and not $P$-minimal. Repeating the calculation above,

$$N_{\alpha-\beta^-,\beta^-} \omega^{\alpha-\beta^-} \wedge \omega^{\beta^-} = 0,$$

so if $\gamma$ is any $P$-root which is not a $P'$-root, then we can write it as $\gamma = \alpha - \beta^-$ and find $\omega^\gamma = a^\gamma \omega^{\beta^-(\gamma)}$, for some functions $a^\gamma$, where here $\beta^- = \beta^-(\gamma)$ is some particular $P$-minimal root associated to each $\gamma$. So the complex manifold $Z''$ is of dimension at most that of $G_{\check{X}}$.

Differentiate the equation $\omega^\gamma = a^\gamma \omega^{\beta^-(\gamma)}$:

$$0 = -(da^\gamma + a^\gamma(\gamma - \beta^-(\gamma))) \wedge \omega^{\beta^-(\gamma)} - \frac{1}{2} N_{\mu+\nu=\gamma} N_{\mu\nu} \omega^{\mu\nu} + \frac{1}{2} a^\gamma N_{\mu+\nu=\beta^-(\gamma)} N_{\mu\nu} \omega^{\mu\nu}.$$

Even when we plug in the equations $\omega^\alpha = 0$ for all $\alpha$ $P$-negative and not $P$-minimal, and $\omega^\gamma = a^\gamma \omega^{\beta^-(\gamma)}$ for any $P$-root which is not a $P'$-root, we find that still all of the expressions $\omega^{\mu\nu}$ become multiples of expressions $\omega^{\mu'\nu'}$. So on the subgroup of $P'$ on which all of these $\omega^\mu$ vanish, we have

$$da^\gamma = -a^\gamma(\gamma - \beta^-(\gamma)),$$

and so if $a^\gamma \ne 0$, as above, we can find a submanifold of $Z''$, say $Z'''$, on which $a^\gamma = 1$, and then $\gamma = \beta^-(\gamma) + b^\gamma_\mu \omega^\mu$. Hence $Z'''$ has dimension smaller than $G_{\check{X}}$, and is foliated by $G_Z$-orbits, with $G_Z$ acting freely, so $G_Z$ has smaller dimension than $G_{\check{X}}$. □

It is not clear whether $Z^{\text{free}}$ bears a holomorphic Cartan geometry modelled on $(\check{X}, G_{\check{X}})$. Nonetheless, the methods of [41] can be used to prove the following. Take any polynomial equation

$$0 = P(c_1, c_2, \dots, c_n)$$

satisfied by the Chern classes of the tangent bundle of $\check{X}$. Then the same equation is satisfied in Dolbeault cohomology on $Z^{\text{free}}$. We won't use this result, so we leave proof to the reader. It seems likely that every free compact complex submanifold is an associated cominuscule subvariety.

We prove theorem 5 on page 41:

*Proof.* Suppose that $V \subset TX$ is the maximal $G$-invariant holomorphic distribution not equal to $TX$, as above. The vector bundle inclusion morphism $TZ \to TX|_Z$ composes with the projection $TX \to TX/V$ to give a holomorphic vector bundle isomorphism $TX/V|_Z = TZ$, because $TZ$ is complementary to $V|_Z$. Since $X$ is homogeneous, $TX$ is spanned by global sections, so $TX/V$ is too, and so $TZ$ is too, so $Z$ is homogeneous. Since $Z \to X$ is free, it is an immersion. Pulling back a positive line bundle from $X$, $Z$ is a projective variety. Every homogeneous projective variety is such a product [5, 52]. □

## 11. Conclusion

The Hasse diagrams of flag varieties are mysterious. We understand the tip of the iceberg, almost literally, as we can predict the box: the component of the highest root. Each component of the Hasse diagram determines an irreducible invariant subbundle of the associated graded vector bundle $\operatorname{gr}_\bullet TX$ of the tangent bundle of $X$, and all irreducible invariant subbundles of $\operatorname{gr}_\bullet TX$ arise uniquely in this way. We can see how complicated the components get, but also see that there appears some attractive regularity in the pictures. We examine the noncompact root edges of the Hasse diagram of $G$ to see how those subbundles arise from the tangent bundle, and its invariant filtration. The invariant exterior differential systems (see [11, 10, 12, 27]) on flag varieties are not yet classified, and we don't know when they are involutive. Their integral manifolds are mysterious but natural submanifolds of flag varieties. It seems that invariant holomorphic Pfaffian systems on smooth complex projective varieties are usually entirely composed of Cauchy characteristics, and so, in some sense, trivial [17]. It might be that the flag varieties are very rare in having interesting exterior differential systems. It would be interesting to consider which minimal rational curves of an irreducible flag variety lie inside associated cominuscule subvarieties.

## Appendix A. Refining the filtrations and gradings

A.1. **The problem.** We want to clarify the filtrations of the tangent bundle of a flag variety by homogeneous vector bundles. There does not appear to be a convenient reference for the definitions of filtrations and gradings in suitable generality, i.e. over ordered groups, so we give them here.

A.2. **Graded vector spaces.** A *graded vector space* $V_\bullet$, graded by a set $A$, is a vector space $V$ which is the direct sum of linear subspaces $V_a$ for $a \in A$. (In practice, we will always grade by a finitely generated free abelian group $A \cong \mathbb{Z}^n$.) If $V_\bullet, W_\bullet$ are graded vector spaces, $(V \oplus W)_\bullet$ is the graded vector space

$$(V \oplus W)_a := (V_a) \oplus (W_a).$$

A *graded subspace* $V_\bullet \subseteq W_\bullet$ is a collection of linear subspaces

$$V_a \subseteq W_a.$$

The quotient $(W/V)_\bullet$ has grading

$$(W/V)_a := W_a/V_a.$$

Each component $V_a$ of a graded vector space is identified with the graded vector space $(V_a)_\bullet$ defined by

$$(V_a)_b = \begin{cases} V_a, & a = b, \\ 0, & a \neq b. \end{cases}$$

A.3. **Grading by groups.** If $A$ is a group, written additively, and $V_\bullet$ is a finite-dimensional graded vector space the *dual space* $V_\bullet^\vee$ is the graded vector space with

$$(V^\vee)_a := (V_{-a})^\vee.$$

If $A$ is a commutative group, or a commutative monoid, the tensor product $(V \otimes W)_\bullet$ of finite-dimensional graded vector spaces is the graded vector space

$$(V \otimes W)_a := \bigoplus_{b+c=a} V_b \otimes W_c.$$

A.4. **Graded Lie algebras.** A *graded Lie algebra* $\mathfrak{g}_\bullet$, graded by a commutative monoid $A$, is a Lie algebra $\mathfrak{g}$ with a grading of vector spaces $\mathfrak{g}_\bullet$ so that

$$[\mathfrak{g}_a, \mathfrak{g}_b] \subseteq \mathfrak{g}_{a+b}$$

for any $a, b \in A$ [39] p. 631, [40] chapter 5, pp. 92–93.

A.5. **Graded modules.** A $\mathfrak{g}_\bullet$*-module* $V_\bullet$ is a graded vector space, graded by a commutative monoid, which is a $\mathfrak{g}$-module so that

$$\mathfrak{g}_a V_b \subseteq V_{a+b}$$

for all $a, b \in A$. An *invariant set* is a subset $B \subseteq A$ so that, if $a \in A$ and $b \in B$ then $a + b \in B$. Let

$$V_a^B := \begin{cases} V_a, & \text{if } a \in B, \\ 0, & \text{otherwise.} \end{cases}$$

for $a \in A$. In particular, $\mathfrak{g}_\bullet^B \subseteq \mathfrak{g}_\bullet$ is a graded ideal. If a Lie group $G$ has Lie algebra $\mathfrak{g}$ equipped with a grading, a *graded $G$-module* is a $G$-module which, as a $\mathfrak{g}$-module, is equipped with a grading to be a $\mathfrak{g}_\bullet$-module. If $V_\bullet \subseteq W_\bullet$ is a graded $G$-module and a $G$-submodule, over the same commutative monoid, with

$$V_b = V \cap W_b$$

for all $b$ then $V_\bullet \subseteq W_\bullet$ is a $\mathfrak{g}_\bullet$*-submodule.* The quotient $(W/V)_\bullet$ is a graded $G$-module. If $V_\bullet$ is a finite dimensional graded $G$-module, and the grading is by a group, the dual space $V_\bullet^\vee$ is also a graded $G$-module. The sum and tensor product of graded $G$-modules, graded by the same group, are graded $G$-modules since

$$A(v \otimes w) = (Av) \otimes w + v \otimes (Aw)$$

for $A \in \mathfrak{g}_\bullet$, $v \in V_\bullet$, $w \in W_\bullet$.

A.6. **Filtered vector spaces.** A *filtered vector space* $V^\bullet$ filtered by a partially ordered set $A$ is a collection of linear subspaces $V^a \subseteq V$, for $a \in A$, of a vector space $V$, so that

- these subspaces span $V$ and
- if $b \geqslant a$ then $V^b \subseteq V^a$.

For any $a \in A$, define partially ordered sets $A_{\geqslant a}$ to be the set of elements $b \in A$ so that $b \geqslant a$, and $A_{>a}$ to be the set of elements $b \in A$ so that $b > a$. Take any vector space $V^\bullet$ filtered over $A$, and define a filtered vector space $V^\bullet_{\geqslant a}$ over $A_{\geqslant a}$ by

$$V^b_{\geqslant a} := V^b.$$

and define a filtered vector space $V^\bullet_{>a}$ over $A_{>a}$ by

$$V^b_{>a} := V^b.$$

If $V^\bullet, W^\bullet$ are filtered vector spaces over the same partially ordered set, let $(V \oplus W)^\bullet$ be the filtered vector space

$$(V \oplus W)^a := (V^a) \oplus (W^a)$$

Every linear subspace $V \subseteq W$ determines a *filtered subspace* $V^\bullet \subseteq W^\bullet$ by

$$V^a := V \cap W^a.$$

The quotient $(W/V)^\bullet$ has filtration

$$(W/V)^a := W^a/V^a.$$

A.7. **Changing partial ordering.** If $A \xrightarrow{\varphi} B$ is an order preserving map of partially ordered sets, and $V^\bullet, W^\bullet$ are filtered by $A, B$, a *morphism* $V^\bullet \xrightarrow{\Phi} W^\bullet$ over $\varphi$ is a collection of linear maps $V^a \xrightarrow{\Phi^a} W^b$ where $b = \varphi(a)$ so that $\Phi^a|_{V^{a'}} = \Phi^{a'}$ if $a \leqslant a'$ in $A$.

If $A \xrightarrow{\varphi} B$ is an order preserving map on partially ordered sets, and $V^\bullet$ is filtered by $A$, define $W^\bullet := \phi_* V^\bullet$ by setting $W^b$ to be the span in $V$ of all subsets $V^a$ for any $a \in A$ with $\varphi(a) \geqslant b$, but set $W^b := 0$ if there is no $a \in A$ with $\varphi(a) \geqslant b$. Define a morphism $V^\bullet \xrightarrow{\Phi^\bullet} W^\bullet$ over $\varphi$ by $\Phi^a(v) = v$ for $v \in V^a$ as a map $V^a \xrightarrow{\Phi^a} W^b$ with $b = \varphi(a)$. In particular, we can extend any filtered vector space defined on a subset $A \subseteq B$ of a partially ordered set to one defined on the whole set $B$.

If $A \xrightarrow{\varphi} B$ is an order preserving map on partially ordered sets, and $W^\bullet$ is filtered by $B$, define $V^\bullet := \phi^* W^\bullet$ by

$$V^a := W^b$$

if $b = \varphi(a)$, with morphism $V^\bullet \xrightarrow{\Phi^\bullet} W^\bullet$ defined by $\Phi^a(v) = v$ on each $V^a$; if $V$ is spanned by these $V^a$, this is a morphism of filtrations.

A.8. **Ordered groups.** An *ordered group* [9] chapter VI, [45] is an abelian group $A$ with a translation invariant partial order. Its *semipositive cone* is

$$A^+ := \{\, a \in A \mid a \geqslant 0 \,\}.$$

The order is determined by the semipositive cone: $a \geqslant b$ just when $a - b \in A^+$. (All of our examples will have $A \cong \mathbb{Z}^n$ a finitely generated abelian group, with distinguished generating set, and order by making the semipositive cone be the finite sums of the distinguished generators. For example, if $a, b, c$ are distinguished generators then $a + b$ is incomparable to $c$.)

Conversely, suppose that $A$ is a group and $A^+ \subseteq A$ is a subset. Then $A$ is an ordered group with $A^+$ as its semipositive cone if and only if all of the following hold: $A^+$ is closed under addition, contains 0, and, for any $a \in A^+$, $-a \in A^+$ if and only if $a = 0$.

Henceforth we will only filter or grade by an ordered group.

A.9. **Filtering by ordered groups.** The *tensor product* $(V \otimes W)^\bullet$ of finite-dimensional filtered vector spaces, filtered by an ordered group, is the filtered vector space

$$(V \otimes W)^a := \sum_{b+c \geqslant a} V^b \otimes W^c.$$

This sum is perhaps not a direct sum, since the summands could overlap. If $V^\bullet$ is a finite-dimensional filtered vector space the dual space $V^\vee$ is also a filtered vector space with

$$(V^\vee)^a := \bigcap_{a+b>0} (V^b)^\perp \subseteq V^\vee.$$

Given a filtered vector space $V^\bullet$, the *associated graded vector space* $\operatorname{gr}_\bullet V$ is

$$\operatorname{gr}_a V := V^a / \left\langle V^b \right\rangle_{b>a}.$$

Note that when defined

$$\operatorname{gr}_\bullet(V \oplus W) = (\operatorname{gr}_\bullet V) \oplus (\operatorname{gr}_\bullet W).$$

If there are gradings so that

$$V^a = \bigoplus_{b \geqslant a} V_b,$$

$$W^a = \bigoplus_{b \geqslant a} W_b,$$

then

$$\operatorname{gr}_\bullet(V \otimes W) = (\operatorname{gr}_\bullet V) \otimes (\operatorname{gr}_\bullet W).$$

Also, if

$$W^a \cap V^b = W^b$$

for all $b > a$ then

$$\operatorname{gr}_\bullet(V/W) = (\operatorname{gr}_\bullet V)/(\operatorname{gr}_\bullet W).$$

A.10. **Filtered Lie algebras.** A *filtered Lie algebra* $\mathfrak{g}^\bullet$ is a Lie algebra $\mathfrak{g}$ with a filtration of vector spaces $\mathfrak{g}^\bullet$ over an ordered group so that

$$\left[\mathfrak{g}^a, \mathfrak{g}^b\right] \subseteq \mathfrak{g}^{a+b}$$

for any $a, b \in A$. The associated graded vector space $\operatorname{gr}_\bullet \mathfrak{g}$ is a graded Lie algebra [63] chapter VIII p. 248: take any $x, y \in \operatorname{gr}_\bullet \mathfrak{g}$, say of grades $a, b$, and lift to elements $X, Y \in \mathfrak{g}$ with $X \in \mathfrak{g}^a$, $Y \in \mathfrak{g}^b$. Then $[X, Y] \in \mathfrak{g}^{a+b}$. If we change the choice of lifts $X, Y$, say to $X + X', Y + Y'$, the differences are zero in the quotient, i.e. $X' \in \sum_{a'>a} \mathfrak{g}^{a'}$ and $Y' \in \sum_{b'>b} \mathfrak{g}^{b'}$ for some $b' > b$, so

$$\left[X + X', Y + Y'\right] = [X, Y] + \varepsilon$$

differs from $[X, Y]$ by

$$\varepsilon = \left[X, Y'\right] + \left[X', Y\right] + \left[X', Y'\right] \in \langle \mathfrak{g}^c \rangle_{c>a+b} .$$

Hence we define

$$[x, y] := [X, Y] + \langle \mathfrak{g}^c \rangle_{c>a+b} ,$$

an antisymmetric bilinear map. The Jacobi identity in $\operatorname{gr}_\bullet \mathfrak{g}$ follows from the Jacobi identity in $\mathfrak{g}^\bullet$. If $\mathfrak{g}^\bullet$ is a filtered Lie algebra and $a \in A$ and $a \geqslant 0$ then $\mathfrak{g}^\bullet_{\geqslant a} \subseteq \mathfrak{g}^\bullet$ is a filtered subalgebra. Similarly, as defined above, $\mathfrak{g}^\bullet_+$ is the filtered Lie algebra with

$$\mathfrak{g}^a_+ := \begin{cases} \mathfrak{g}^a, & a > 0, \\ \mathfrak{g}_+, & \text{otherwise,} \end{cases}$$

where

$$\mathfrak{g}_+ = \left\langle \mathfrak{g}^b \right\rangle_{b>0} .$$

A.11. **Filtered modules.** A $\mathfrak{g}^\bullet$*-module* $V^\bullet$ is a filtered vector space which is a $\mathfrak{g}$-module so that

$$\mathfrak{g}^a V^b \subseteq V^{a+b}$$

for all $a, b \in A$. If a Lie group $G$ has Lie algebra $\mathfrak{g}$ equipped with a filtration, a *filtered $G$-module* is a $G$-module which is equipped with a $G$-invariant filtration as a $\mathfrak{g}^\bullet$-module. Henceforth we only consider filtered $G$-modules. For any $a \geqslant 0$ and any $b \in A$, if $V^\bullet$ is a filtered $G$-module then $V_{\geqslant b}$ is a filtered $\mathfrak{g}^\bullet_{\geqslant a}$ module. Similarly $V_{>b}$ is a filtered $\mathfrak{g}^\bullet_{\geqslant a}$ module.

If $\mathfrak{g}^\bullet = \mathfrak{g}^\bullet_{\geqslant 0}$, i.e. $\mathfrak{g}$ is *semipositively filtered*, then $\mathfrak{g}^\bullet_{\geqslant a}, \mathfrak{g}^\bullet_{>a}$ are invariant under brackets with $\mathfrak{g}^\bullet$. If in addition $G$ is connected then both are $G$-modules.

Suppose that $G$ is a connected complex Lie group and that $\mathfrak{g}$ is semipositively filtered. Then $\mathfrak{g}^\bullet_+$ acts on each filtered $G$-module $V^\bullet$ moving strictly upward its filtration degree. So $\mathfrak{g}^\bullet_+$ acts trivially on the associated graded $\operatorname{gr}_\bullet V$. Hence the action of $\mathfrak{g}$ reduces to an action of the quotient, i.e. of $\mathfrak{g}_0$. The Lie subalgebra $\mathfrak{g}^\bullet_+ \subseteq \mathfrak{g}^\bullet$ is an ideal. Hence the normal subgroup $G_+ \subseteq G$ with Lie algebra $\mathfrak{g}^\bullet_+$ acts trivially on $\operatorname{gr}_\bullet V$. (If $V^\bullet$ is an algebraic $G$-module, then the smallest Zariski closed normal subgroup of $G$ whose Lie algebra contains $\mathfrak{g}^\bullet_+$ acts trivially on $\operatorname{gr}_\bullet V$.) So we can quotient to an action of $G/G+$ on $\operatorname{gr}_\bullet V$.

If $V^\bullet \subseteq W^\bullet$ is a filtered $G$-module and a $G$-submodule with filtration

$$V^a = V \cap W^a$$

then $V^\bullet \subseteq W^\bullet$ is a *$G$-submodule.* The quotient $(W/V)^\bullet$ is a filtered $G$-module. If $V^\bullet$ is a filtered $G$-module, the dual space $V^\vee$ is also a filtered $G$-module. The sum and tensor product of filtered $G$-modules are filtered $G$-modules since

$$A(v \otimes w) = (Av) \otimes w + v \otimes (Aw)$$

for $A \in \mathfrak{g}^\bullet$, $v \in V^\bullet$, $w \in W^\bullet$.

A.12. **Augmentation.** An *augmented ordered group $A$* is an ordered group with a group epimorphism $A \to \mathbb{Z}$, denoted

$$a \in A \overset{\#}{\longmapsto} a^\# \in \mathbb{Z},$$

preserving $\leqslant$. In all of our examples below, $A$ is freely generated by a distinguished set of generators, and we set $a^\# = 1$ on these generators to augment.

A.13. **Underlying integer filtration.** To each filtered vector space $V^\bullet$ over an augmented ordered group, we can associate the *underlying integer filtration* ${}^\# V^\bullet$: the filtered vector space filtered by the integers, defined by taking ${}^\# V^j$ to be the span of the union of the $V^a$ for $a^\# \geqslant j$:

$${}^\# V^j := \langle V^a \rangle_{a^\# \geqslant j} .$$

Each filtered Lie algebra $\mathfrak{g}^\bullet$ has underlying integer filtration ${}^\#\mathfrak{g}^\bullet$ an integer filtered Lie algebra.

A.14. **Underlying integer grading.** Take a graded vector space $V_\bullet$. For any integer $j$, let

$${}^\# V_j := \bigoplus_{a^\# = j} V_a,$$

grading by the integers, the *underlying integer graded vector space.*

A.15. **Flag varieties and filtrations.** We will refine the integer filtration defined above. Take a flag variety $(X, G)$. Pick a point $x_0 \in X$. Let $P := G^{x_0}$. Pick a Cartan subgroup of $P$. Let $\Lambda$ be the root lattice, $\Lambda_c \subseteq \Lambda$ the subgroup generated by the $P$-compact roots and

$$\bar{\Lambda} := \Lambda / \Lambda_c.$$

For each $\lambda \in \Lambda$, let $\bar{\lambda}$ be its image

$$\bar{\lambda} := \lambda + \Lambda_c \in \bar{\Lambda}.$$

The images of the $P$-noncompact roots are generators of the free abelian group $\bar{\Lambda}$, so $\bar{\Lambda}$ is an ordered group: for any two roots $\lambda, \mu$, declare $\bar{\lambda} \leqslant \bar{\mu}$ if $\mu - \lambda$, expressed as an integer combination of simple roots

$$\mu - \lambda = \sum n_i \alpha_i$$

has $n_i \geqslant 0$ for every $P$-noncompact root $\alpha_i$. But recall that every root is always either a nonnegative linear combination of simple roots or a nonpositive one. So $\bar{\lambda} \leqslant \bar{\mu}$ just when $\mu - \lambda$ is a nonnegative integral linear combination of the images of the noncompact simple roots.

We refine the $P$-filtration: instead of grading by $P$-height, we grade by the value of $\bar{\lambda}$, as above. Augment by the sum of coefficients $\sum_i n_i$ of the $P$-noncompact roots. The induced integer filtration is precisely the $P$-height filtration.

Denote by $\mathfrak{g}_\alpha$ the $\alpha$-root space of $\mathfrak{g}$. Define a graded Lie algebra $\mathfrak{g}_\bullet$: for $[\bar{\alpha}] \neq 0$, set

$$\mathfrak{g}_{\bar{\alpha}} := \bigoplus_{\beta}^{\bar{\beta} = \bar{\alpha}} \mathfrak{g}_\beta$$

and then

$$\mathfrak{g}_0 := \mathfrak{h} \oplus \bigoplus_{\bar\beta=0} \mathfrak{g}_{\bar\beta}.$$

Define the induced filtered Lie algebra $\mathfrak{g}^\bullet$

$$\mathfrak{g}^{\bar\alpha} := \bigoplus_{\bar\beta \geqslant \bar\alpha} \mathfrak{g}_\beta = \bigoplus_{\bar\beta \geqslant \bar\alpha} \mathfrak{g}_{\bar\beta}.$$

Note that $\mathfrak{p} = \mathfrak{g}^{\bar 0}$, so $\mathfrak{p}$ is both filtered and graded, and $\mathfrak{g}_+ \subseteq \mathfrak{p}$ is a filtered and graded Lie subalgebra. The underlying integer filtrations and gradings are by $P$-height, so precisely those we discussed in §8.1 on page 24.

In particular, for any flag variety, this induces a filtered $P$-module structure on $\mathfrak{g}, \mathfrak{p}, \mathfrak{p}^\vee, \mathfrak{g}/\mathfrak{p}, \mathfrak{g}_+ = (\mathfrak{g}/\mathfrak{p})^\vee$. Each $P$-module is also a module over the Cartan subgroup, since the Cartan subgroup sits in $P$. We are mostly interested in the flag varieties $(X, G)$ where $G$ acts effectively on $X$, which implies that $G$ is in its adjoint form $G = \operatorname{Ad} G$. In particular, we are not interested in studying phenomena on a flag variety in which it differs from the associated effective flag variety. We will not assume that $G$ is in adjoint form, but rather we will consider only those $P$-modules which are trivial on the center of $G$, i.e. have weights only on the root lattice of $G$. Every such $P$-module is a sum of weight spaces

$$V = \bigoplus V_\lambda,$$

where all of the weights are in the root lattice, hence $V$ is graded into $G_0$-modules

$$V = \bigoplus V_{\bar\lambda},$$

and filtered by the associated $P$-invariant filtration. The associated graded is acted on trivially by $\mathfrak{g}_+$, as in §A.11 on page 50, so is a $P/G_+$-module. By Langlands decomposition, $P/G_+ = G_0$.

A.16. **Hasse diagrams of homogeneous vector bundles.** The *Hasse diagram* of a $P$-module $V$ has a vertex, a point of the plane, for each weight of $V$, and an edge marked $\ell$ between two vertices $\lambda, \mu$ just when $\lambda - \mu$ is a $P$-compact simple root labelled $\ell$ in the Dynkin diagram of $(X, G)$. The *grade* of a weight is the sum of its coefficients as a sum of simple roots. Larger grade weights lie on higher horizontal lines. When we take $X = G/B$ with $B$ the Borel subgroup, i.e. there are no $P$-compact roots, then this agrees with the Hasse diagrams of Plotkin et al. [49], except for their approach to the zero element of the weight lattice. The *Hasse diagram* of a homogeneous vector bundle $\mathbf{V} = G \times^P V$ is that of the $P$-module $V$.

*Example* 22. Take $(X, G)$ to be $G = B_4$, $X = G/P$, $P$ the Borel subgroup, the Hasse diagrams we have been drawing are of the parabolic subalgebra $\mathfrak{p}$:

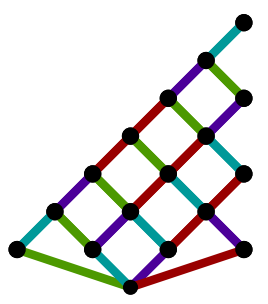

If instead we draw the Hasse diagram of $\mathfrak{g}$, we get the same picture drawn right side up and also upside down, since $\mathfrak{g} = \mathfrak{g}_- \oplus \mathfrak{h} \oplus \mathfrak{g}_+$ when $P$ is a Borel subgroup, i.e.

$\mathfrak{g}_0 = \mathfrak{h}$:

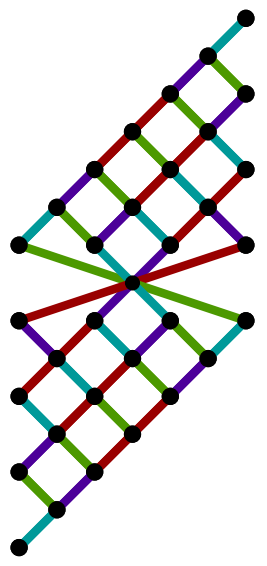

For simplicity, we will henceforth not draw the zero weight or the edges touching it. Due to the obvious reflection symmetry in the diagrams of $\mathfrak{g}$, the diagram of $\mathfrak{p}$ contains the same information as that of $\mathfrak{g}$. Therefore we have opted throughout, except in this section, to draw only the Hasse diagram of the parabolic subalgebra $\mathfrak{p}$.

Consider 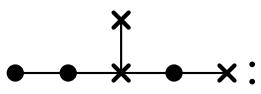: we pass from the Hasse diagram of $\mathfrak{p}$ to that of $\mathfrak{g}/\mathfrak{p}$ by turning it upside down and backward, and removing the $P$-compact roots:

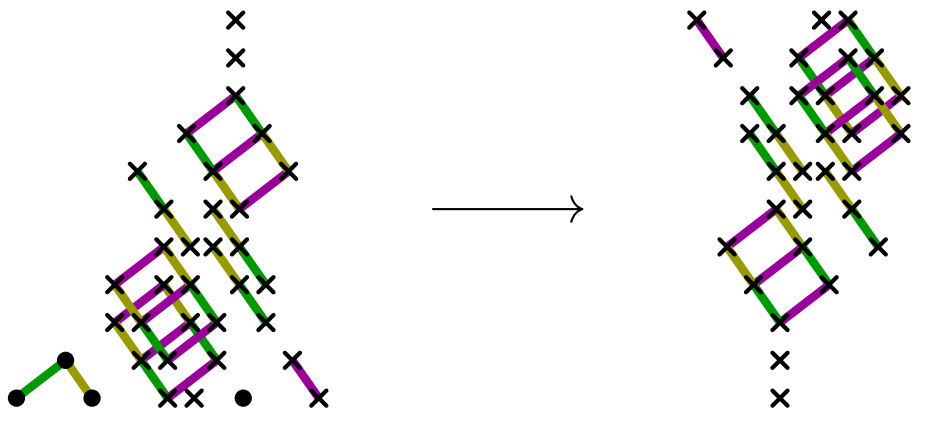

Again, since we can read the Hasse diagram of $\mathfrak{g}/\mathfrak{p}$ directly from the Hasse diagram of $\mathfrak{p}$, we will only draw the latter diagram henceforth.

A.17. **The associated graded as a module.** Suppose that $(X, G)$ is a flag variety. Pick a point $x_0 \in X$, let $P := G^{x_0}$. Consider $\mathfrak{g}^\bullet$ as a filtered $P$-module. The root spaces of the noncompact simple roots of $G$ move up the filtration, as

$$\mathfrak{g}^{\bar{\alpha}} V^{\bar{\beta}} \subseteq V^{\bar{\alpha}+\bar{\beta}}$$

Consequently, on the associated graded Lie algebra $\operatorname{gr}_\bullet \mathfrak{g}$, the root spaces of the noncompact simple roots of $G$ act trivially on each quotient

$$\operatorname{gr}_{\bar{\alpha}} \mathfrak{g} = \mathfrak{g}^{\bar{\alpha}} / \left\langle \mathfrak{g}^{\bar{\beta}} \right\rangle_{\bar{\beta} > \bar{\alpha}}.$$

Every element of $\mathfrak{g}_+$ has a positive grading, so moves up the grading of any element of $\mathfrak{g}^\bullet$. So $\mathfrak{g}_+$ also acts trivially on the associated graded $\operatorname{gr}_\bullet \mathfrak{g}$. Hence $G_+ = e^{\mathfrak{g}_+}$ acts trivially on the associated graded $\operatorname{gr}_\bullet \mathfrak{g}$. So $P$ acts, but the action descends to $G_0 \cong P/G_+$.

A.18. **The Hasse diagram and irreducibles.**

**Theorem 6.** *Suppose that $(X, G)$ is a flag variety. Pick a point $x_0 \in X$, let $P := G^{x_0}$. Let $G_0 \subseteq P$ be a reductive Levi factor containing the Cartan subgroup. Pick a Cartan subgroup of $G$ in $G_0$. The associated graded Lie algebra $\operatorname{gr}_\bullet \mathfrak{g}$ is the direct sum of $P$-modules*

$$\operatorname{gr}_\bullet \mathfrak{g} = \bigoplus_{\bar{\alpha}} \operatorname{gr}_{\bar{\alpha}} \mathfrak{g}.$$

*For $\alpha \neq 0$, these $P$-modules are irreducible. The subgroup $G_+ \subseteq P$ acts trivially on the associated graded, so each is a module of $G_0 \cong P/G_+$. Each of these $\operatorname{gr}_{\bar{\alpha}} \mathfrak{g}$, for* $\bar{\alpha} \neq 0$, *is the sum of the root spaces of a unique component of the Hasse diagram of* $\mathfrak{g}$. *Conversely every component of the Hasse diagram of* $\mathfrak{g}$ *occurs uniquely in this way. In particular, the box is the component of the highest root, and its associated* $P$*-module is* $\mathfrak{z}$, *the center of the unipotent radical* $\mathfrak{g}_+$ *of* $\mathfrak{p}$.

*Proof.* Our proof follows Wolf [62] p. 300 §8.13.3, who says that it is an unpublished result of Kostant; also see [38] p. 5 Proposition 2.4. Let $V := \operatorname{gr}_{\bar{\alpha}}\mathfrak{g}$. Clearly $V$ is a $P$-module, with trivial $G_+$-action, so a $G_0$-module, graded by weight:

$$V_\bullet = \bigoplus_\lambda V_\lambda.$$

As $G_0$-modules, $\operatorname{gr}_\bullet\mathfrak{g}$ is isomorphic to $\mathfrak{g}$, so each $V_\lambda$ is a $G_0$-submodule of $\mathfrak{g}$. Every root $\lambda$, as a weight, has a weight space of multiplicity one [19] p. 198, spanned by a root vector. Decompose $V$ into irreducible $G_0$-modules

$$V = \bigoplus_{i \in I} U_i.$$

Each $U_i$ has multiplicity one at its weights. Hence any distinct $U_i, U_j$ have disjoint weights. Each $U_i$ is a $G_0$-submodule of $\mathfrak{g}$, so its weight vectors are root vectors of $\mathfrak{g}$ and $\mathfrak{g}_0$ acts on each $U_i$ by Lie bracket.

Denote the $G_0$-highest weight of each $G_0$-irreducible $U_i$ by $\lambda_i$. Recall that these are $G_0$-integral and $G_0$-dominant, hence $G$-integral and $P$-compact dominant. To be precise: take a $P$-compact simple root $\alpha$. If $\langle \lambda_i, \alpha \rangle < 0$ we apply a root vector of the root $\alpha$ to get

$$\mathfrak{g}_{\lambda_i} \to \mathfrak{g}_{\lambda_i + \alpha}$$

to arrive at a higher weight in the same $G_0$-module [53] p. 29 Proposition 3. Since $\lambda_i$ is $G_0$-highest,

$$\langle \lambda_i, \alpha \rangle \geqslant 0$$

for all $P$-compact simple roots $\alpha$; indeed

$$2\frac{\lambda_i \cdot \alpha}{\alpha^2}$$

is a nonnegative integer [53] p. 60 theorem 3.

We won't need this, but we recall that $0 \leqslant \langle \lambda_i, \lambda_j \rangle$. In other words, since the simple roots are at least at right angles, the highest weights are at most at right angles [8] p. 168, just before Proposition 28. Indeed, the weights we are considering are positive roots, so nonnegative coefficient sums of simple roots. More generally, the fundamental weights of any reductive linear algebraic group are linear combinations of the simple roots, with coefficients the entries of the inverse Cartan matrix [51] p. 94, and these are positive [48] p. 295, [51] p. 95. We see the fundamental weights written in these linear combinations in [8] pp. 250–275 Planche I, part VI of each example. As above, the weights have nonnegative inner product with simple roots, and hence with one another.

Suppose that there are at least two irreducible $G_0$-submodules $U_i \neq U_j$. All weights of $V = \operatorname{gr}_{\bar{\alpha}}\mathfrak{g}$ have the same image $\bar{\alpha} \in \bar{\Lambda}$, so $\beta := \lambda_j - \lambda_i$ is in the lattice $\Lambda_c$ generated by the $P$-compact roots. We can suppose, after perhaps permuting indices $i$ and $j$, that $\beta$ is a sum of $P$-compact simple roots with nonnegative coefficients, at least one coefficient being positive. So

$$\langle \lambda_i, \lambda_j \rangle = \langle \lambda_i, \lambda_i + \beta \rangle = |\lambda_i|^2 + \langle \lambda_i, \beta \rangle > 0,$$

So $\beta$ is a $G$-root [53] p. 29 Proposition 3, but $P$-compact, hence is a $G_0$-root. Lie bracketing by a root vector of the root $-\beta$ isomorphically maps

$$\mathfrak{g}_{\lambda_j} \to \mathfrak{g}_{\lambda_i}$$

in the same $G_0$-module. But then one of these $\lambda_i, \lambda_j$ is not the highest $G_0$-weight of its associated $G_0$-module $U_i, U_j$, a contradiction. We conclude that there is only one $U_i$ factor in $V$, i.e. $V$ is irreducible.

Take a root $\alpha$ for which $V = \mathfrak{g}_{\bar{\alpha}}$. We can pick out $\alpha$ uniquely by asking that it is the lowest weight with given value of $\bar{\alpha}$. Being an irreducible $G_0$-module, the weights of $V$ have connected $G_0$-Hasse diagram, which is precisely the component of $\alpha$ in the Hasse diagram of $\mathfrak{g}$, by definition, since the simple roots of $G_0$ are precisely the $P$-compact simple roots of $G$.

On the other hand, take a component of the Hasse diagram of $\mathfrak{g}$. All of its roots $\alpha$ have the same value $\bar{\alpha}$, so all lie inside a single component of the $G_0$-Hasse diagram of $V = \mathfrak{g}_{\bar{\alpha}}$. □

**Corollary 1.** *Suppose that $(X, G)$ is a flag variety. Pick a point $x_0 \in X$, let $P := G^{x_0}$. Let $G_0 \subseteq P$ be a reductive Levi factor containing the Cartan subgroup. Pick a Cartan subgroup of $G$ in $G_0$.*

*Consider $V^{\bullet} = \mathfrak{p}^{\bullet}, (\mathfrak{g}/\mathfrak{p})^{\bullet}, (\mathfrak{p}^{\vee})^{\bullet}$ as filtered $P$-modules. Then the associated graded vector space $\operatorname{gr}_{\bullet} V$ is the direct sum of $P$-modules*

$$\operatorname{gr}_{\bullet} V = \bigoplus_{\bar{\alpha}} V_{\bar{\alpha}}.$$

*For $\alpha \neq 0$, these $P$-modules are irreducible: each of these $V_{\bar{\alpha}}$ is the sum of the root spaces of a unique component of the Hasse diagram of $V$, and hence of that of $\mathfrak{g}$. Conversely every component of the Hasse diagram of $V$ occurs uniquely in this way.*

*Proof.* Each associated graded $P$-module $\operatorname{gr}_{\bullet} V$ is isomorphic as a $G_0$-module to $V$, and each is then a $G_0$-submodule of $\mathfrak{g}$, so the same argument applies. □

**Corollary 2.** *Suppose that $(X, G)$ is a flag variety. Filter the tangent bundle $TX$ by its homogeneous vector subbundles. The associated graded vector bundle $\operatorname{gr}_{\bullet} TX$ is the direct sum of homogeneous vector subbundles, one for each component of the Hasse diagram of $\mathfrak{g}/\mathfrak{p}$ (as defined in §8.2 on page 24).*

If $K \subseteq G$ is a maximal compact subgroup, note that $T_{\bullet} X$ is isomorphic to $TX$ as a $K$-homogeneous vector bundle, so we are perhaps approaching an understanding of the tangent bundles of flag varieties.

*Example* 23. The $E_6$-flag variety with Dynkin diagram

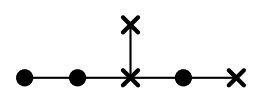

has Hasse diagram

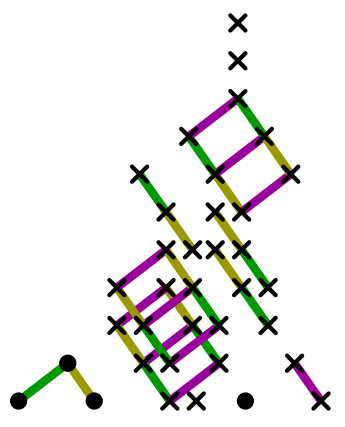

We see 10 components that have crosses on them. Hence the tangent bundle of this flag variety has associated graded vector bundle the direct sum of 10 irreducible homogeneous vector subbundles.

# Appendix B. The Dynkin diagrams of flag varieties and their associated cominuscule subvarieties

Table 6: The associated cominuscule subvariety of every irreducible flag variety. Some redundant entries are included where they might clarify. Compact simple roots are marked ×, noncompact •, and we mark by ∘ any node which could be either a × or •.

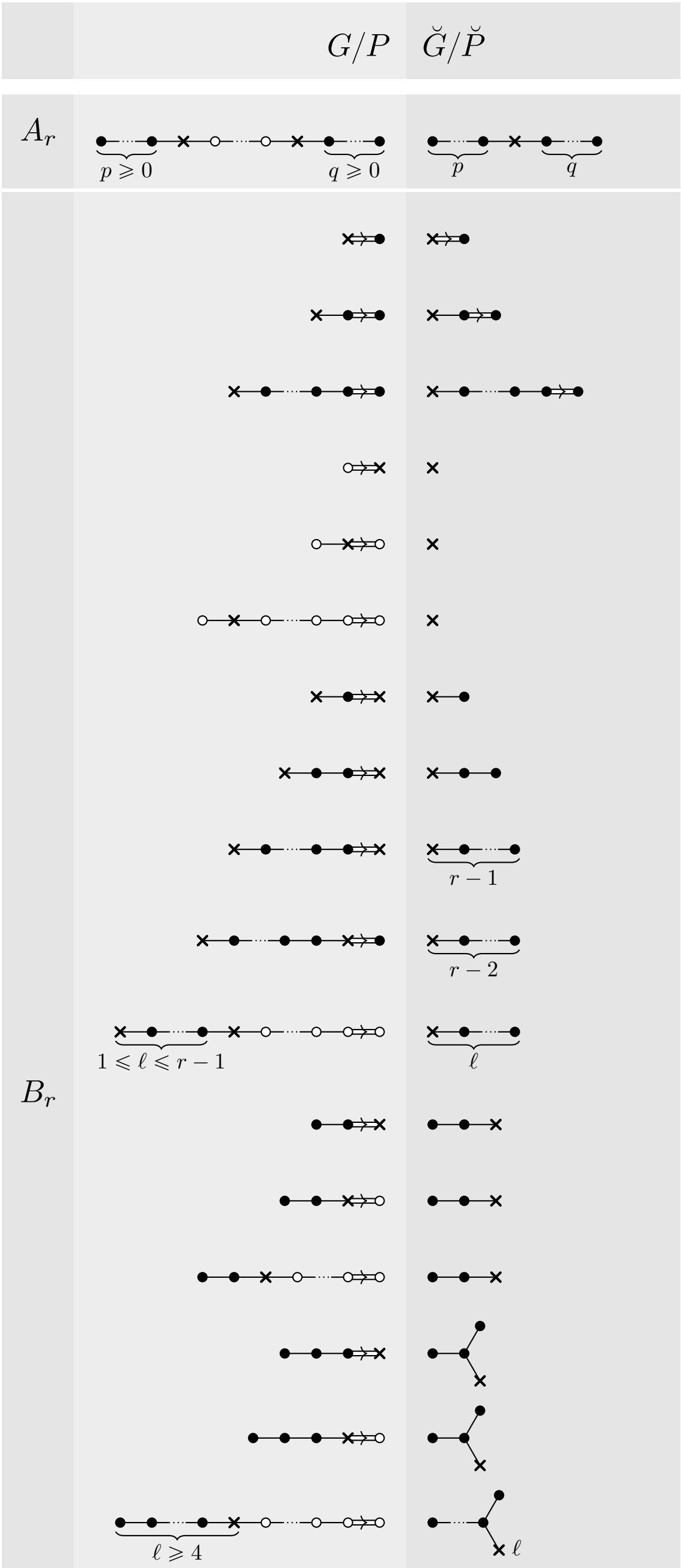

continued…

. . . continued

| | $G/P$ | $\breve{G}/\breve{P}$ |
|---|---|---|
| $C_r$ | 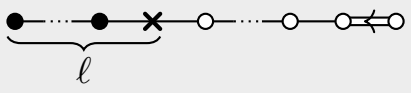  | 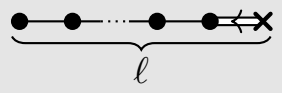  |
| $D_r$ | 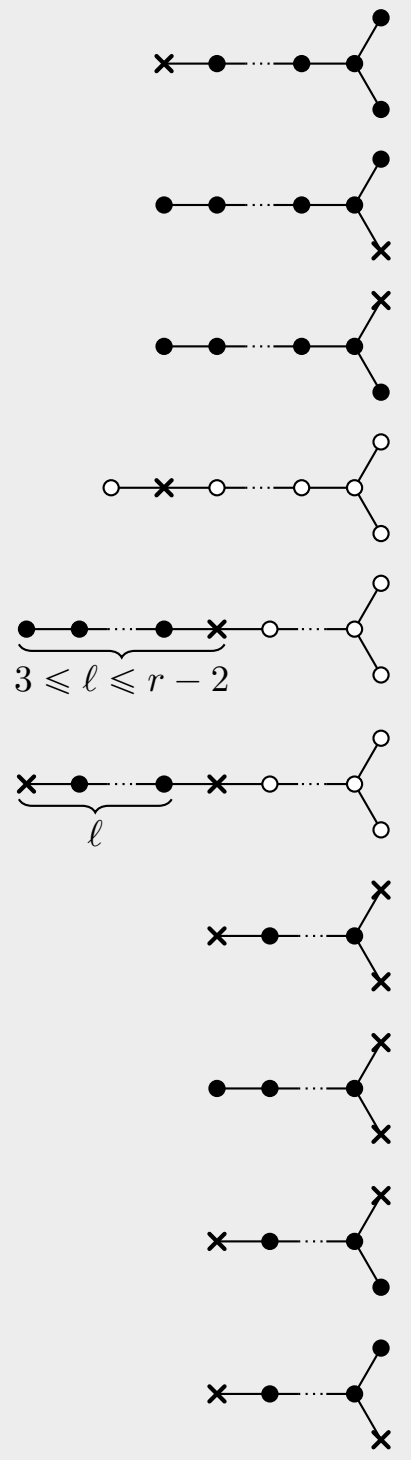  | 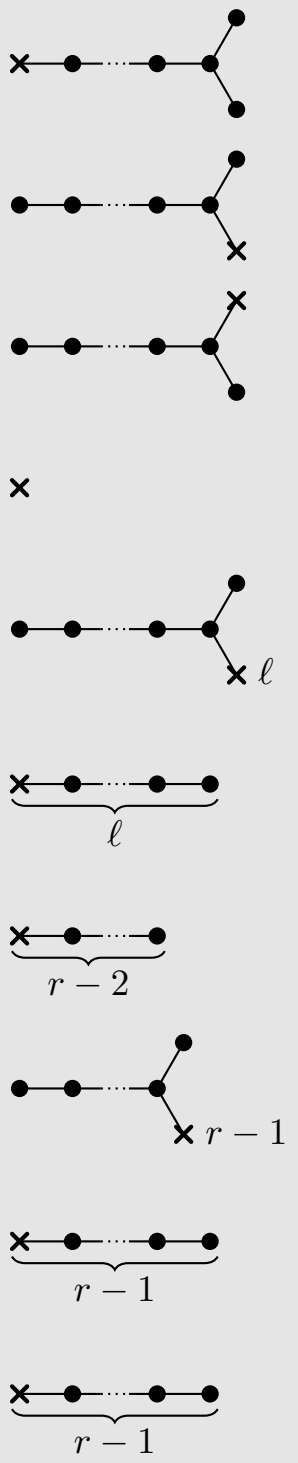  |
| $E_6$ | 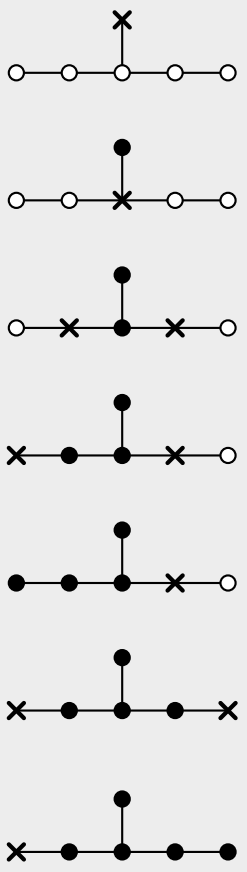 | 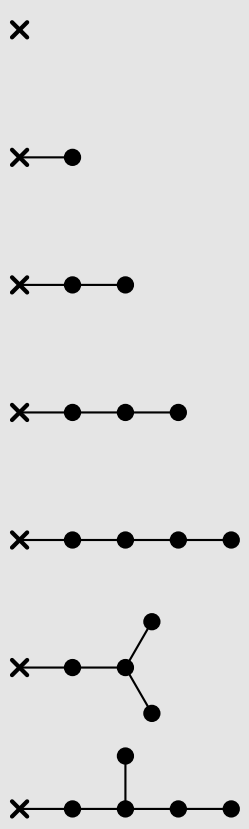 |

continued. . .

. . . continued

| | $G/P$ | $\breve{G}/\breve{P}$ |
|---|---|---|
| $E_7$ | 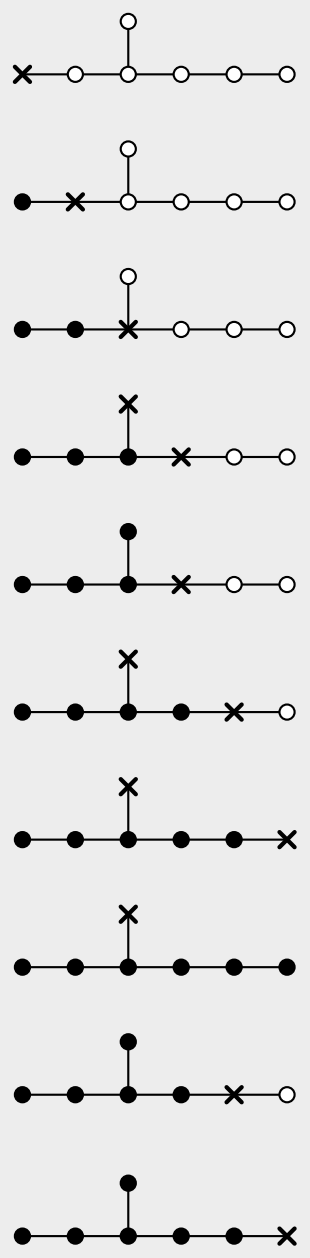 | 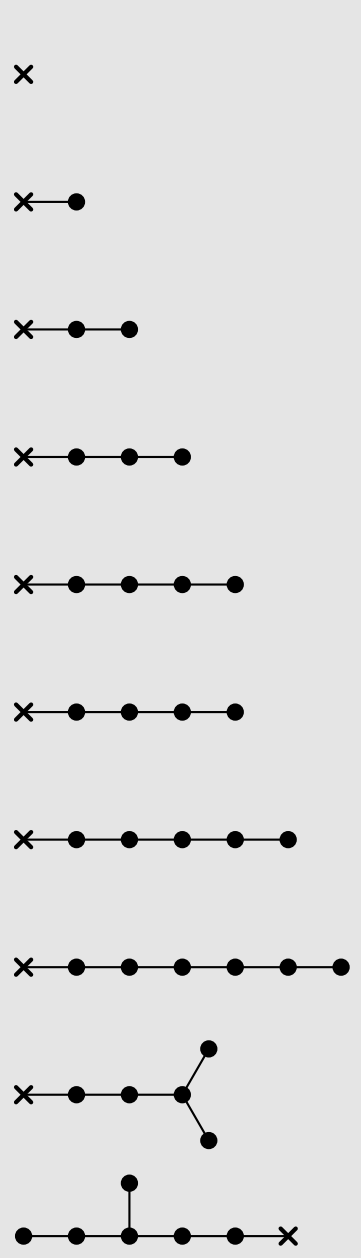 |
| $E_8$ | 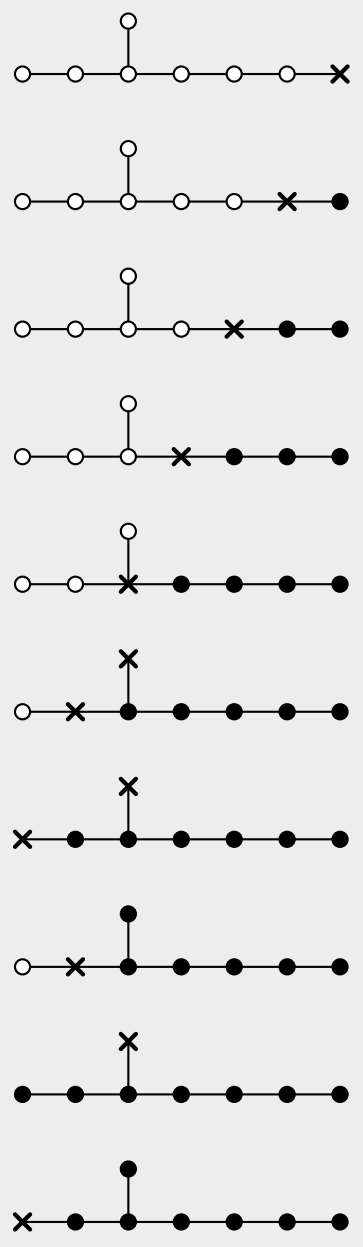 | 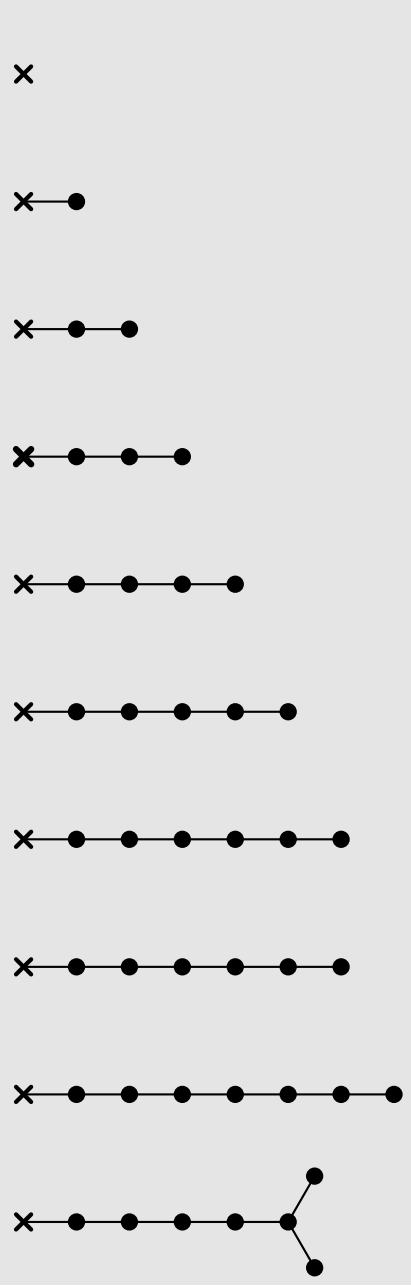 |

continued. . .

. . . continued

| | $G/P$ | $\check{G}/\check{P}$ |
|---|---|---|
| $F_4$ |  |  |
| |  |  |
| |  |  |
| |  |  |
| $G_2$ |  |  |
| |  |  |
| |  |  |

## Appendix C. The Hasse diagrams of the exceptional simple Lie groups

The lowest nodes sit in the Dynkin diagram, which we could draw in a plane perpendicular to the picture.

Table 7: Exceptional Hasse diagrams

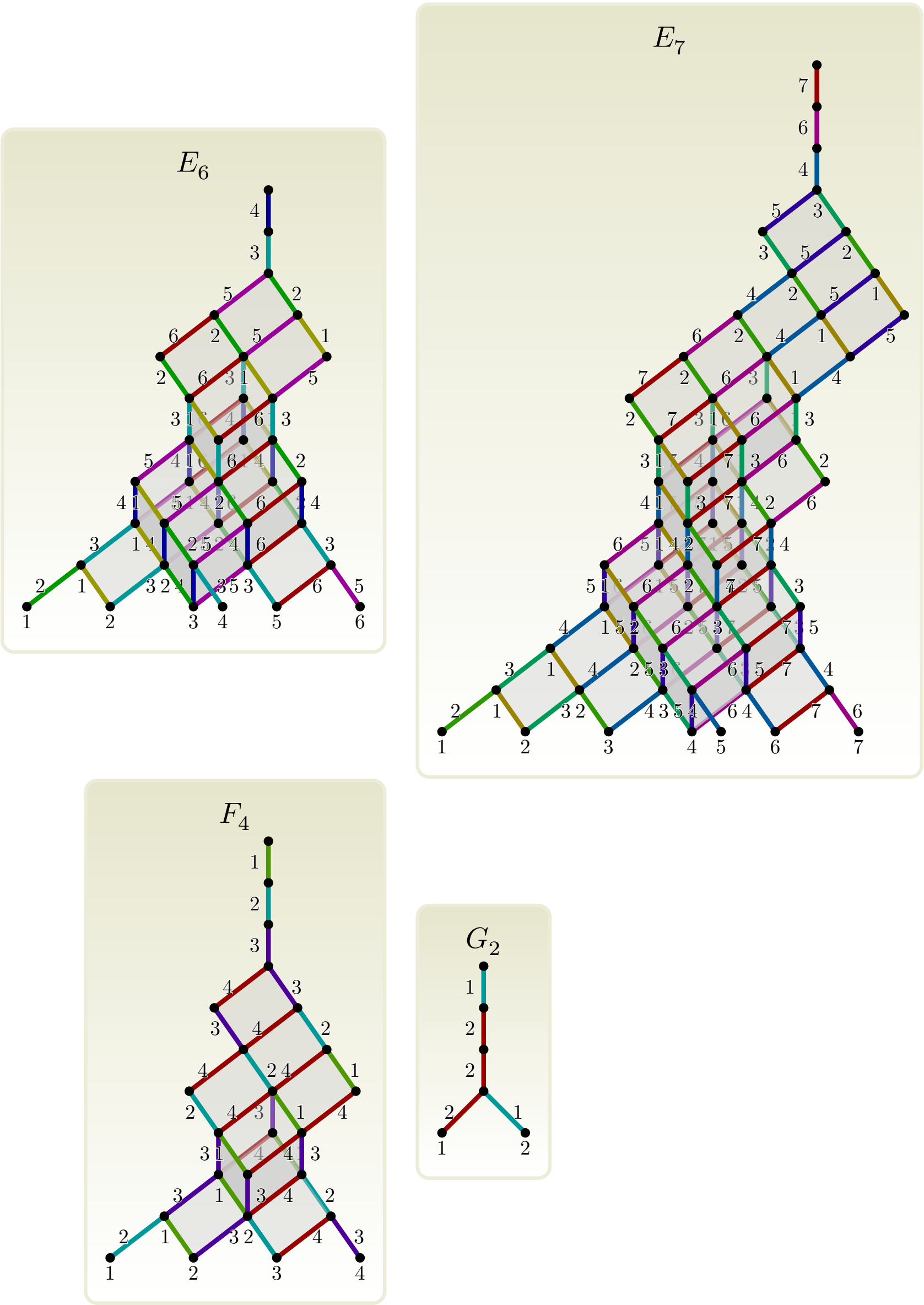

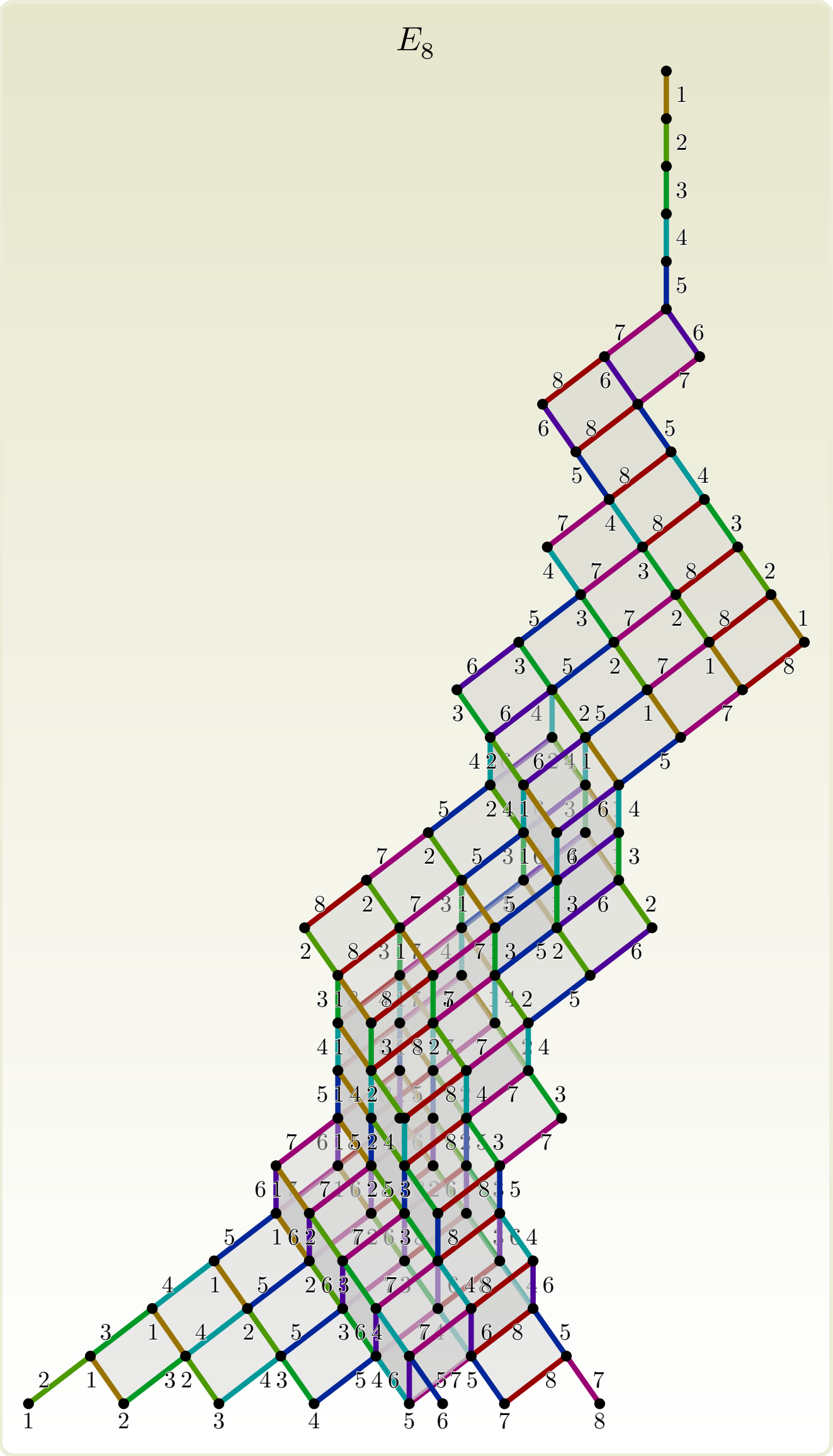
$E_8$

## Appendix D. The Hasse diagrams of the cominuscule varieties

Table 8: Hasse diagrams of the cominuscule varieties

| $G$ | $G/P$ | dim | description |
|---|---|---|---|
| $A_r$ |  | $k(r+1-k)$ | Grassmannian of $k$-planes in $\mathbb{C}^{r+1}$ |
| | 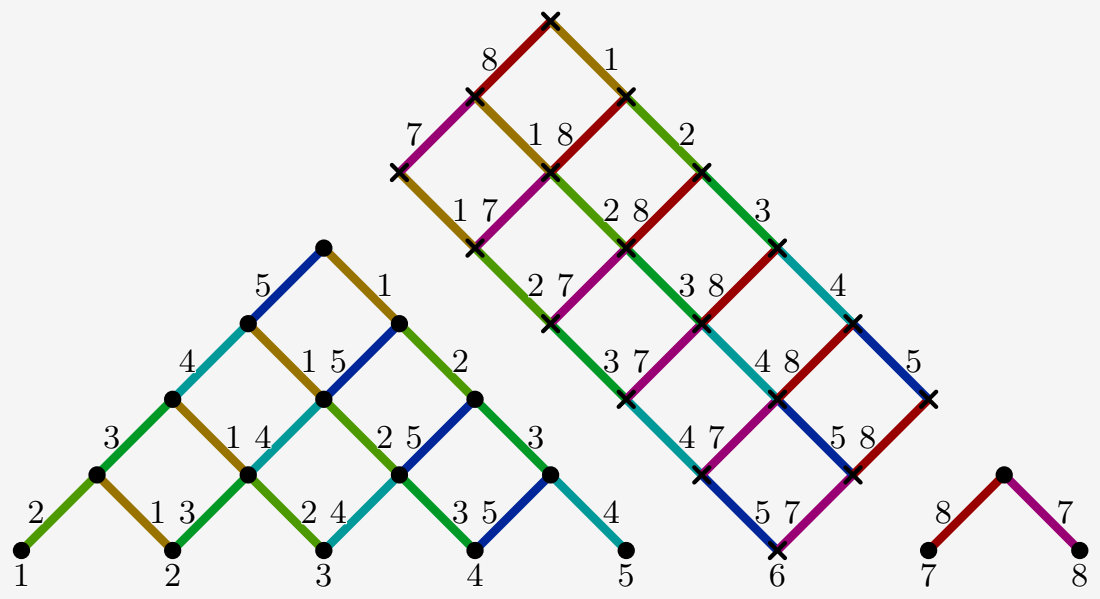 | | |
| $B_r$ |  | $2r-1$ | quadric hypersurface in $\mathbb{P}^{2r}$ |
| | 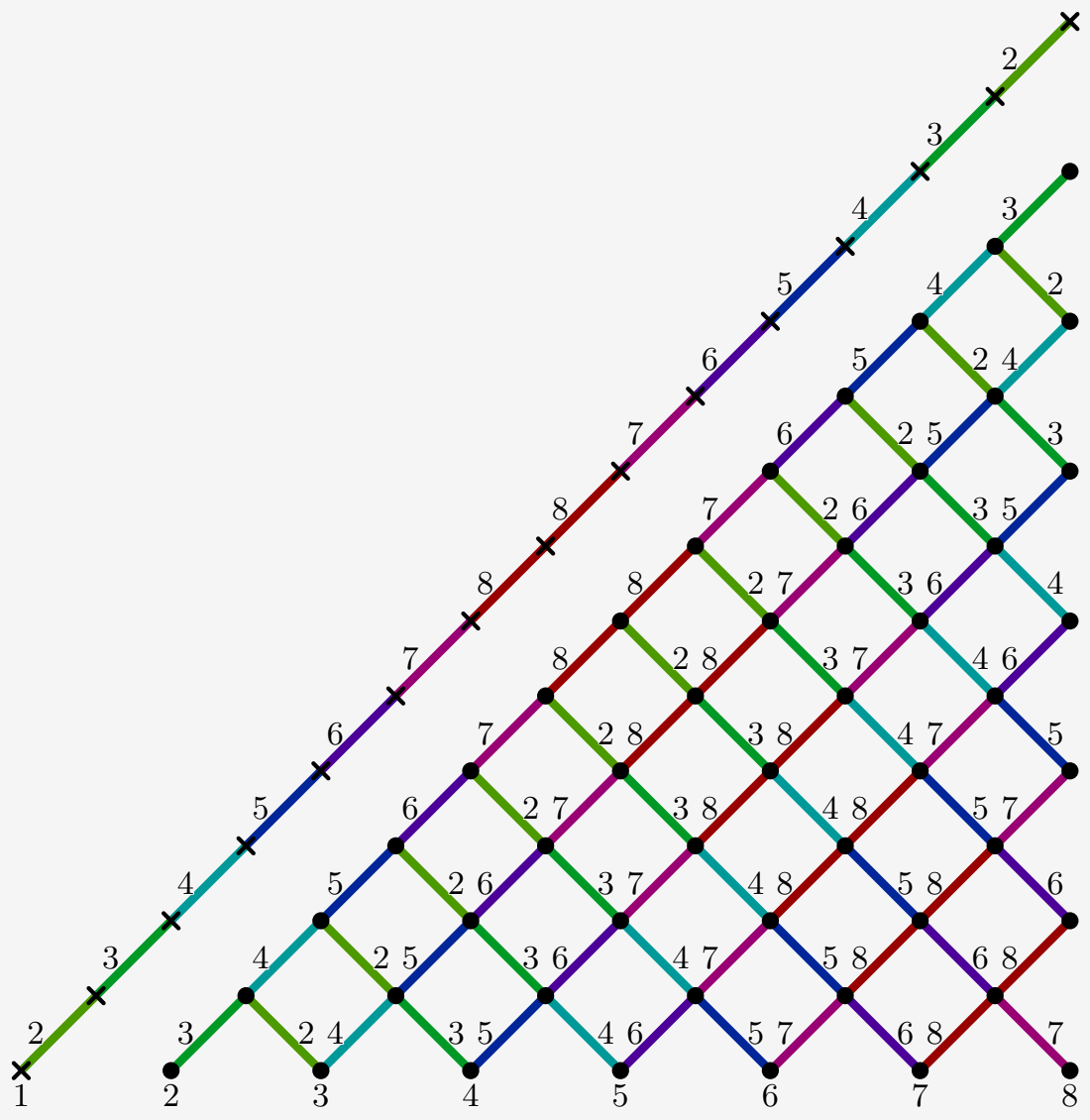 | | |

| $G$ | $G/P$ | dim | description |
|---|---|---|---|
| $C_r$ |  | $\frac{r(r+1)}{2}$ | Lagrangian $r$-planes in $\mathbb{C}^{2r}$ |
| | 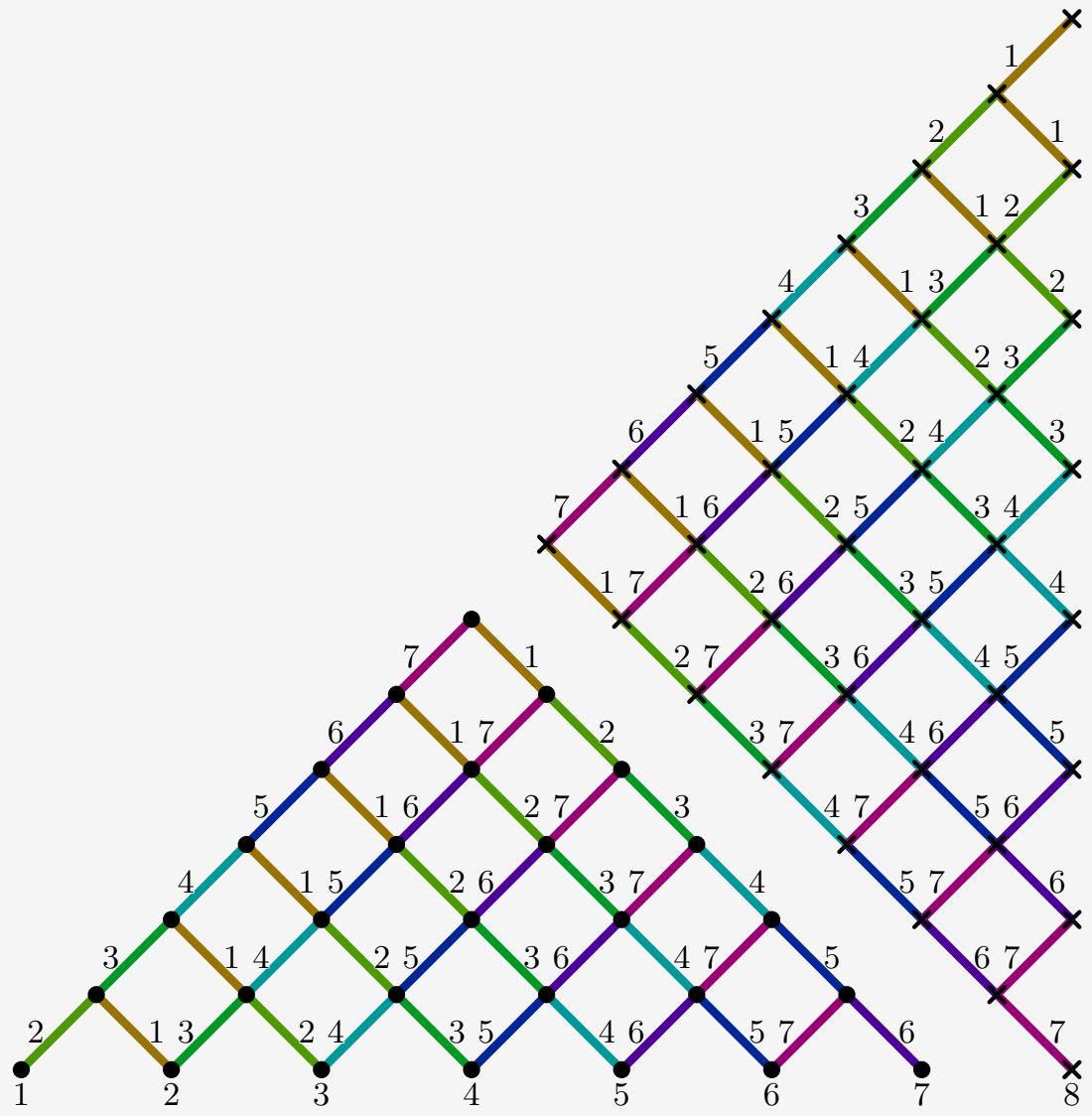 | | |
| $D_r$ | 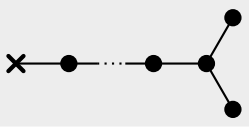 | $2r-2$ | quadric hypersurface in $\mathbb{P}^{2r-1}$ |
| | 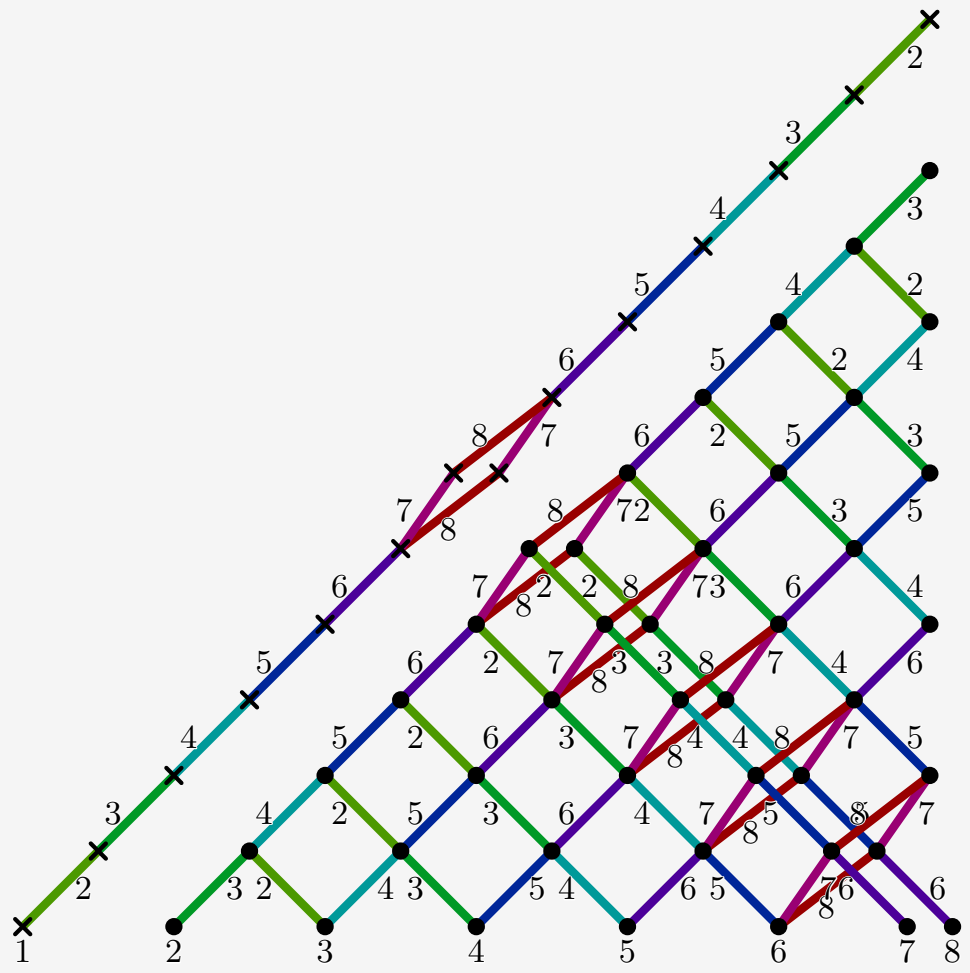 | | |

| $G$ | $G/P$ | dim | description |
|---|---|---|---|
| $D_r$ | 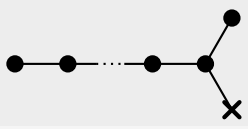 | $\frac{r(r-1)}{2}$ | null $r$-planes in $\mathbb{C}^{2r}$ |
| | 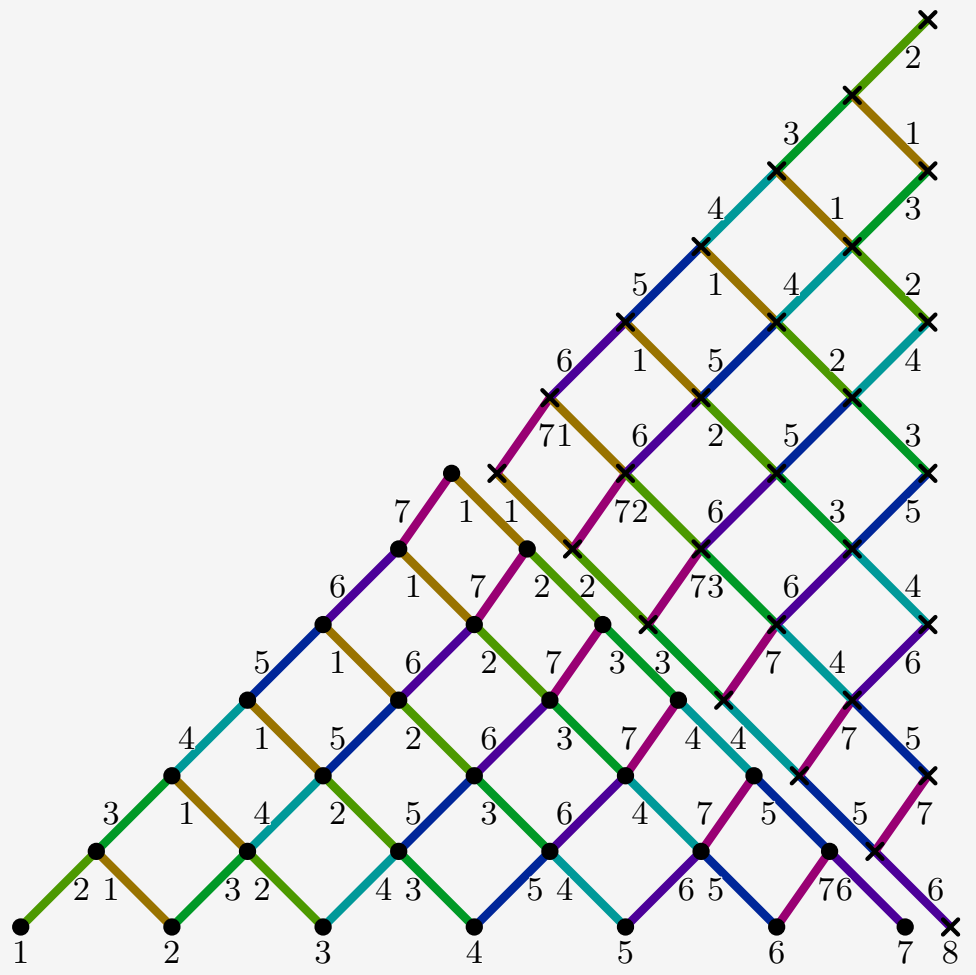 | | |
| $E_6$ | 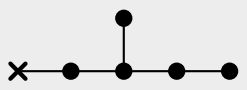 | 16 | complexified Cayley plane |
| | 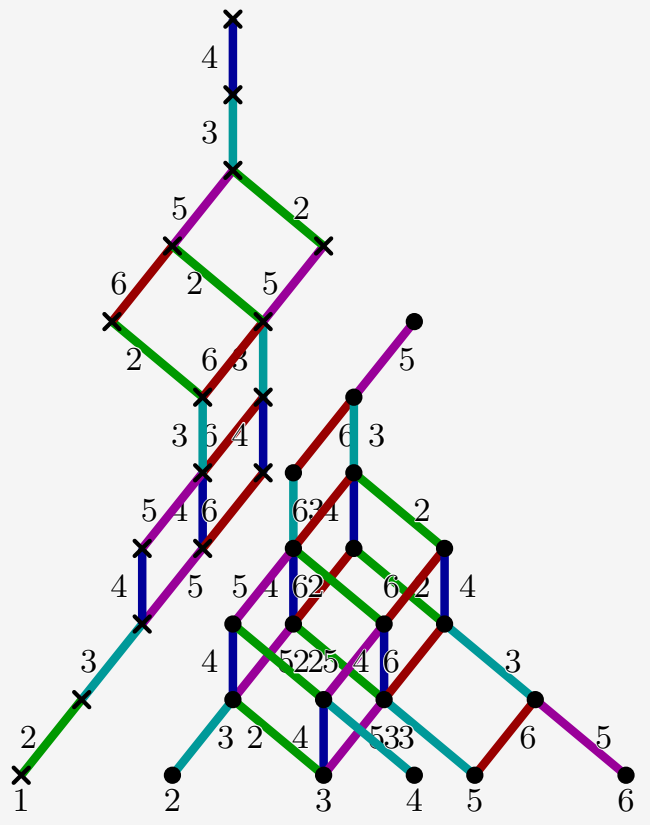 | | |

| $G$ | $G/P$ | dim | description |
|---|---|---|---|
| $E_7$ | 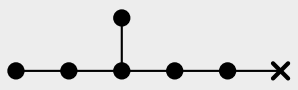 | 27 | null octave 3-planes in octave 6-space |

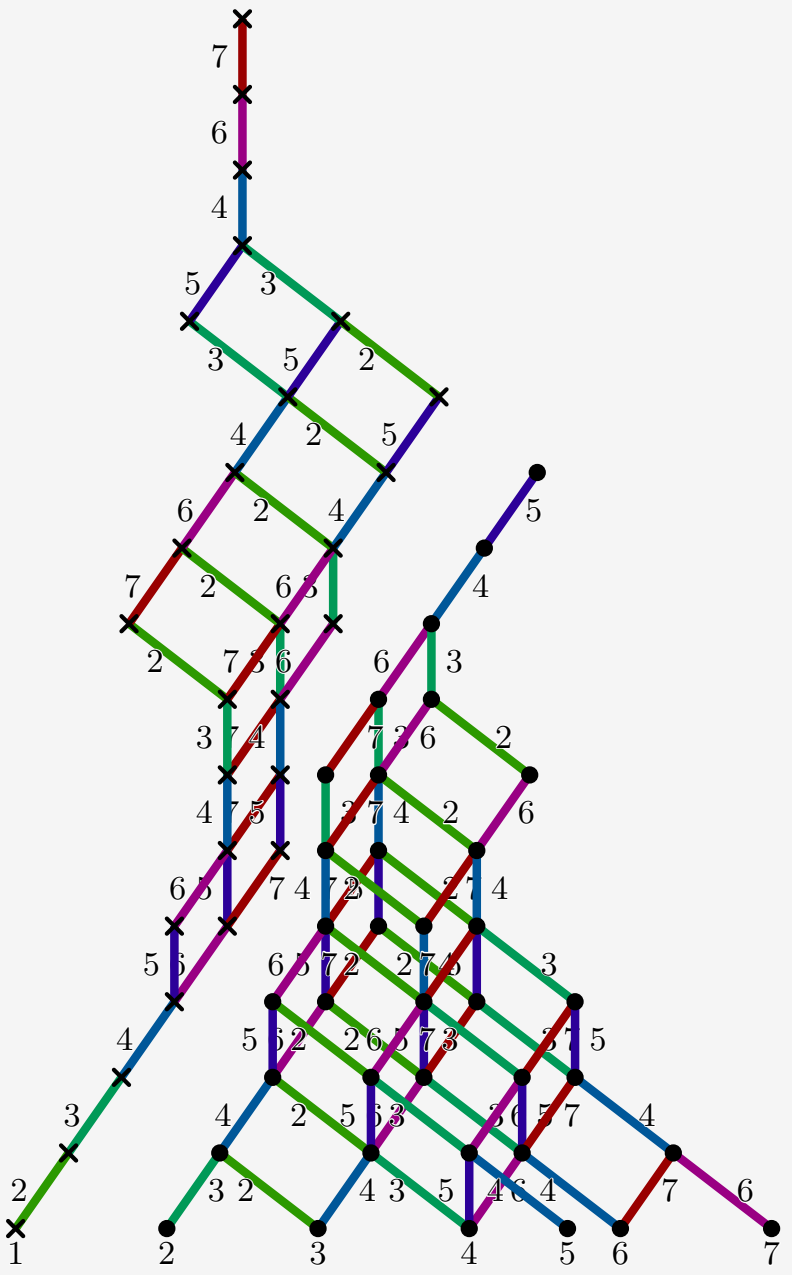

School of Mathematical Sciences, University College Cork, Cork, Ireland
*Email address*: b.mckay@ucc.ie